
\ifx\shlhetal\undefinedcontrolsequence\let\shlhetal\relax\fi

\input amstex
\expandafter\ifx\csname mathdefs.tex\endcsname\relax
  \expandafter\gdef\csname mathdefs.tex\endcsname{}
\else \message{Hey!  Apparently you were trying to
  \string\input{mathdefs.tex} twice.   This does not make sense.} 
\errmessage{Please edit your file (probably \jobname.tex) and remove
any duplicate ``\string\input'' lines}\endinput\fi




\catcode`\X=12\catcode`\@=11

\def\n@wcount{\alloc@0\count\countdef\insc@unt}
\def\n@wwrite{\alloc@7\write\chardef\sixt@@n}
\def\n@wread{\alloc@6\read\chardef\sixt@@n}
\def\r@s@t{\relax}\def\v@idline{\par}\def\@mputate#1/{#1}
\def\l@c@l#1X{\firstpart.#1}\def\gl@b@l#1X{#1}\def\t@d@l#1X{{}}

\def\crossrefs#1{\ifx\all#1\let\tr@ce=\all\else\def\tr@ce{#1,}\fi
   \n@wwrite\cit@tionsout\openout\cit@tionsout=\jobname.cit 
   \write\cit@tionsout{\tr@ce}\expandafter\setfl@gs\tr@ce,}
\def\setfl@gs#1,{\def\@{#1}\ifx\@\empty\let\next=\relax
   \else\let\next=\setfl@gs\expandafter\xdef
   \csname#1tr@cetrue\endcsname{}\fi\next}
\def\m@ketag#1#2{\expandafter\n@wcount\csname#2tagno\endcsname
     \csname#2tagno\endcsname=0\let\tail=\all\xdef\all{\tail#2,}
   \ifx#1\l@c@l\let\tail=\r@s@t\xdef\r@s@t{\csname#2tagno\endcsname=0\tail}\fi
   \expandafter\gdef\csname#2cite\endcsname##1{\expandafter
     \ifx\csname#2tag##1\endcsname\relax?\else\csname#2tag##1\endcsname\fi
     \expandafter\ifx\csname#2tr@cetrue\endcsname\relax\else
     \write\cit@tionsout{#2tag ##1 cited on page \folio.}\fi}
   \expandafter\gdef\csname#2page\endcsname##1{\expandafter
     \ifx\csname#2page##1\endcsname\relax?\else\csname#2page##1\endcsname\fi
     \expandafter\ifx\csname#2tr@cetrue\endcsname\relax\else
     \write\cit@tionsout{#2tag ##1 cited on page \folio.}\fi}
   \expandafter\gdef\csname#2tag\endcsname##1{\expandafter
      \ifx\csname#2check##1\endcsname\relax
      \expandafter\xdef\csname#2check##1\endcsname{}%
      \else\immediate\write16{Warning: #2tag ##1 used more than once.}\fi
      \multit@g{#1}{#2}##1/X%
      \write\t@gsout{#2tag ##1 assigned number \csname#2tag##1\endcsname\space
      on page \number\count0.}%
   \csname#2tag##1\endcsname}}

\def\multit@g#1#2#3/#4X{\def\t@mp{#4}\ifx\t@mp\empty%
      \global\advance\csname#2tagno\endcsname by 1 
      \expandafter\xdef\csname#2tag#3\endcsname
      {#1\number\csname#2tagno\endcsnameX}%
   \else\expandafter\ifx\csname#2last#3\endcsname\relax
      \expandafter\n@wcount\csname#2last#3\endcsname
      \global\advance\csname#2tagno\endcsname by 1 
      \expandafter\xdef\csname#2tag#3\endcsname
      {#1\number\csname#2tagno\endcsnameX}
      \write\t@gsout{#2tag #3 assigned number \csname#2tag#3\endcsname\space
      on page \number\count0.}\fi
   \global\advance\csname#2last#3\endcsname by 1
   \def\t@mp{\expandafter\xdef\csname#2tag#3/}%
   \expandafter\t@mp\@mputate#4\endcsname
   {\csname#2tag#3\endcsname\lastpart{\csname#2last#3\endcsname}}\fi}
\def\t@gs#1{\def\all{}\m@ketag#1e\m@ketag#1s\m@ketag\t@d@l p
\let\realscite\scite
\let\realstag\stag
   \m@ketag\gl@b@l r \n@wread\t@gsin
   \openin\t@gsin=\jobname.tgs \re@der \closein\t@gsin
   \n@wwrite\t@gsout\openout\t@gsout=\jobname.tgs }
\outer\def\localtags{\t@gs\l@c@l}
\outer\def\globaltags{\t@gs\gl@b@l}
\outer\def\newlocaltag#1{\m@ketag\l@c@l{#1}}
\outer\def\newglobaltag#1{\m@ketag\gl@b@l{#1}}

\newif\ifpr@ 
\def\m@kecs #1tag #2 assigned number #3 on page #4.%
   {\expandafter\gdef\csname#1tag#2\endcsname{#3}
   \expandafter\gdef\csname#1page#2\endcsname{#4}
   \ifpr@\expandafter\xdef\csname#1check#2\endcsname{}\fi}
\def\re@der{\ifeof\t@gsin\let\next=\relax\else
   \read\t@gsin to\t@gline\ifx\t@gline\v@idline\else
   \expandafter\m@kecs \t@gline\fi\let \next=\re@der\fi\next}
\def\pretags#1{\pr@true\pret@gs#1,,}
\def\pret@gs#1,{\def\@{#1}\ifx\@\empty\let\n@xtfile=\relax
   \else\let\n@xtfile=\pret@gs \openin\t@gsin=#1.tgs \message{#1} \re@der 
   \closein\t@gsin\fi \n@xtfile}

\newcount\sectno\sectno=0\newcount\subsectno\subsectno=0
\newif\ifultr@local \def\ultralocal{\ultr@localtrue}
\def\firstpart{\number\sectno}
\def\lastpart#1{\ifcase#1 \or a\or b\or c\or d\or e\or f\or g\or h\or 
   i\or k\or l\or m\or n\or o\or p\or q\or r\or s\or t\or u\or v\or w\or 
   x\or y\or z \fi}

\def\resetall{\global\advance\sectno by 1\subsectno=0
   \gdef\firstpart{\number\sectno}\r@s@t}
\def\resetsub{\global\advance\subsectno by 1
   \gdef\firstpart{\number\sectno.\number\subsectno}\r@s@t}
\def\newsection#1\par{\resetall\vskip0pt plus.3\vsize\penalty-250
   \vskip0pt plus-.3\vsize\bigskip\bigskip
   \message{#1}\leftline{\bf#1}\nobreak\bigskip}
\def\subsection#1\par{\ifultr@local\resetsub\fi
   \vskip0pt plus.2\vsize\penalty-250\vskip0pt plus-.2\vsize
   \bigskip\smallskip\message{#1}\leftline{\bf#1}\nobreak\medskip}


\newdimen\marginshift

\newdimen\margindelta
\newdimen\marginmax
\newdimen\marginmin

\def\margininit{       
\marginmax=3 true cm                  
				      
\margindelta=0.1 true cm              
\marginmin=0.1true cm                 
\marginshift=\marginmin
}    

\def\t@gsjj#1,{\def\@{#1}\ifx\@\empty\let\next=\relax\else\let\next=\t@gsjj
   \def\@@{p}\ifx\@\@@\else
   \expandafter\gdef\csname#1cite\endcsname##1{\citejj{##1}}
   \expandafter\gdef\csname#1page\endcsname##1{?}
   \expandafter\gdef\csname#1tag\endcsname##1{\tagjj{##1}}\fi\fi\next}
\newif\ifshowstuffinmargin
\showstuffinmarginfalse
\def\jjtags{\ifx\shlhetal\relax 
  \else
\ifx\shlhetal\undefinedcontrolseq
\else
\showstuffinmargintrue
\ifx\all\relax\else\expandafter\t@gsjj\all,\fi\fi \fi
}

\def\tagjj#1{\realstag{#1}\oldmginpar{\zeigen{#1}}}
\def\citejj#1{\rechnen{#1}\mginpar{\zeigen{#1}}}     

\def\rechnen#1{\expandafter\ifx\csname stag#1\endcsname\relax ??\else
                           \csname stag#1\endcsname\fi}

\newdimen\theight

\def\marginfont{\sevenrm}

\def\trymarginbox#1{\setbox0=\hbox{\marginfont\hskip\marginshift #1}%
		\global\marginshift\wd0 
		\global\advance\marginshift\margindelta}

\def \oldmginpar#1{%
\ifvmode\setbox0\hbox to \hsize{\hfill\rlap{\marginfont\quad#1}}%
\ht0 0cm
\dp0 0cm
\box0\vskip-\baselineskip
\else 
             \vadjust{\trymarginbox{#1}%
		\ifdim\marginshift>\marginmax \global\marginshift\marginmin
			\trymarginbox{#1}%
                \fi
             \theight=\ht0
             \advance\theight by \dp0    \advance\theight by \lineskip
             \kern -\theight \vbox to \theight{\rightline{\rlap{\box0}}%
\vss}}\fi}

\newdimen\upordown
\global\upordown=8pt
\font\tinyfont=cmtt8 
\def\mginpar#1{\smash{\hbox to 0cm{\kern-10pt\raise7pt\hbox{\tinyfont #1}\hss}}}
\def\mginpar#1{{\hbox to 0cm{\kern-10pt\raise\upordown\hbox{\tinyfont #1}\hss}}\global\upordown-\upordown}


\def\t@gsoff#1,{\def\@{#1}\ifx\@\empty\let\next=\relax\else\let\next=\t@gsoff
   \def\@@{p}\ifx\@\@@\else
   \expandafter\gdef\csname#1cite\endcsname##1{\zeigen{##1}}
   \expandafter\gdef\csname#1page\endcsname##1{?}
   \expandafter\gdef\csname#1tag\endcsname##1{\zeigen{##1}}\fi\fi\next}
\def\verbatimtags{\showstuffinmarginfalse
\ifx\all\relax\else\expandafter\t@gsoff\all,\fi}
\def\zeigen#1{\hbox{$\scriptstyle\langle$}#1\hbox{$\scriptstyle\rangle$}}


\def\(#1){\edef\dot@g{\ifmmode\ifinner(\hbox{\noexpand\etag{#1}})
   \else\noexpand\eqno(\hbox{\noexpand\etag{#1}})\fi
   \else(\noexpand\ecite{#1})\fi}\dot@g}

\newif\ifbr@ck
\def\eat#1{}
\def\[#1]{\br@cktrue[\br@cket#1'X]}
\def\br@cket#1'#2X{\def\temp{#2}\ifx\temp\empty\let\next\eat
   \else\let\next\br@cket\fi
   \ifbr@ck\br@ckfalse\br@ck@t#1,X\else\br@cktrue#1\fi\next#2X}
\def\br@ck@t#1,#2X{\def\temp{#2}\ifx\temp\empty\let\neext\eat
   \else\let\neext\br@ck@t\def\temp{,}\fi
   \def\teemp{#1}\ifx\teemp\empty\else\rcite{#1}\fi\temp\neext#2X}
\def\resetbr@cket{\gdef\[##1]{[\rtag{##1}]}}
\def\references{\resetbr@cket\newsection References\par}

\newtoks\symb@ls\newtoks\s@mb@ls\newtoks\p@gelist\n@wcount\ftn@mber
    \ftn@mber=1\newif\ifftn@mbers\ftn@mbersfalse\newif\ifbyp@ge\byp@gefalse
\def\defm@rk{\ifftn@mbers\n@mberm@rk\else\symb@lm@rk\fi}
\def\n@mberm@rk{\xdef\m@rk{{\the\ftn@mber}}%
    \global\advance\ftn@mber by 1 }
\def\rot@te#1{\let\temp=#1\global#1=\expandafter\r@t@te\the\temp,X}
\def\r@t@te#1,#2X{{#2#1}\xdef\m@rk{{#1}}}
\def\b@@st#1{{$^{#1}$}}\def\str@p#1{#1}
\def\symb@lm@rk{\ifbyp@ge\rot@te\p@gelist\ifnum\expandafter\str@p\m@rk=1 
    \s@mb@ls=\symb@ls\fi\write\f@nsout{\number\count0}\fi \rot@te\s@mb@ls}
\def\byp@ge{\byp@getrue\n@wwrite\f@nsin\openin\f@nsin=\jobname.fns 
    \n@wcount\currentp@ge\currentp@ge=0\p@gelist={0}
    \re@dfns\closein\f@nsin\rot@te\p@gelist
    \n@wread\f@nsout\openout\f@nsout=\jobname.fns }
\def\m@kelist#1X#2{{#1,#2}}
\def\re@dfns{\ifeof\f@nsin\let\next=\relax\else\read\f@nsin to \f@nline
    \ifx\f@nline\v@idline\else\let\t@mplist=\p@gelist
    \ifnum\currentp@ge=\f@nline
    \global\p@gelist=\expandafter\m@kelist\the\t@mplistX0
    \else\currentp@ge=\f@nline
    \global\p@gelist=\expandafter\m@kelist\the\t@mplistX1\fi\fi
    \let\next=\re@dfns\fi\next}
\def\symbols#1{\symb@ls={#1}\s@mb@ls=\symb@ls} 
\def\bigsymbol{\textstyle}
\symbols{\bigsymbol\ast,\dagger,\ddagger,\sharp,\flat,\natural,\star}
\def\ftnumbers{\ftn@mberstrue} \def\ftsymbols{\ftn@mbersfalse}
\def\paginal{\byp@ge} \def\resetftnumbers{\ftn@mber=1}
\def\ftnote#1{\defm@rk\expandafter\expandafter\expandafter\footnote
    \expandafter\b@@st\m@rk{#1}}

\long\def\jump#1\endjump{}
\def\ssum{\mathop{\lower .1em\hbox{$\textstyle\Sigma$}}\nolimits}

\def\qed{\nobreak\kern 1em \vrule height .5em width .5em depth 0em}
\def\newneq{\hbox{\rlap{\hbox to 1\wd9{\hss$=$\hss}}\raise .1em 
   \hbox to 1\wd9{\hss$\scriptscriptstyle/$\hss}}}
\def\subsetne{\setbox9 = \hbox{$\subset$}\mathrel{\hbox{\rlap
   {\lower .4em \newneq}\raise .13em \hbox{$\subset$}}}}
\def\supsetne{\setbox9 = \hbox{$\subset$}\mathrel{\hbox{\rlap
   {\lower .4em \newneq}\raise .13em \hbox{$\supset$}}}}

\def\vbar{\mathchoice{\vrule height6.3ptdepth-.5ptwidth.8pt\kern-.8pt}
   {\vrule height6.3ptdepth-.5ptwidth.8pt\kern-.8pt}
   {\vrule height4.1ptdepth-.35ptwidth.6pt\kern-.6pt}
   {\vrule height3.1ptdepth-.25ptwidth.5pt\kern-.5pt}}
\def\f@dge{\mathchoice{}{}{\mkern.5mu}{\mkern.8mu}}
\def\b@c#1#2{{\rm \mkern#2mu\vbar\mkern-#2mu#1}}
\def\b@b#1{{\rm I\mkern-3.5mu #1}}
\def\b@a#1#2{{\rm #1\mkern-#2mu\f@dge #1}}
\def\bb#1{{\count4=`#1 \advance\count4by-64 \ifcase\count4\or\b@a A{11.5}\or
   \b@b B\or\b@c C{5}\or\b@b D\or\b@b E\or\b@b F \or\b@c G{5}\or\b@b H\or
   \b@b I\or\b@c J{3}\or\b@b K\or\b@b L \or\b@b M\or\b@b N\or\b@c O{5} \or
   \b@b P\or\b@c Q{5}\or\b@b R\or\b@a S{8}\or\b@a T{10.5}\or\b@c U{5}\or
   \b@a V{12}\or\b@a W{16.5}\or\b@a X{11}\or\b@a Y{11.7}\or\b@a Z{7.5}\fi}}

\catcode`\X=11 \catcode`\@=12




\let\thischap\jobname

\def\partof#1{\csname returnthe#1part\endcsname}
\def\chapof#1{\csname returnthe#1chap\endcsname}

\def\setchapter#1,#2,#3;{%
  \expandafter\def\csname returnthe#1part\endcsname{#2}%
  \expandafter\def\csname returnthe#1chap\endcsname{#3}%
}

\setchapter 300a,A,II.A;
\setchapter 300b,A,II.B;
\setchapter 300c,A,II.C;
\setchapter 300d,A,II.D;
\setchapter 300e,A,II.E;
\setchapter 300f,A,II.F;
\setchapter 300g,A,II.G;
\setchapter  E53,B,N;
\setchapter  88r,B,I;
\setchapter  600,B,III;
\setchapter  705,B,IV;
\setchapter  734,B,V;

\def\cprefix#1{
\edef\theotherpart{\partof{#1}}\edef\theotherchap{\chapof{#1}}%
\ifx\theotherpart\thispart
   \ifx\theotherchap\thischap 
    \else 
     \theotherchap%
    \fi
   \else 
     \theotherchap\fi}

\def\sectioncite[#1]#2{%
     \cprefix{#2}#1}

\edef\thispart{\partof{\thischap}}
\edef\thischap{\chapof{\thischap}}

\def\lastpage of '#1' is #2.{\expandafter\def\csname lastpage#1\endcsname{#2}}


\def\spuriousreset{}


\expandafter\ifx\csname citeadd.tex\endcsname\relax
\expandafter\gdef\csname citeadd.tex\endcsname{}
\else \message{Hey!  Apparently you were trying to
\string\input{citeadd.tex} twice.   This does not make sense.} 
\errmessage{Please edit your file (probably \jobname.tex) and remove
any duplicate ``\string\input'' lines}\endinput\fi

\sectno=-1   
\localtags
\jjtags
\NoBlackBoxes
\define\mr{\medskip\roster}
\define\sn{\smallskip\noindent}
\define\mn{\medskip\noindent}
\define\bn{\bigskip\noindent}
\define\ub{\underbar}
\define\wilog{\text{without loss of generality}}
\define\rest{\restriction}
\define\ermn{\endroster\medskip\noindent}
\define\dbca{\dsize\bigcap}
\define\dbcu{\dsize\bigcup}
\define \nl{\newline}
\magnification=\magstep 1
\documentstyle{amsppt}

{    
\catcode`@11

\ifx\alicetwothousandloaded@\relax
  \endinput\else\global\let\alicetwothousandloaded@\relax\fi

\gdef\subjclass{\let\savedef@\subjclass
 \def\subjclass##1\endsubjclass{\let\subjclass\savedef@
   \toks@{\def\usualspace{{\rm\enspace}}\eightpoint}%
   \toks@@{##1\unskip.}%
   \edef\thesubjclass@{\the\toks@
     \frills@{{\noexpand\rm2000 {\noexpand\it Mathematics Subject
       Classification}.\noexpand\enspace}}%
     \the\toks@@}}%
  \nofrillscheck\subjclass}
} 


\expandafter\ifx\csname alice2jlem.tex\endcsname\relax
  \expandafter\xdef\csname alice2jlem.tex\endcsname{\the\catcode`@}
\else \message{Hey!  Apparently you were trying to
\string\input{alice2jlem.tex}  twice.   This does not make sense.}
\errmessage{Please edit your file (probably \jobname.tex) and remove
any duplicate ``\string\input'' lines}\endinput\fi

\expandafter\ifx\csname bib4plain.tex\endcsname\relax
  \expandafter\gdef\csname bib4plain.tex\endcsname{}
\else \message{Hey!  Apparently you were trying to \string\input
  bib4plain.tex twice.   This does not make sense.}
\errmessage{Please edit your file (probably \jobname.tex) and remove
any duplicate ``\string\input'' lines}\endinput\fi

\def\renewcommand{\newcommand}	       
\edef\cite{\the\catcode`@}%
\catcode`@ = 11
\let\@oldatcatcode = \cite
\chardef\@letter = 11
\chardef\@other = 12
%
%
%
%
\def\@innerdef#1#2{\edef#1{\expandafter\noexpand\csname #2\endcsname}}%
%
%
\@innerdef\@innernewcount{newcount}%
\@innerdef\@innernewdimen{newdimen}%
\@innerdef\@innernewif{newif}%
\@innerdef\@innernewwrite{newwrite}%
%
%
%
\def\@gobble#1{}%
%
%
%
\ifx\inputlineno\@undefined
   \let\@linenumber = \empty 
\else
   \def\@linenumber{\the\inputlineno:\space}%
\fi
%
%
%
\def\@futurenonspacelet#1{\def\cs{#1}%
   \afterassignment\@stepone\let\@nexttoken=
}%
\begingroup 
\def\\{\global\let\@stoken= }%
\\ 
\endgroup
\def\@stepone{\expandafter\futurelet\cs\@steptwo}%
\def\@steptwo{\expandafter\ifx\cs\@stoken\let\@@next=\@stepthree
   \else\let\@@next=\@nexttoken\fi \@@next}%
\def\@stepthree{\afterassignment\@stepone\let\@@next= }%
%
%
%
\def\@getoptionalarg#1{%
   \let\@optionaltemp = #1%
   \let\@optionalnext = \relax
   \@futurenonspacelet\@optionalnext\@bracketcheck
}%
%
%
\def\@bracketcheck{%
   \ifx [\@optionalnext
      \expandafter\@@getoptionalarg
   \else
      \let\@optionalarg = \empty
      \expandafter\@optionaltemp
   \fi
}%
\def\@@getoptionalarg[#1]{%
   \def\@optionalarg{#1}%
   \@optionaltemp
}%
%
%
%
\def\@nnil{\@nil}%
\def\@fornoop#1\@@#2#3{}%
\def\@for#1:=#2\do#3{%
   \edef\@fortmp{#2}%
   \ifx\@fortmp\empty \else
      \expandafter\@forloop#2,\@nil,\@nil\@@#1{#3}%
   \fi
}%
\def\@forloop#1,#2,#3\@@#4#5{\def#4{#1}\ifx #4\@nnil \else
       #5\def#4{#2}\ifx #4\@nnil \else#5\@iforloop #3\@@#4{#5}\fi\fi
}%
\def\@iforloop#1,#2\@@#3#4{\def#3{#1}\ifx #3\@nnil
       \let\@nextwhile=\@fornoop \else
      #4\relax\let\@nextwhile=\@iforloop\fi\@nextwhile#2\@@#3{#4}%
}%
%
%
%
\@innernewif\if@fileexists
\def\@testfileexistence{\@getoptionalarg\@finishtestfileexistence}%
\def\@finishtestfileexistence#1{%
   \begingroup
      \def\extension{#1}%
      \immediate\openin0 =
         \ifx\@optionalarg\empty\jobname\else\@optionalarg\fi
         \ifx\extension\empty \else .#1\fi
         \space
      \ifeof 0
         \global\@fileexistsfalse
      \else
         \global\@fileexiststrue
      \fi
      \immediate\closein0
   \endgroup
}%
%
%
%
%
\def\bibliographystyle#1{%
   \@readauxfile
   \@writeaux{\string\bibstyle{#1}}%
}%
\let\bibstyle = \@gobble
%
%
\let\bblfilebasename = \jobname
\def\bibliography#1{%
   \@readauxfile
   \@writeaux{\string\bibdata{#1}}%
   \@testfileexistence[\bblfilebasename]{bbl}%
   \if@fileexists
      \nobreak
      \@readbblfile
   \fi
}%
\let\bibdata = \@gobble
%
%
\def\nocite#1{%
   \@readauxfile
   \@writeaux{\string\citation{#1}}%
}%
\@innernewif\if@notfirstcitation
%
%
\def\cite{\@getoptionalarg\@cite}%
%
%
\def\@cite#1{%
   \let\@citenotetext = \@optionalarg
   \printcitestart
   \nocite{#1}%
   \@notfirstcitationfalse
   \@for \@citation :=#1\do
   {%
      \expandafter\@onecitation\@citation\@@
   }%
   \ifx\empty\@citenotetext\else
      \printcitenote{\@citenotetext}%
   \fi
   \printcitefinish
}%
\newif\ifweareinprivate
\weareinprivatetrue
\ifx\shlhetal\undefinedcontrolseq\weareinprivatefalse\fi
\ifx\shlhetal\relax\weareinprivatefalse\fi
\def\@onecitation#1\@@{%
   \if@notfirstcitation
      \printbetweencitations
   \fi
   \expandafter \ifx \csname\@citelabel{#1}\endcsname \relax
      \if@citewarning
         \message{\@linenumber Undefined citation `#1'.}%
      \fi
     \ifweareinprivate
      \expandafter\gdef\csname\@citelabel{#1}\endcsname{%
\strut 
\vadjust{\vskip-\dp\strutbox
\vbox to 0pt{\vss\parindent0cm \leftskip=\hsize 
\advance\leftskip3mm
\advance\hsize 4cm\strut\openup-4pt 
\rightskip 0cm plus 1cm minus 0.5cm ?  #1 ?\strut}}
         {\tt
            \escapechar = -1
            \nobreak\hskip0pt\pfeilsw
            \expandafter\string\csname#1\endcsname
             \pfeilso
            \nobreak\hskip0pt
         }%
      }%
     \else  
      \expandafter\gdef\csname\@citelabel{#1}\endcsname{%
            {\tt\expandafter\string\csname#1\endcsname}
      }%
     \fi  
   \fi
   \csname\@citelabel{#1}\endcsname
   \@notfirstcitationtrue
}%
%
%
\def\@citelabel#1{b@#1}%
%
%
\def\@citedef#1#2{\expandafter\gdef\csname\@citelabel{#1}\endcsname{#2}}%
%
%
%
\def\@readbblfile{%
   \ifx\@itemnum\@undefined
      \@innernewcount\@itemnum
   \fi
   \begingroup
      \def\begin##1##2{%
         \setbox0 = \hbox{\biblabelcontents{##2}}%
         \biblabelwidth = \wd0
      }%
      \def\end##1{}
      %
      %
      \@itemnum = 0
      \def\bibitem{\@getoptionalarg\@bibitem}%
      \def\@bibitem{%
         \ifx\@optionalarg\empty
            \expandafter\@numberedbibitem
         \else
            \expandafter\@alphabibitem
         \fi
      }%
      \def\@alphabibitem##1{%
         \expandafter \xdef\csname\@citelabel{##1}\endcsname {\@optionalarg}%
         \ifx\biblabelprecontents\@undefined
            \let\biblabelprecontents = \relax
         \fi
         \ifx\biblabelpostcontents\@undefined
            \let\biblabelpostcontents = \hss
         \fi
         \@finishbibitem{##1}%
      }%
      \def\@numberedbibitem##1{%
         \advance\@itemnum by 1
         \expandafter \xdef\csname\@citelabel{##1}\endcsname{\number\@itemnum}%
         \ifx\biblabelprecontents\@undefined
            \let\biblabelprecontents = \hss
         \fi
         \ifx\biblabelpostcontents\@undefined
            \let\biblabelpostcontents = \relax
         \fi
         \@finishbibitem{##1}%
      }%
      \def\@finishbibitem##1{%
         \biblabelprint{\csname\@citelabel{##1}\endcsname}%
         \@writeaux{\string\@citedef{##1}{\csname\@citelabel{##1}\endcsname}}%
         \ignorespaces
      }%
      %
      %
      \let\em = \bblem
      \let\newblock = \bblnewblock
      \let\sc = \bblsc
      \frenchspacing
      \clubpenalty = 4000 \widowpenalty = 4000
      \tolerance = 10000 \hfuzz = .5pt
      \everypar = {\hangindent = \biblabelwidth
                      \advance\hangindent by \biblabelextraspace}%
      \bblrm
      \parskip = 1.5ex plus .5ex minus .5ex
      \biblabelextraspace = .5em
      \bblhook
      \input \bblfilebasename.bbl
   \endgroup
}%
%
%
\@innernewdimen\biblabelwidth
\@innernewdimen\biblabelextraspace
%
%
%
\def\biblabelprint#1{%
   \noindent
   \hbox to \biblabelwidth{%
      \biblabelprecontents
      \biblabelcontents{#1}%
      \biblabelpostcontents
   }%
   \kern\biblabelextraspace
}%
%
%
%
\def\biblabelcontents#1{{\bblrm [#1]}}%
%
%
\def\bblrm{\rm}%
%
%
\def\bblem{\it}%
%
%
\def\bblsc{\ifx\@scfont\@undefined
              \font\@scfont = cmcsc10
           \fi
           \@scfont
}%
%
%
\def\bblnewblock{\hskip .11em plus .33em minus .07em }%
%
%
\let\bblhook = \empty
%
%
%
\def\printcitestart{[}
\def\printcitefinish{]}
\def\printbetweencitations{, }
\def\printcitenote#1{, #1}
%
%
%
\let\citation = \@gobble
%
%
%
\@innernewcount\@numparams
%
%
\def\newcommand#1{%
   \def\@commandname{#1}%
   \@getoptionalarg\@continuenewcommand
}%
%
%
\def\@continuenewcommand{%
   \@numparams = \ifx\@optionalarg\empty 0\else\@optionalarg \fi \relax
   \@newcommand
}%
%
%
\def\@newcommand#1{%
   \def\@startdef{\expandafter\edef\@commandname}%
   \ifnum\@numparams=0
      \let\@paramdef = \empty
   \else
      \ifnum\@numparams>9
         \errmessage{\the\@numparams\space is too many parameters}%
      \else
         \ifnum\@numparams<0
            \errmessage{\the\@numparams\space is too few parameters}%
         \else
            \edef\@paramdef{%
               \ifcase\@numparams
                  \empty  No arguments.
               \or ####1%
               \or ####1####2%
               \or ####1####2####3%
               \or ####1####2####3####4%
               \or ####1####2####3####4####5%
               \or ####1####2####3####4####5####6%
               \or ####1####2####3####4####5####6####7%
               \or ####1####2####3####4####5####6####7####8%
               \or ####1####2####3####4####5####6####7####8####9%
               \fi
            }%
         \fi
      \fi
   \fi
   \expandafter\@startdef\@paramdef{#1}%
}%
%
%
%
%
\def\@readauxfile{%
   \if@auxfiledone \else 
      \global\@auxfiledonetrue
      \@testfileexistence{aux}%
      \if@fileexists
         \begingroup
            \endlinechar = -1
            \catcode`@ = 11
            \input \jobname.aux
         \endgroup
      \else
         \message{\@undefinedmessage}%
         \global\@citewarningfalse
      \fi
      \immediate\openout\@auxfile = \jobname.aux
   \fi
}%
%
%
\newif\if@auxfiledone
\ifx\noauxfile\@undefined \else \@auxfiledonetrue\fi
%
%
%
%
\@innernewwrite\@auxfile
\def\@writeaux#1{\ifx\noauxfile\@undefined \write\@auxfile{#1}\fi}%
%
%
%
\ifx\@undefinedmessage\@undefined
   \def\@undefinedmessage{No .aux file; I won't give you warnings about
                          undefined citations.}%
\fi
%
%
\@innernewif\if@citewarning
\ifx\noauxfile\@undefined \@citewarningtrue\fi
%
%
%
\catcode`@ = \@oldatcatcode

\def\pfeilso{\leavevmode
            \vrule width 1pt height9pt depth 0pt\relax
           \vrule width 1pt height8.7pt depth 0pt\relax
           \vrule width 1pt height8.3pt depth 0pt\relax
           \vrule width 1pt height8.0pt depth 0pt\relax
           \vrule width 1pt height7.7pt depth 0pt\relax
            \vrule width 1pt height7.3pt depth 0pt\relax
            \vrule width 1pt height7.0pt depth 0pt\relax
            \vrule width 1pt height6.7pt depth 0pt\relax
            \vrule width 1pt height6.3pt depth 0pt\relax
            \vrule width 1pt height6.0pt depth 0pt\relax
            \vrule width 1pt height5.7pt depth 0pt\relax
            \vrule width 1pt height5.3pt depth 0pt\relax
            \vrule width 1pt height5.0pt depth 0pt\relax
            \vrule width 1pt height4.7pt depth 0pt\relax
            \vrule width 1pt height4.3pt depth 0pt\relax
            \vrule width 1pt height4.0pt depth 0pt\relax
            \vrule width 1pt height3.7pt depth 0pt\relax
            \vrule width 1pt height3.3pt depth 0pt\relax
            \vrule width 1pt height3.0pt depth 0pt\relax
            \vrule width 1pt height2.7pt depth 0pt\relax
            \vrule width 1pt height2.3pt depth 0pt\relax
            \vrule width 1pt height2.0pt depth 0pt\relax
            \vrule width 1pt height1.7pt depth 0pt\relax
            \vrule width 1pt height1.3pt depth 0pt\relax
            \vrule width 1pt height1.0pt depth 0pt\relax
            \vrule width 1pt height0.7pt depth 0pt\relax
            \vrule width 1pt height0.3pt depth 0pt\relax}

\def\pfeilsw{ \leavevmode 
            \vrule width 1pt height0.3pt depth 0pt\relax
            \vrule width 1pt height0.7pt depth 0pt\relax
            \vrule width 1pt height1.0pt depth 0pt\relax
            \vrule width 1pt height1.3pt depth 0pt\relax
            \vrule width 1pt height1.7pt depth 0pt\relax
            \vrule width 1pt height2.0pt depth 0pt\relax
            \vrule width 1pt height2.3pt depth 0pt\relax
            \vrule width 1pt height2.7pt depth 0pt\relax
            \vrule width 1pt height3.0pt depth 0pt\relax
            \vrule width 1pt height3.3pt depth 0pt\relax
            \vrule width 1pt height3.7pt depth 0pt\relax
            \vrule width 1pt height4.0pt depth 0pt\relax
            \vrule width 1pt height4.3pt depth 0pt\relax
            \vrule width 1pt height4.7pt depth 0pt\relax
            \vrule width 1pt height5.0pt depth 0pt\relax
            \vrule width 1pt height5.3pt depth 0pt\relax
            \vrule width 1pt height5.7pt depth 0pt\relax
            \vrule width 1pt height6.0pt depth 0pt\relax
            \vrule width 1pt height6.3pt depth 0pt\relax
            \vrule width 1pt height6.7pt depth 0pt\relax
            \vrule width 1pt height7.0pt depth 0pt\relax
            \vrule width 1pt height7.3pt depth 0pt\relax
            \vrule width 1pt height7.7pt depth 0pt\relax
            \vrule width 1pt height8.0pt depth 0pt\relax
            \vrule width 1pt height8.3pt depth 0pt\relax
            \vrule width 1pt height8.7pt depth 0pt\relax
            \vrule width 1pt height9pt depth 0pt\relax
      }


\def\widestnumber#1#2{}

\def\citewarning#1{\ifx\shlhetal\relax 
    \else
    \par{#1}\par
    \fi
}

\def\rm{\fam0 \tenrm}

\def\fakesubhead#1\endsubhead{\bigskip\noindent{\bf#1}\par}



%
%
%

%

\font\textrsfs=rsfs10
\font\scriptrsfs=rsfs7
\font\scriptscriptrsfs=rsfs5

\newfam\rsfsfam
\textfont\rsfsfam=\textrsfs
\scriptfont\rsfsfam=\scriptrsfs
\scriptscriptfont\rsfsfam=\scriptscriptrsfs

\edef\oldcatcodeofat{\the\catcode`\@}
\catcode`\@11

\def\Cal@@#1{\noaccents@ \fam \rsfsfam #1}

\catcode`\@\oldcatcodeofat


\expandafter\ifx \csname margininit\endcsname \relax\else\margininit\fi

\long\def\red#1\endred{}
\long\def\green#1\endgreen{}
\long\def\blue#1\endblue{}
\long\def\private#1\endprivate{}

\def\endred{ \unmatched endred! }
\def\endgreen{ \unmatched endgreen! }
\def\endblue{ \unmatched endblue! }
\def\endprivate{ \unmatched endprivate! }

\ifx\latexcolors\undefinedcs\def\latexcolors{}\fi

\def\emptycs{}
\def\evaluatelatexcolors{%
        \ifx\latexcolors\emptycs\else
        \expandafter\xxevaluate\latexcolors\xxfertig\evaluatelatexcolors\fi}
\def\xxevaluate#1,#2\xxfertig{\setupthiscolor{#1}%
        \def\latexcolors{#2}}


\font\smallfont=cmsl7
\def\rutgerscolor{\ifmmode\else\endgraf\fi\smallfont
\advance\leftskip0.5cm\relax}
\def\setupthiscolor#1{\edef\tmptmpcs{\noexpand\bgroup\noexpand\rutgerscolor
\noexpand\def\noexpand\currentcolor{#1}%
\noexpand}%
\expandafter\let\csname#1\endcsname\tmptmpcs
\def\tmptmpcs{\checkColorUnmatched{#1}\popthecolor}
\expandafter\let\csname end#1\endcsname\tmptmpcs}

\def\checkColorUnmatched#1{\def\expectcolor{#1}%
    \ifx\expectcolor\currentcolor   
    \else \edef\failhere{\noexpand\tryingToClose '\currentcolor' with end\expectcolor}\failhere\fi}

\def\currentcolor{???}

\def\popthecolor{\ifmmode\else\endgraf\fi\egroup}

\expandafter\def\csname#1\endcsname{}

\evaluatelatexcolors

 \let\outerhead\head
 \def\head{\innerhead}
 \let\innerhead\outerhead

 \let\outersubhead\subhead
 \def\subhead{\innersubhead}
 \let\innersubhead\outersubhead

 \let\outersubsubhead\subsubhead
 \def\subsubhead{\innersubsubhead}
 \let\innersubsubhead\outersubsubhead

 \let\outerproclaim\proclaim
 \def\proclaim{\innerproclaim}
 \let\innerproclaim\outerproclaim

 %
 %
 %
 %

\def\demo#1{\medskip\noindent{\it #1.\/}}
\def\enddemo{\smallskip}

\def\remark#1{\medskip\noindent{\it #1.\/}}
\def\endremark{\smallskip}

\pageheight{8.5truein}
\topmatter
\title{Advances in cardinal arithmetic} \endtitle
\author {Saharon Shelah \thanks {\null\newline 
I thank Alice Leonhardt for typing (and retyping) the manuscript so nicely and 
accurately. \null\newline
Partially supported by the BSF, Publ. 420} \endthanks} \endauthor 

\affil{Institute of Mathematics\\
 The Hebrew University\\
 Jerusalem, Israel
 \medskip
 Rutgers University\\
 Mathematics Department\\
 New Brunswick, NJ  USA} \endaffil
\endtopmatter
\document

\newpage











\head {Annotated Content} \endhead
 \spuriousreset
\bn
\S1 $\quad I[\lambda]$ is quite large 
\mr
\item "{{}}"  [If cf$\kappa = \kappa,\kappa^+ < \text{ cf}\lambda = \lambda$ then there is
a stationary subset $S$ of $\{\delta < \lambda:\text{cf}(\delta) = \kappa\}$  
in $I[\lambda]$.  Moreover, we can find $\bar C = \langle C_\delta :\delta  
\in  S\rangle $,  $C_\delta $ a club of $\lambda$, otp$(C_\delta) = 
\kappa$,  guessing clubs and for each $\alpha < \lambda$ we have:  
$\{C_\delta \cap \alpha:\alpha \in \text{ nacc } C_\delta\}$ has
cardinality $< \lambda$.]
\endroster
\bn
\S2 $\quad$ Measuring ${\Cal S}_{<\kappa }(\lambda)$ 
\mr
\item "{{}}"  [We prove that e.g. there is a stationary subset of  
${\Cal S}_{<\aleph_1}(\lambda)$ of cardinality cf$({\Cal
S}_{<\aleph_1}(\lambda),\subseteq)$.]
\endroster
\bn
\S3 $\quad$ Nice filters revisited 
\mr
\item "{{}}"  [We prove the existence of nice filters when instead being normal filters on  
$\omega_1$ they are normal filters with larger domains, which can increase 
during a play.  They can help us transfer situation on $\aleph _1$-complete 
filters to normal ones].
\endroster
\bn
\S4 $\quad$ Ranks 
\mr
\item "{{}}"  [We reconsider ranks and niceness of normal filters,
such that we can pass say from $pp_{\Gamma(\aleph_1)}(\mu)$ (where cf$\mu = \aleph_1)$ to  
pp$_{\text{normal}}(\mu)$.]
\endroster
\bn
\S5 $\quad$ More on ranks and higher objects
\bn
\S6 $\quad$ Hypotheses 
\mr
\item "{{}}"  [We consider some weakenings of G.C.H. and their 
consequences.  Most have not been proved independent of ZFC.]
\endroster
\newpage

\head{\S1  $I[\lambda]$ is Quite Large and Guessing Clubs} \endhead  \resetall \sectno=1
 \spuriousreset
\bigskip

On $I[\lambda]$ see 
\cite{Sh:108}, \cite{Sh:88a}, \cite[\S4]{Sh:351} (but this section is 
self-contained; see Definition \scite{1.1} and Claim \scite{1.2}
below).  
We shall prove that for regular 
$\kappa,\lambda$, such that $\kappa^+ < \lambda$,  there
is a stationary $S \subseteq \{\delta < \lambda:\text{cf}(\delta) 
= \kappa\}$ in $I[\lambda]$.
We then investigate ``guessing clubs" in (ZFC).
\bigskip

\definition{\stag{1.1} Definition}  For a regular uncountable cardinal
$\lambda,I[\lambda]$ is the family of $A \subseteq \lambda$ such that 
$\{\delta \in A:\delta = \text{ cf}(\delta)\}$ is not stationary and for some  
$\langle {\Cal P}_\alpha:\alpha < \lambda \rangle$ we have: 
\mr
\item "{$(a)$}"   ${\Cal P}_\alpha$ is a family of $< \lambda$ subsets of $\alpha$ 
\sn
\item "{$(b)$}"   for every limit 
$\alpha \in A$ of cofinality $< \alpha$ there is  
$x \subseteq \alpha$, otp$(x) < \alpha = \sup(x)$ 
such that $\zeta < \alpha \Rightarrow x \cap \zeta 
\in \{{\Cal P}_\gamma:\gamma < \alpha\}$.
\endroster
\enddefinition
\bigskip

\demo{\stag{1.1A} Observation}  In Definition \scite{1.1} we 
can weaken $(b)$ to:  
\mr
\item "{{}}"  for some club $E$ of $x$ for every limit $\alpha \in A
\cap E$ of cofinality $< \alpha \ldots$.
\endroster
\enddemo
\bigskip

\demo{Proof}  Just replace 
${\Cal P}_\alpha$ by $\{x \cap \alpha:x \in \cup \{{\Cal P}_\beta:\beta \le
\text{ Min}(E \backslash \alpha + 1)\}\}$. 
\enddemo
\bn
We know (see \cite{Sh:108}, \cite{Sh:88a} or below)
\proclaim{\stag{1.2} Claim}  Let $\lambda > \aleph_0$ be regular.
\nl
1)  $A \in I[\lambda]$ \ub{iff} 
(note: by $(c)$ below the set of inaccessibles in
$A$ is not stationary and) 
there is $\langle C_\alpha:\alpha < \lambda \rangle$ such that: 
\mr
\item "{$(a)$}"   $C_\alpha$ is a closed subset of $\alpha$ 
\sn
\item "{$(b)$}"  if $\alpha^* \in { \text{\rm nacc\/}}(C_\alpha)$ 
then $C_{\alpha^\ast} =  C_\alpha \cap \alpha$ 
(${\text{\rm nacc\/}}$ stands for ``non-accumulation") 
\sn
\item "{$(c)$}"  for some club $E$ of $\lambda$, for every 
$\delta \in A \cap E$, we have:  ${ \text{\rm cf\/}}(\delta) < \delta$
and $\delta = \sup(C_\delta)$, and ${\text{\rm cf\/}}(\delta ) = 
{ \text{\rm otp\/}}(C_\delta)$ 
\sn
\item "{$(d)$}" ${\text{\rm nacc\/}}(C_\alpha)$ is a set of successor ordinals.
\ermn
2)  $I[\lambda]$ is a normal ideal.
\endproclaim
\bigskip

\demo{Proof}  1) \ub{The ``if" part}: 

Assume $\langle C_\beta:\beta < \lambda \rangle$ satisfy $(a), (b), (c)$
with a club $E$ for $(c)$.  For each limit $\alpha  < \lambda$ choose a 
club $e_\alpha$ of order type cf$(\alpha)$.  
We define, for  $\alpha < \lambda$: 

$$
{\Cal P}_\alpha =: \{C_\beta:\beta \le \alpha\} \cup \{e_\beta:\beta  
\le \alpha\} \cup \{e_\gamma \cap \alpha:\gamma \le \text{ Min}(E\backslash (\alpha +1)\}.
$$
\mn
It is easy to check that $\langle {\Cal P}_\alpha:\alpha < \lambda \rangle$
exemplify  $``A \in I[\lambda]"$.
\medskip
\noindent
\ub{The ``only if" part}: 

Let $\bar{\Cal P} = 
\langle {\Cal P}_\alpha:\alpha < \lambda \rangle$ exemplify  $``A 
\in  I[\lambda]"$ (by Definition \scite{1.1}).  Without loss of generality 
\medskip
$(\ast)$ if  $C \in {\Cal P}_\alpha$,  and $\zeta \in C$  then  
$C \backslash \zeta \in {\Cal P}_\alpha$ 
and $C \cap \zeta \in {\Cal P}_\alpha$ 
\medskip
For each limit $\beta < \lambda$ let $e_\beta $ be a club of $\beta$
satisfying otp$(e_\beta) = \text{ cf}(\beta)$ and cf$(\beta) < \beta \Rightarrow \text{ cf}(\beta) < 
\min(e_\beta)$.  Let $\langle \gamma _i:i < \lambda \rangle $  be strictly 
increasing continuous, each  $\gamma _i$ a non-successor ordinal  $< \lambda $,
$\gamma _0 = 0$,  and $\gamma_{i+1} - \gamma_i \ge \aleph_0 +
|\dbcu_{\alpha \le \gamma_i} 
{\Cal P}_\alpha | + |\gamma _i|$  and  $\gamma _i \in  A 
\Rightarrow \text{ cf}(\gamma_i) < \gamma_i$. \nl
(Why?  Let $E'$ be a club of $\lambda$ such that $\gamma \in E \cap A
\Rightarrow \text{ cf}(\gamma) < \gamma$, and then choose $\gamma_i \in
E$ by induction on $i < \lambda$.) 

Let $F_i$ be a one to one function from $(\dbcu_{\alpha \le \gamma_i}
{\Cal P}_\alpha) \times \gamma_i$ into $\{\zeta  + 1:\gamma_i < \zeta + 1 <
\gamma_{i+1}\}$.  
Now we choose $C_\alpha \subseteq \alpha$ as follows. First, for
$\aleph=0$ let $C_\alpha = \emptyset$.  Second, assume $\alpha$ 
is a successor ordinal, let  $i(\alpha )$  be such that  
$\gamma _{i(\alpha )} < \alpha  < \gamma _{i(\alpha )+1}$.\  If  $\alpha  
\notin \text{ Rang}(F_{i(\alpha)})$, let $C_\alpha = \emptyset$.  
If $\alpha  = F_{i(\alpha)}(x,\beta)$ hence necessarily 
$x \in \dbcu_{\epsilon \le \gamma_{i(\alpha)}} {\Cal P}_\epsilon,
\beta < \gamma_{i(\alpha)})$ and $x,\beta$ are unique.  Let  
$C_\alpha $ be the closure (in the order topology) of $C^-_\alpha$,  
which is defined as:

$$
\bigl\{ F_j(x \cap \zeta,\beta):\text{ the sequence } (j,\zeta,\beta)
\text{ satisfies } (*)^{x,\beta}_{j,\zeta} \text{ below}\bigr\}
\text{ where}
$$
\mr
\item "{$\boxtimes^{x,\beta}_{j,\zeta}(i)$}"  $\zeta \in x$
\sn
\item "{${{}}(ii)$}"  $\text{otp}(x \cap \zeta) \in e_\beta$,
\sn
\item "{${{}}(iii)$}"  $j < i(\alpha) \text{ is minimal such that } x \cap
\zeta \in \dbcu_{\epsilon \le \gamma_j} {\Cal P}_\epsilon$
\sn
\item "{${{}}(iv)$}"  $\text{if } \xi \in x \cap \zeta, \text{ otp}(x \cap
\xi) \in e_\beta \text{ then}$ \nl  
$(\exists j(1) < j)[x \cap \xi \in 
\dbcu_{\epsilon \le \gamma_{j(1)}} {\Cal P}_\epsilon]$
\sn
\item "{${{}}(v)$}"  $\beta < \text{ Min}(x)$.
\ermn
Third, for $\alpha < \lambda$ limit, choose $C_\alpha$: if possible,  
nacc$(C_\alpha)$ is a set of successor ordinals, $C_\alpha$ is a club of  
$\alpha,[\beta \in \text{ nacc}(C_\alpha) \Rightarrow C_\beta  = \beta  \cap  
C_\alpha ]$;  if this is impossible, let  $C_\delta = \emptyset$.  
Lastly, let $C_0 = \emptyset$ and let $E =: \{\gamma_i:i$ is a limit
ordinal $< \lambda\}$.
\nl
Now we can check the condition in \scite{1.2}(1).
\mn
Note that for  $\alpha$ successor $C^-_\alpha = \text{ nacc}(C_\alpha)$.
\sn
\ub{Clause (a)}:  $C_\alpha $ a closed subset of  $\alpha$. 

If $\alpha = 0$ trivial as $C_\alpha  = \emptyset$ and if $\alpha$ is a limit 
ordinal, this is immediate by the definition.  So let  $\alpha $  be a 
successor ordinal, hence, by the choice of $\langle \gamma_i:i <
\lambda \rangle$ as an increasing continuous sequence of nonsuccessor
ordinals with $\gamma_0=0$, clearly $i(\alpha)$ is well defined,  
$\gamma_{i(\alpha )} < \alpha < \gamma_{i(\alpha )+1}$.  Now if  
$\alpha \notin \text{ Rang}(F_{i(\alpha)})$ then $C_\alpha = 
\emptyset$ and we are done so for some $x,\beta$ we have $\alpha = 
F_{i(\alpha)}(x,\beta)$ hence necessarily $x \in \dbcu_{\epsilon \le
\gamma_{i(\alpha)}}{\Cal P}_\epsilon$ and $\beta <
\gamma_{i(\alpha)}$.  By the definition of  $C_\alpha $ (the closure in the 
order topology on $\alpha$, of the set of $C^-_\alpha$ i.e. the set of  
$F_j(x \cap \zeta,\beta)$ for the pair $(j,\zeta)$ satisfying 
$\boxtimes^{x,\beta}_{j,\zeta}$ it 
suffices to show $C^-_\alpha \subseteq \alpha$, i.e. 
\mr
\item "{$(*)$}"  if the pair $(j,\zeta)$ satisfies
$\boxtimes^{x,beta}_{j,\zeta}$ then $F_j(x \cap \zeta,\beta) < \alpha$. 
\ermn
So assume $(j,\zeta)$ satisfies $\boxtimes^{x,\beta}_{j,\zeta}$ 
but by clause (iii) we know that
$j < i(\alpha )$ and so Rang$(F_j) \subseteq  \gamma _{j+1} \subseteq  
\gamma _{i(\alpha )} < \alpha$ as required.
\mn
\ub{Clause $(b)$}:  If $\alpha^* \in \text{ nacc}(C_\alpha)$ then  
$C_{\alpha^\ast} = C_\alpha \cap \alpha^\ast$. 

If it is enough to show  $C^-_{\alpha^*} = \alpha^* \cap C^-_\alpha$ 
and as $C^-_\alpha  = \text{ nacc}(C_\alpha)$, we have $\alpha^*
\in C^-_\alpha$.  As $\alpha^* \in C^-_\alpha$ necessarily for some  
$\zeta,j$  satisfying $\boxtimes^{x,\beta}_{j,\zeta}$ 
we have $\alpha^* = F_j(x \cap \zeta,\beta)$.  
By the choice of  $F_j$ necessarily $\alpha^*$ is a 
successor ordinal and  $\gamma _j < \alpha ^\ast  < \gamma _{j+1}.$ 

Now any member $\alpha(1)$ of $\alpha^* \cap C^-_\alpha $ has the 
form  $F_{j(1)}(x \cap  \zeta (1),\beta )$  with  $j(1)$, $\zeta (1)$  
satisfying $\boxtimes^{x,\beta}_{j,\zeta}$; 
clearly  $\gamma _{j(1)} < \alpha (1) = F_{j(\ast )}(x 
\cap  \zeta (1),\beta ) < \gamma _{j(1)+1}$ and  $\gamma _j < \alpha ^\ast  = 
F_j(x \cap  \zeta ,\beta ) < \gamma _{j+1}$.  But  $\alpha (1) < \alpha ^\ast $
(being in  $\alpha ^\ast  \cap  C^-_\alpha )$  so necessarily  $j(1) + 1 \leq  
j$.  So  $j(1)$, $\zeta (1)$  satisfy $(i)-(v)$ with  $x$  
replaced by $x \cap \zeta$, i.e., satisfy
$\boxtimes^{x,\beta}_{j,\zeta}$; recall by $\alpha^* = 
F_j(x \cap  \zeta ,\beta )$,  so  
$F_{j(x)}(x \cap  \zeta (1),\beta ) \in C^-_{\alpha^*}$.  So  
$\alpha ^\ast  \cap  C^-_\alpha  \subseteq  C^-_{\alpha ^\ast }$;  similarly  
$C^-_{\alpha^*} \subseteq \alpha^* \cap C^-_\alpha $,  so we get the 
desired equality.
\mn
\ub{Clause $(c)$}:  We shall show that $E = \{\gamma_i:i$ is a limit
ordinal $< \lambda\}$ is as required in closed (c).

Clearly $E$ is a club of $\lambda$.  So assume that 
$\delta \in A \cap E$ we should prove:  cf$(\delta) < \delta, 
\delta = \sup(C_\delta)$, cf$(\delta) = \text{ otp}(C_\delta)$. \nl
Now $\delta \in E \cap A \Rightarrow \delta > \text{ cf}(\delta)$ 
holds as we assume $\gamma_i \in A \Rightarrow 
\text{ cf}(\gamma_i) < \gamma _i$.  
As $\delta \in E$, by $E$'s definition for some limit ordinal $i(*)$
we have  $\delta  = \gamma_{i(*)}$.  
By the choice of  $C_\delta $ it is enough to find a set  
$C$  closed unbounded in  $\delta $  of order type cf$(\delta)$ such that  
$\alpha \in \text{ nacc}(C) \Rightarrow \alpha$ successor 
$\and C_\alpha = C \cap \alpha$. \nl
By the choice of $\bar{\Cal P}$,  for some $x \subseteq \delta$, 
otp$(x) < \delta = \sup(x)$ and $\dsize \bigwedge_{\zeta < \delta} x \cap
\zeta \in \dbcu_{\gamma < \delta} {\Cal P}_\gamma$. \nl
By $(\ast)$ above also $\xi \in x \and \bar S \in x \backslash \xi
\Rightarrow x \cap \zeta \backslash \xi \in \dbcu_{\gamma < \delta}
{\Cal P}_\gamma$ so without loss of generality otp$(x) < \text{
Min}(x)$.  Let $\beta = \text{ otp}(x)$,  so we know that $\beta$ is a
limit ordinal, moreover cf$(\beta) = \text{ cf}(\delta)$. 
Remember $e_\beta$ is a club of $\beta$ of order type cf$(\beta)$ which
is cf$(\delta)$.  Let

$$
y =: \{\zeta \in x:\text{ otp}(x \cap \zeta) \in e_\beta\}.
$$
\mn
Clearly  $y$  is a subset of  $x$  of order type otp$(e_\beta ) = 
\text{ cf}(\delta)$.  Define $h:y \rightarrow i(*)$ by $h(\zeta) =
\text{ Min}\{j:x \cap \zeta \in \dbcu_{\epsilon \le \gamma_j} {\Cal
P}_\epsilon\}$,  so by $(\ast )$ we know that  
$h$  is non-decreasing, and by the choice of $x,\dsize
\bigwedge_{\zeta \in y} \gamma _{h(\zeta)} < \delta$, equivalently 
$\dsize \bigwedge_{\zeta \in y} h(\zeta) < i(\ast)$. \nl
Let  $z = \{\zeta  \in  y:for$ every  $\xi  \in  y \cap  \zeta $  we have  
$h(\xi ) < h(\zeta )\}$.  Let  $C^- = \left\lbrace F_{h(\zeta )}(x \cap  
\zeta ,\beta ):\zeta  \in  z\right\rbrace $;  it satisfies:  $C^- \subseteq  
\delta  = \sup {}^\alpha \delta_\alpha$ and it is easy to check, as in
the proof of clause (c) that $[\alpha \in C^- \Rightarrow  C^-_\alpha  = C^- 
\cap \alpha]$.  So by the choice of $C^-$ its closure in $\delta$  
is as required.
\mn
\ub{Clause $(d)$}:  nacc$(C_\alpha)$ is a set of successor ordinals. 

Check.
\enddemo
\bigskip

\remark{Remark}  1) We could also 
strengthen $(*)$ to make $z \cap \zeta \in {\Cal P}_{h(\zeta)}$. 
\nl
2) By Definition \scite{1.1} we know that  
$I[\lambda]$ is an ideal; by \scite{1.2}(1) we know that $I[\lambda]$
includes the ideal of non-stationary subsets of $\lambda$.  By the last
phrase and Definition \scite{1.1}, clearly $I[\lambda]$ is normal.
\hfill$\square_{\scite{1.2}}$
\endremark
\bigskip

\proclaim{\stag{1.3} Claim}  If $\kappa,\lambda$ are regular, 
$S \subseteq \{\delta < \lambda:{\text{\rm cf\/}}(\delta) = \kappa\},
S \in I[\lambda],S$ stationary, $\kappa^+ < \lambda$ \ub{then} 
we can find  $\bar {\Cal P} = 
\langle {\Cal P}_\alpha:\alpha < \lambda \rangle$ such that for  
$\delta(\ast) =: \kappa$ we have: 
\mr
\item "{$\oplus^{\lambda,\delta(*)}_{{\Cal P}_S}(i)$}"  ${\Cal P}_\alpha$ is a 
family of closed subsets of  $\alpha,|{\Cal P}_\alpha| < \lambda$
\sn
\item "{$(ii)$}"   ${\text{\rm otp\/}}(C) \le \delta(*)$ 
for $C \in \dbcu_\alpha {\Cal P}_\alpha$
\sn
\item "{$(iii)$}"  for some club $E$ of $\lambda$,  we have:        
\nl
$[\alpha  \notin  E \Rightarrow {\Cal P}_\alpha  = \emptyset ]$  and        
\nl
$[\alpha  \in  E \Rightarrow  
(\forall C \in {\Cal P}_\alpha)({\text{\rm otp\/}}(C) \le \delta(*))]$        
\nl
$[\alpha \in E \backslash (S \cap { \text{\rm acc\/}}(E)) 
\Rightarrow (\forall C \in {\Cal P}_\alpha)[{\text{\rm otp\/}}(C) 
< \delta(*)]$        
\nl
$[\alpha \in S \cap { \text{\rm acc\/}}(E) \Rightarrow (\exists ! C \in  
{\Cal P}_\alpha)({\text{\rm otp\/}}(C) = \delta(*))]$        
\nl
$[\alpha \in S \cap { \text{\rm acc\/}}(E) \and C \in {\Cal P}_\alpha \and
{ \text{\rm otp\/}}(C) = \delta(*) \Rightarrow \alpha = \sup(C))]$     
\sn
\item "{$(iv)$}"   $C \in {\Cal P}_\alpha \and \beta \in { \text{\rm nacc\/}}
(C) \Rightarrow \beta \cap C \in {\Cal P}_\beta$     
\sn
\item "{$(v)$}"  for any club $E'$ of $\lambda$ 
for some $\delta \in S \cap E'$ and 
$C \in {\Cal P}_\delta$ we have $C \subseteq E' \and { \text{\rm otp\/}}(C) = 
\delta(*)$.
\endroster
\endproclaim
\bigskip

\demo{Proof}  Let $\langle C_\alpha:\alpha < \lambda \rangle$ witness $``S \in 
I[\lambda]"$ be as in 
\scite{1.2}(1); without loss of generality otp$(C_\alpha) \le
\delta(*)$.  
For any club $E$, consisting of limit ordinals for simplicity, let 
us define  ${\Cal P}^\alpha_E$ by induction on $\alpha < \lambda$:

$$
\align
{\Cal P}^\alpha_E =: &\{ \alpha \cap g \ell(C_\beta,E):\alpha \in E
\text{ and } \alpha \le \beta  < \text{ Min}[E 
\backslash (\alpha +1)]\} \\    
  & \cup \{C \cup \{\beta\}:\beta \in E \cap \alpha,
C \in {\Cal P}^\beta_E \text{ and otp}(C) < \delta(*)\}
\endalign
$$
\mn
where 

$$
g\ell(C_\beta,E) =: \{ \sup(E \cap (\gamma +1)):\gamma \in C_\beta
\text{ and } \gamma > \text{ Min}(E)\}.
$$
\mn
Note that $|{\Cal P}^\alpha_E| \le |\text{Min}(E \backslash (\alpha +1)| < 
\lambda$. \nl
We can prove that for some club $E$ of $\lambda$ the sequence   
$\langle {\Cal P}^\alpha _E:\alpha < \lambda \rangle$ is as
required except possibly clause $(v)$ which can be corrected gotten by
a right of $E$ (just by trying
successively $\kappa^+$ clubs $E_\zeta$ (for $\zeta < \kappa^+$) decreasing 
with $\zeta$,  see \cite{Sh:365}).  Note that clause (iv) 
guaranteed by demanding  $E$  to 
consist of limit ordinals only and the second set in the union defining  
${\Cal P}^\alpha_E$.  \hfill$\square_{\scite{1.3}}$
\enddemo
\bn
The following lemma gives sufficient condition for the existence of ``quite 
large" stationary sets in  $I[\lambda ]$  of almost any fixed cofinality.
\proclaim{\stag{1.4} Lemma}  Suppose 
\mr
\item "{$(i)$}"   $\lambda > \kappa > \aleph_0,\lambda$ and $\kappa$ are regular 
\sn
\item "{$(ii)$}"  $\bar{\Cal P} = \langle {\Cal P}_\alpha:\alpha < \kappa \rangle$,  
${\Cal P}_\alpha$ a family of $< \lambda$ closed subsets of $\alpha$ 
\sn
\item "{$(iii)$}"  $I_{\bar{\Cal P}} =: \{S \subseteq \kappa:\text{for
some club } E \text{ of } \kappa \text{ for no }  \delta \in S \cap E$        
is there a club $C$ of $\delta$, such that $C \subseteq E$ and        
$[\alpha  \in { \text{\rm nacc\/}}(C) \Rightarrow C \cap \alpha \in
\dbcu_{\beta < \alpha} {\Cal P}_\beta]\}$ is a proper ideal on $\kappa$.
\ermn
\ub{Then} there is $S^\ast \in I[\lambda]$ such that for stationarily many  
$\delta < \lambda$ of cofinality $\kappa,S^\ast \cap \delta$ is 
stationary in $\delta$,  moreover for some club $E$ of $\delta$ of order
type $\kappa$

$$
\{{\text{\rm otp\/}}(\alpha \cap E):\alpha \in E \backslash S^\ast \} 
\in I_{\Cal P}.
$$
\endproclaim
\bigskip

\remark{\stag{1.4A} Remark}  1) The 
``for stationarily many" in the conclusion can be 
strengthened to: a set whose complement is in the ideal defined in 
\cite[\S2]{Sh:371}. \nl
2) So if $\kappa^\sigma < \lambda$ then we can have $\{i <
\kappa:\text{cf}(i) = \sigma\} \in I_{\bar{\Cal P}}$.
\endremark
\bigskip

\demo{Proof}  Let $\chi$ be regular large enough, $N^*$ be an elementary 
submodel of $({\Cal H}(\chi ),\in,<^*_\chi)$ of cardinality  $\lambda $  such 
that  $(\lambda  + 1) \subseteq  N^\ast $,  $\bar{\Cal P} \in  N$.  Let  
$\bar C = \langle C_i:i < \lambda \rangle $  list  $N^\ast  \cap  \{A 
\subseteq  \lambda :|A| < \kappa \}$  and let 

$$
\align
S^\ast  = \{ \delta  < \lambda:&\text{ cf}(\delta ) < \kappa \text{
and for some } A \subseteq  \delta \text{ satisfying } \delta =
\sup(A), \text{ we have} \\
  &\text{ otp}(A) < \kappa 
\text{ and } (\forall \alpha < \delta)[A \cap \alpha \in \{C_i:i < \delta\}]\}.
\endalign
$$
\mn
Clearly $S^* \in I[\lambda]$; so we should only find enough $\delta <
\lambda$ of cofinality $\kappa$ as required in the conclusion of 
\scite{1.4}.  So let $E^*$ be a club of $\lambda$ and we shall prove
that such $\delta \in E^*$ exists.  
We can choose $M_\zeta$ by induction on $\zeta \le \kappa$  
such that:
\mr
\item "{$(a)$}"   $M_\zeta \prec ({\Cal H}(\chi),\in,<^* _\chi)$ 
\sn
\item "{$(b)$}"  $\|M_\zeta \| < \lambda,M_\zeta \cap \lambda$ an ordinal 
\sn
\item "{$(c)$}"  $M_\zeta $ is increasing continuous 
\sn
\item "{$(d)$}"  $N,\kappa,\bar{\Cal P},\bar C,E^*$ belongs to $M_0$ 
\sn
\item "{$(e)$}"  $\langle M_\epsilon:\epsilon \le \zeta \rangle  
\in M_{\zeta +1}$.
\ermn
Let $\delta_\zeta  = \sup(M_\zeta \cap \lambda)$, clearly
$\delta_\zeta \in E^*$ for every $\zeta \le \kappa$ and
$\langle \delta_\zeta:\zeta \le \kappa \rangle$ is a (strictly) increasing 
continuous, so  $\delta  =: \delta _\kappa $ has cofinality  $\kappa $.  Hence 
there is a (strictly) increasing continuous sequence  
$\langle \alpha _\zeta:\zeta < \kappa \rangle \in N^*$ with limit  
$\delta$, and clearly $E = \{\zeta < \kappa:\alpha_\zeta = 
\delta_\zeta$ and $\zeta$ is a limit ordinal$\}$ 
is a club of $\kappa$.  We know that 

$$
\align
T =: \{ \zeta < \kappa:&\zeta \in E \text{ and for some club } C \text{ of }
\zeta,C \subseteq E \text{ and} \\
  & \dsize \bigwedge_{\epsilon < \zeta}[C \cap \epsilon \in \dbcu_{\xi
< \zeta} {\Cal P}_\xi]\}.
\endalign
$$
\mn
is stationary; moreover,  $\kappa \backslash T \in I_{\bar{\Cal P}}$ 
(see assumption (iii)) and clearly  $T \subseteq E$. \nl
Clearly it suffices to show 
\mr
\item "{$(*)$}"   $\zeta \in T \Rightarrow \delta_\zeta \in S^*$.
\ermn
Suppose $\zeta  \in  T$,  so there is  $C$,  a club of  $\zeta $  such that  
$C \subseteq  E$  and  $\dsize \bigwedge_{\epsilon < \zeta} 
[C \cap \epsilon \in \dbcu_{\xi < \zeta} {\Cal P}_\xi]$. 
Let $C^\ast  = \{\delta _\epsilon :\epsilon  
\in  C\}$,  so  $C^\ast $ is a club of  $\delta _\zeta $ of order type  $\leq  
\zeta  < \kappa $  (which $is  < \delta _0 \leq   \delta _\zeta )$.  It 
suffices to show for  $\xi  \in  C$  that  $\{\delta _\epsilon :\epsilon  \in  
\xi  \cap  C\} \in  \{C_i:i < \delta _\zeta \}.$
\nl
For this end we shall show 
\mr
\item "{$(\alpha)$}"   $\{\delta_\epsilon :\epsilon  \in  C \cap  \xi \} \in  \{C_i:i < 
\lambda \}$ 
\sn
\item "{$(\beta)$}"  $\{\delta_\epsilon:\epsilon \in C \cap  \xi \} \in  
M_{\xi +1}$.
\ermn
This suffices as  $\langle C_i:i < \lambda \rangle  \in  M_0 \prec  M_{\xi +1}$
and  $M_{\xi +1} \cap  \{C_i:i < \lambda \} = \{C_i:i \in  \lambda  \cap  
M_{\xi +1}\} = \{C_i:i < \delta _{\xi +1}\}.$
\enddemo
\bigskip

\demo{\ub{Proof of $(\alpha)$}}   Remember 
$\langle \alpha_\epsilon:\epsilon < \kappa \rangle \in N^*$.  
Also $\bar{\Cal P} = \langle {\Cal P}_\epsilon :\epsilon  < 
\kappa \rangle \in N^\ast $
hence $\dbcu_{\epsilon < \kappa} {\Cal P}_\epsilon  \subseteq N^\ast$
(as $\kappa < \lambda,|{\Cal P}_\epsilon| < \lambda,\lambda +1
\subseteq N,\bar{\Cal P} \in N^*$ so now for $\xi \in C$ we have
$C \cap  \xi \in \dbcu_{\epsilon  < \kappa} {\Cal
P}_\epsilon$;  hence $C \cap  \xi  \in  N^\ast $.  Together  
$\{\alpha _\epsilon :\epsilon  \in  \xi  \cap  C\} \in  N^\ast $;  as  
$\epsilon  \in  C \Rightarrow  \epsilon  \in  E \Rightarrow  \alpha _\epsilon  
= \delta _\epsilon $ (as  $C \subseteq  E$  and the definition of  $E)$,  and 
the definition of $\langle C_i:i < \lambda \rangle $,  we are done.
\enddemo
\bigskip

\demo{\ub{Proof of $(\beta)$}}   We know $\bar {\Cal P} \in  M_0$;  as  
$|{\Cal P}_\epsilon| < \lambda,\kappa < \lambda$ clearly
$|\dbcu_{\epsilon < \kappa} {\Cal P}_\epsilon| < \lambda$ 
so as $M_\epsilon  \cap  \lambda $  
is an ordinal, clearly  $\dbcu_{\epsilon < \kappa} {\Cal P}_\epsilon
\subseteq M_0$.  So for $\epsilon < \zeta$ we have $C \cap  \epsilon \in
\dbcu_{\gamma < \zeta} {\Cal P}_\gamma  \subseteq  M_0 
\subseteq  M_{\xi +1}$.  As  
$\langle M_i:i \le \xi \rangle  \in  M_{\xi +1}$ clearly  
$\langle \delta _i:i \le \xi \rangle  \in  M_{\xi +1}$ hence by the previous 
sentence also 
$\langle \delta_i:i \in C \cap  \xi \rangle  \in  M_{\xi +1}$,  as 
required. \hfill$\square_{\scite{1.4}}$
\enddemo
\bigskip

\demo{\stag{1.5} Conclusion}  If $\kappa $,  $\lambda$ are regular,  
$\kappa^+ < \lambda$  \ub{then} there is a stationary  
$S \subseteq \{\delta < \lambda:\text{cf}(\delta)  = \kappa\}$  
in $I[\lambda]$.
\enddemo
\bigskip

\demo{Proof}  If $\lambda = \kappa^{++}$ - use \cite[4.1]{Sh:351}.  
So assume $\lambda > \kappa^{++}$.  
By \cite[4.1]{Sh:351} the pair $(\kappa,\kappa^{++})$  satisfies 
the assumption of \scite{1.3} for 
$S = \{\delta < \kappa^{++}:\text{cf}(\delta) = 
\kappa\}$; (i.e. $\kappa$, $\lambda $  there stands for $\kappa$, 
$\kappa ^{++}$ here).  Hence the conclusion of \scite{1.3} holds for some  
$\bar{\Cal P} = \langle {\Cal P}_\alpha:\alpha < \kappa^{++}\rangle$,  
$|{\Cal P}_\alpha| < \kappa^{++}$.  Now apply \scite{1.4} with  
$(\kappa ^{++},\lambda)$ here standing for $(\kappa,\lambda)$ there (we 
have just proved  
$I_{\bar{\Cal P}}$ is a proper ideal, so assumption (ii) holds).  Note: 
\mr
\item "{$(*)$}"  $\{\delta < \kappa^{++}:\text{cf}(\delta) = \kappa\}
\notin I_{\bar{\Cal P}}$.
\ermn
Now the conclusion of \scite{1.4} 
(see the moreover and choice of  $\bar{\Cal P}$  
i.e. $(*)$) gives the desired conclusion.   \hfill$\square_{\scite{1.5}}$
\enddemo
\bigskip

\demo{\stag{1.6} Conclusion}  If $\lambda > \kappa$ are 
uncountable regular, $\kappa^+ < \lambda$, \ub{then} for 
some stationary  $S \subseteq \{\delta < \lambda:\text{cf}(\delta)  
= \kappa \}$ and some $\bar{\Cal P} = 
\langle {\Cal P}_\alpha:\alpha < \lambda \rangle$ we have:  
$\oplus^{\lambda,\kappa}_{{\Cal P},S}$ from the conclusion of 
\scite{1.3} holds.
\enddemo
\bigskip

\demo{Proof}   As $\kappa$ is regular apply \scite{1.5} and then
\scite{1.3}. 
\hfill$\square_{\scite{1.6}}$
\enddemo
\bn
Now \scite{1.6} 
was a statement I have long wanted to know, still sometimes we want to 
have  $``C_\delta \subseteq E$, otp$(C) = \delta (\ast)"$,
$\delta(\ast)$ not a regular cardinal.  We 
shall deal with such problems.
\proclaim{\stag{1.7} Claim}  Suppose 
\mr
\widestnumber\item{$(iii)$}
\item "{$(i)$}"  $\lambda > \kappa > \aleph_0,\lambda$ and  
$\kappa$ are regular cardinals 
\sn
\item "{$(ii)$}"   $\bar{\Cal P}_\ell = 
\langle {\Cal P}_{\ell,\alpha}:\alpha < \kappa \rangle$ 
for $\ell  = 1,2$,  where  ${\Cal P}_{1,\alpha }$ is a 
family of $< \lambda$ closed subsets of  $\alpha $,  ${\Cal P}_{2,\alpha }$ 
is a family of $\le \lambda$ clubs of $\alpha$ and  $[C \in {\Cal
P}_{2,\alpha} \and \beta  \in  C \Rightarrow  C \cap \beta \in
\dbcu_{\gamma < \alpha} {\Cal P}_{1,\gamma}]$ 
\sn
\item "{$(iii)$}"   $I_{\bar{\Cal P}_1,\bar{\Cal P}_2} =: \{ S \subseteq
\kappa:\text{ for some club } E \text{ of } \kappa$  for no      
$\delta  \in  S \cap  E$  is there  $C \in  {\Cal P}_{2,\alpha}$, $C
\subseteq E\}$  is a proper ideal on $\kappa$.
\ermn
\ub{Then} we can find $\bar{\Cal P}^\ast_\ell  = 
\langle {\Cal P}^\ast_{\ell,\alpha}:\alpha  < \lambda \rangle$ for $\ell 
= 1,2$ such that: 
\mr
\item "{$(A)$}"  ${\Cal P}^*_{1,\alpha}$ is 
a family of $< \lambda$ closed subsets of $\alpha$ 
\sn
\item "{$(B)$}"   $\beta \in { \text{\rm nacc\/}}(C) \and C \in 
{\Cal P}^*_{1,\alpha} \Rightarrow  C \cap \beta \in {\Cal P}^\ast_{1,\beta}$ 
\sn
\item "{$(C)$}"   ${\Cal P}^*_{2,\delta}$ is 
a family of $\le \lambda$ clubs of  
$\delta$ (for $\delta$ limit $< \lambda$ such that)
$[\beta \in { \text{\rm nacc\/}}(C) \and C \in 
{\Cal P}^\ast_{2,\delta} \Rightarrow  C \cap  \beta  \in  
{\Cal P}^\ast_{1,\beta}]$ 
\sn
\item "{$(D)$}"   for every club $E$ of $\lambda$ for 
some strictly increasing continuous sequence \nl
$\langle \delta _\zeta:\zeta \le \kappa \rangle$ of ordinals $<
\lambda$ we have $\{ \zeta < \kappa:\zeta \text{ limit, and for some }C \in  
{\Cal P}_{2,\zeta}$ we have: \nl     
$\{\delta _\epsilon:\epsilon \in C\} \in {\Cal P}^\ast_{2,\delta_\zeta}$
(hence $[\xi \in { \text{\rm nacc\/}}(C) \Rightarrow  
\{\delta_\epsilon:\epsilon \in C \cap \xi\} \in 
{\Cal P}^\ast _{1,\delta _\xi}]\} \equiv \kappa
\text{ mod }\ I_{\bar{\Cal P}_1,\bar{\Cal P}_2}$ 
\sn
\item "{$(E)$}"   we have $e_\delta$ a club of 
$\delta$ of order type ${\text{\rm cf\/}}(\delta)$  
for any limit $\delta < \lambda$; such that for any  $C \in
\dbcu_{\alpha < \lambda}{\Cal P}^*_{2,\alpha}$ for some $\delta  < \lambda,  
{\text{\rm cf\/}}(\delta) = \kappa$ and $C' \in \dbcu_{\beta < \kappa} {\Cal
P}_{2,\beta}$ we have $C = \{\gamma \in e_\delta:{\text{\rm otp\/}}
(e_\delta \cap \gamma) \in C'\}$.
\ermn
\endproclaim
\bigskip

\demo{Proof}  Same proof as \scite{1.4}.  
(Note that without loss of generality  $[C \in  
{\Cal P}_{1,\alpha} \and \beta  < \alpha < \kappa 
\Rightarrow C \cap \beta \in {\Cal P}_{1,\beta}])$.
\enddemo
\bigskip

\demo{\stag{1.8} Conclusion}  If $\delta(\ast)$ is a limit ordinal and
$\lambda = \text{ cf}(\lambda) > |\delta(\ast)|^+$ \ub{then} 
we can find $\bar{\Cal P}^\ast_\ell  =
\langle {\Cal P}^\ast_{\ell,\alpha }:\alpha < \lambda \rangle $ for $\ell 
= 1,2$ and stationary $S \subseteq \{\delta <
\lambda:\text{cf}(\delta) = \text{ cf}(\delta(*))\}$ such that: 
\smallskip

$\oplus^{\lambda,\delta(*)}_{\bar{\Cal P}^*_1,\bar{\Cal P}^*_2}
\quad (A) \quad  {\Cal P}^\ast_{1,\alpha}$ is a family of $< \lambda$ closed 
subsets of $\alpha$ each of \nl

\hskip68pt  order type $< \delta(*)$ \nl

\hskip45pt $(B) \quad \beta  \in \text{ nacc}(C) \and 
C \in {\Cal P}^*_{1,\alpha} \Rightarrow  C \cap  \beta 
\in {\Cal P}^*_{1,\beta}$  \nl

\hskip45pt $(C) \quad {\Cal P}^*_{2,\delta}$ is a family 
of $\le \lambda$ clubs of $\delta$ \nl

\hskip68pt (yes, maybe $= \lambda$) of order type \nl

\hskip68pt  $\delta(*)$, and $[\beta \in \text{ nacc}(C) \and 
C \in {\Cal P}^*_{2,\delta} \Rightarrow C 
\cap \beta \in {\Cal P}^\ast_{1,\beta}]$ \nl

\hskip45pt $(D) \quad$ for every club $E$ of 
$\lambda$ for some $\delta  \in  E \cap  S$,  \nl

\hskip68pt cf$(\delta) = \text{ cf}(\delta(\ast))$ and there is $C \in {\Cal
P}^\ast_{2,\beta}$ such that $C \subseteq E$.
\enddemo
\bigskip

\demo{Proof}  If $\lambda = |\delta(*)|^{++}$ (or any successor of regulars) 
use \cite[ChIII,6.4]{Sh:e}(2) or \cite[2.14]{Sh:365}(2)((c)+(d)). 

If $\lambda > |\delta(*)|^{++}$ let $\kappa  = |\delta (\ast )|^{++}$ 
and let $S_1 = \{\delta  < \kappa^{++}:\text{cf}(\delta) 
= \text{ cf}(\delta(*))\}$;  
applying the previous sentence we get  $\bar{\Cal P}^\ast _1$,  
$\bar{\Cal P}^\ast _2$ satisfying $\oplus^{\kappa^{++},\delta(\ast)}
_{\bar{\Cal P}^*_1,\bar{\Cal P}^*_2,S_1}$, hence satisfying the assumption
of \scite{1.7} so we can apply \scite{1.7}.  \hfill$\square_{\scite{1.8}}$
\enddemo
\bigskip

\definition{\stag{1.9} Definition}   ${}^+
\oplus^{\lambda,\delta(*)}_{\bar{\Cal P}_1,\bar{\Cal P}_{2,S}}$ is defined
as in \scite{1.8} except that we replace $(C)$ by 
\mr
\item "{$(C)^+$}"   ${\Cal P}^\ast_{2,\delta}$ is a family 
of $< \lambda$ clubs of $\delta$ of order type $\delta(\ast)$.
\endroster
\enddefinition
\bigskip

\remark{\stag{1.9A} Remark}  Note that if 
${\Cal P}_\alpha = {\Cal P}_{1,\alpha} \cup  
{\Cal P}_{2,\alpha}$,  $|{\Cal P}_{2,\alpha}| \le 1$, ${\Cal
P}_{1,\alpha} = \{C \in {\Cal P}_\alpha:\text{otp}(C) < \delta(*)\},
{\Cal P}_{2,\alpha} = \{C \in {\Cal P}_\alpha:\text{otp}(C) = \delta(*)\}$
then ${}^+ \oplus^{\lambda,\delta(*)}_{\bar{\Cal P}_1,{\bar\Cal P}_{2,S}}
\Leftrightarrow \oplus^{\lambda,\delta(\ast)}_{\bar{\Cal P}_S}$ mod.
\endremark
\bigskip

\proclaim{\stag{1.10} Claim}   Suppose $\lambda = { \text{\rm cf\/}}
(\lambda) > |\delta(\ast)|^+$,  
$\delta(\ast)$ a limit ordinal,  additively indecomposable (i.e. $\alpha 
< \delta (\ast ) \Rightarrow \alpha + \alpha < \delta(\ast))$,  
$\oplus^{\lambda,\delta(*)}_{\bar{\Cal P}_1,\bar{\Cal P}_{2,S}}$ 
from \scite{1.8} and 
\mr
\item "{$(*)$}"   $\alpha \in S \Rightarrow |{\Cal P}_{2,\alpha}| \le
|\alpha|$.
\ermn
(Note: a non-stationary subset of $S$ does not count; e.g. for  $\lambda $
successor cardinal the $\alpha$ with $|\alpha |^+ < \lambda $.  Note:  
${}^+\oplus^{\lambda,\delta(\ast)}_{\bar{\Cal P}_1,\bar{\Cal P}_{2,S}}$ holds
by $(\ast )$ and if $\lambda$  is successor then  
${}^+ \oplus^{\lambda,\delta(\ast)}_{\bar{\Cal P}_1,\bar{\Cal
P}_{2,S}}$ suffice).
\nl
\ub{Then} for some stationary  $S_1 \subseteq S$ and $\bar{\Cal P} =
\langle {\Cal P}_\alpha:\alpha  < \lambda \rangle $ we have:  
${\Cal P}_\alpha \subseteq {\Cal P}_{1,\alpha} \cup {\Cal P}_{2,\alpha}$ and: 
\mr
\item "{${}^*\otimes^{\lambda,\delta(*)}_{\bar{\Cal P},S_1}$}"
 (i)  ${\Cal P}_\alpha $ is a family of closed subsets of $\alpha $, 
$|{\Cal P}_\alpha | < \lambda$ \nl
\hskip40pt  $(ii) \,{ \text{\rm otp\/}}C < \delta(*)$ if 
$C \in {\Cal P}_\alpha,\alpha \notin S_1$  \nl
\hskip40pt  $(iii)$  if $\alpha \in S_1$ then:  
${\Cal P}_\alpha = \{C_\alpha\},{\text{\rm otp\/}}(C_\alpha) = 
\delta(\ast)$, \nl

\hskip42pt  $C_\alpha$ a club of $\alpha$ disjoint to $S_1$ \nl
\hskip40pt  $(iv)$   $C \in {\Cal P}_\alpha \and \beta 
\in { \text{\rm nacc\/}}(C) \Rightarrow \beta \cap C \in {\Cal P}_\beta$ \nl
\hskip40pt $(v)$  for any club $E$ of $\lambda$ for some 
$\delta \in S_1$ we have $C_\delta \subseteq E$.
\endroster
\endproclaim
\bigskip

\remark{\stag{1.10A} Remark}   Note there 
are two points we gain: for  $\alpha \in S_1$,  
${\Cal P}_\alpha $ is a singleton (similarly to \scite{1.3} where we have  
$(\exists^{\le 1} C \in  {\Cal P}_\delta)[\text{otp}(C) = \delta(*)])$, and an 
ordinal $\alpha$ cannot have a double role $-C_\alpha$ a guess (i.e.  
$\alpha \in S_1)$ and $C_\alpha$ is a proper initial segment of such  
$C_\delta$.  When $\delta(*)$ is a regular cardinal this is easier.
\endremark
\bigskip

\demo{Proof}  Let ${\Cal P}_{2,\alpha } = 
\{C_{\alpha ,i}:i < \alpha \}$  (such a 
list exists as we have assumed  $|{\Cal P}_{2,\alpha }| \le |\alpha|$,  we 
ignore the case ${\Cal P}_{2,\alpha} = \emptyset )$.  Now 
\mr
\item "{$(*)_0$}"  for some $i < \lambda$ for every club 
$E$ of $\lambda$ for some $\delta \in S \cap E$ we have 
$C_{\delta ,i }\backslash E$  is 
bounded in  $\alpha$ \nl
[Why?  If not, for every $i < \lambda$ there is a club  $E_i$ of  $\lambda $
such that for no  $\delta  \in  S \cap E$ is $C_{\delta,i} \backslash E$  
bounded in $\alpha$.  Let $E^\ast = \{j < \lambda:j$ a limit ordinal, 
$j \in \dbca_{i < j} E_i\}$,  it is a club of  $\lambda $,  hence for
some  $\delta  \in  S \cap  E^\ast$ and  $C \in  {\Cal P}_{2,\delta}$
we have  $C \subseteq E^\ast$.  So for some $i < \alpha,C =
C_{\delta,i}$,  so $C \subseteq E^\ast \subseteq E_i \cup i$ hence $C_{\delta ,i}\backslash i 
\subseteq E_i$,  contradicting the choice of  $E_i$.]. 
\sn
\item "{$(*)_1$}"   for some $i < \lambda$ and 
$\gamma < \delta(\ast)$,  letting
$C_\delta =: C_{\delta,i} \backslash 
\{\zeta \in C_{\delta,i}:\text{otp}(\zeta  
\cap  C_{\delta ,i}) < \gamma \}$  we have: for every club  
$E$  of  $\lambda$ for some  $\delta \in  S \cap  E$  we have:  $C_\delta \subseteq E$ \nl   
[Why?  Let  $i(\ast )$  be as in $(\ast)_0$, and for each  $\gamma  < 
\delta (\ast )$  suppose  $E_\gamma $ exemplify the failure of $(\ast )_1$ for 
$i(\ast )$  and  $\gamma $,  now  $\dbca_{\gamma < \delta(*)}
E_\gamma$ is a club of  $\lambda $  exemplifying the failure 
of $(\ast)_0$ for  $i(\ast )$  contradiction.  So for some  $\gamma  < 
\delta (\ast)$  we succeed.] 
\sn
\item "{$(*)_2$}"   Without 
loss of generality  $|{\Cal P}_{2,\alpha}| \le 1$,  so 
let  ${\Cal P}_{2,\alpha} = \{C_\alpha \}$ \nl   
[Why? \ Let $i,\gamma$ and $C_\delta $ (for $\delta \in S$)  be as in
$(\ast )_1$ and use  ${\Cal P}'_{1,\alpha} = \{ C \backslash \{\zeta
\in C:\text{otp}(\zeta \cap C) < \gamma \}:C \in {\Cal P}_{1,\alpha}\},
{\Cal P}'_{2,i} = \{C_\delta \}$.]
\sn
\item "{$(*)_3$}"   for some 
$h:\lambda \rightarrow |\delta(\ast)|^+$,  for every 
$\alpha \in S$ we have $h(\alpha) \notin \{h(\beta):\beta \in C_\alpha\}$   
\nl
[Why?  Choose $h(\alpha)$ by induction on $\alpha$.] 
\sn
\item "{$(*)_4$}"   for some $\beta < |\delta(*)|^+$ for every 
club $E$ of $\lambda$,  for some  
$\delta \in S \cap  h^{-1}(\{\beta\}),C_\delta \subseteq E$
\nl   
[Why?  If for each $\beta$ there is a counterexample  $E_\beta $ then  
$\cap \{E_\beta :\beta < |\delta(\ast)|^+\}$ is a counterexample for 
$(\ast)_2$.] 
\ermn
Now we have gotten the desired conclusion. \hfill$\square_{\scite{1.10}}$
\enddemo
\bigskip

\proclaim{\stag{1.11} Claim}   If $S \subseteq \{\delta <
\lambda:{\text{\rm cf\/}}(\delta) = 
\kappa\},S \in  I[\lambda],\kappa^+ < \lambda  
= { \text{\rm cf\/}}(\lambda)$,  \ub{then} for some stationary  $S_1
\subseteq S$ and $\bar{\Cal P}_1$ we have  ${}^\ast
\oplus^{\lambda,\delta(\ast)}_{{\Cal P}_{1,S_1}}$.
\endproclaim
\bigskip

\demo{Proof}   Same proof as \scite{1.3} (plus $(\ast)_3,(\ast)_4$ in
the proof of \scite{1.8}). \hfill$\square_{\scite{1.11}}$
\enddemo
\bigskip

\proclaim{\stag{1.12} Claim}   Assume 
$\lambda = \mu^+$, $|\delta(\ast)| < \mu$ and ${\text{\rm
cf\/}}(\delta(*)) \ne { \text{\rm cf\/}}(\mu)$. 

\ub{Then} we can find stationary  $S \subseteq \{\delta <
\lambda:{\text{\rm cf\/}}(\delta) = { \text{\rm cf\/}}(\delta)(*)\}$ 
and $\bar{\Cal P}$ such that ${}^* \otimes^{\lambda,\delta(*)}
_{\bar{\Cal P},S}$.
\endproclaim
\bigskip

\remark{Remark}   This strengthens \scite{1.8}.
\endremark
\bigskip

\demo{Proof}   \ub{Case $(\alpha)$.$\mu$ regular}. 

By \cite[Ch.III,6.4]{Sh:e}(2), \cite[2.14]{Sh:365}(2)((c)+(d)).
\mn
\ub{Case $\beta$. $\mu$ singular}. 

Let $\theta =: \text{ cf}(\mu),\sigma =: |\delta(*)|^+ + \theta ^+$ and  
$\mu = \dsize \sum_{\zeta < \theta} \mu_\zeta,\langle \mu _\zeta :\zeta  < 
\theta \rangle$ strictly increasing,  $\mu _0 > \sigma $  and for each  
$\alpha < \lambda$ let $\alpha = \dbcu_{\zeta < \theta}
A_{\alpha,\zeta}$,  $\langle A_{\alpha ,\zeta }:\zeta  < 
\theta \rangle$ increasing,  $|A_{\alpha ,\zeta }| \leq  \mu _\zeta$. 

By \scite{1.6} there is a sequence 
$\bar{\Cal P} = \langle {\Cal P}_\alpha :\alpha <
\lambda \rangle$ and stationary $S_1 \subseteq \{\delta < \lambda:
\text{ cf}(\delta) = \sigma \}$  such that
$\oplus^{\lambda,\sigma}_{\bar{\Cal P},S_1}$ of \scite{1.3} holds.  
$Let \cup\{{\Cal P}_\alpha:\alpha  < \lambda \} 
\cup  \{\emptyset \}$  be  $\{C_\alpha :\alpha  < \lambda \}$ such that  
$C_\alpha  \subseteq \alpha$, $[\alpha \in  S_1 \Rightarrow  C_\alpha
\in  {\Cal P}_\alpha  \and \text{ otp}(C_\alpha) = \sigma]$ and  
$[\alpha  \notin  S_1 \Rightarrow  \text{ otp}(C_\alpha) < \sigma]$.  
For some club $E^*_1$ of $\lambda,[\alpha \in E^*_1 
\Rightarrow \dbcu_{\beta < \alpha} {\Cal P}_\beta = 
\{C_\beta:\beta < \alpha\}]$. 

Looking again at $\oplus^{\lambda,\sigma}_{\bar{\Cal P},S_1}$,  
we can assume  $S_1 \subseteq E^\ast_1 \and (\forall \delta)[\delta
\in S_1 \Rightarrow  C_\delta \subseteq E^\ast_1]\}$,  hence 
\mr
\item "{$(*)$}"   $\delta  \in  S_1 \and \alpha  \in \text{ nacc } 
C_\delta  \Rightarrow  \alpha  \cap  C_\delta  \in  \{C_\beta :\beta  < 
\text{ Min}(C_\delta \backslash (\alpha +1))\}.$
\ermn
So as we can replace every  $C_\alpha $ by  $\{\beta \in  
C_\alpha:\text{otp}(C_\alpha  \cap  \beta )\}$  is even, 
without loss of generality [because we can replace every $C_\alpha$ by
$\{\beta \in C_\alpha:\text{otp}(\beta \cap C_\alpha)$ is even$\}$,
\wilog \, (check)]
\mr
\item "{$(\ast)^+$}"    $\delta  \in  S_1 \and \alpha  \in 
\text{ nacc } C_\delta  \Rightarrow  
\alpha  \cap  C_\delta  \in  \{C_\beta :\beta  < \alpha \}$. 
\ermn
Without loss of generality  $[\beta  \in  A_{\alpha ,\zeta } \Rightarrow  
C_\beta \subseteq A_{\alpha,\zeta}]$  (just note  
$|C_\beta | \le \sigma < \mu _\zeta)$ and $\alpha \in  
A_{\beta ,\zeta } \Rightarrow  A_{\alpha ,\zeta } \subseteq
A_{\beta ,\zeta }$.  For  $\alpha  \in  S_1$ let  $C_\alpha  = 
\{\beta _{\alpha ,\epsilon }:\epsilon  < \sigma \}  
(\beta _{\alpha ,\epsilon }$ increasing in  $\epsilon )$  and let  
$\beta ^\ast _{\alpha ,\epsilon } \in  
[\beta _{\alpha ,\epsilon },\beta _{\alpha ,\epsilon +1})$  be mimimal such 
that  $C_\alpha  \cap  \beta _{\alpha ,\epsilon +1} = 
C_{\beta ^\ast _{\alpha ,\epsilon }}$ (exists as  $\delta  \in  S_1 
\Rightarrow  C_\delta  \subseteq  E^\ast _1)$.  Without loss of generality 
every  $C_\alpha $ is an initial segment of some  $C_\beta $,  $\beta  \in  
S_1$ (if not, we redefine it  as $\emptyset$). 
\mr
\item "{$(\ast)_1$}"   there are $\gamma = \gamma (\ast ) <  \theta $  and stationary  
$S_2 \subseteq S_1$ such that for every club  $E$  of  $\lambda $,  
for some  $\delta  \in  S_2$ we have:  $C_\delta  \subseteq E$,  
and for arbitrarily large  $\epsilon  < \sigma $,  
$\beta ^\ast _{\delta ,\epsilon} \in A_{\beta _{\delta ,\epsilon
+1},\gamma }$. \nl   
[Why?  If not, for every  $\gamma  < \theta $  (by trying  $\gamma (\ast ) = 
\gamma )$  there is a club  $E_\gamma $ of  $\lambda $  exemplifying the 
failure of $(\ast )_1$ for  $\gamma $.  Let  
$E = \dbca_{\gamma < \theta} E_\gamma \cap E^\ast_1$, so $E$ is a club
of $\lambda$,  hence    
$$
S' =: \{ \delta:\delta < \lambda,\delta \in S_1 ( \text{so cf}(\delta) 
= \sigma) \text{ and } C_\delta \subseteq E\}
$$
is a stationary subset of  $\lambda $.  For each  $\delta  \in  S'$  and  $\epsilon 
< \sigma $  for some  $\gamma  = \gamma (\delta ,\epsilon ) < \theta $  we have
$\beta^\ast _{\delta ,\epsilon } \in  
A_{\beta_{\delta ,\epsilon +1},\gamma }$,  but as $\sigma = \text{ cf}(\sigma) \ne  
\text{ cf}(\theta) = \theta$ for some $\gamma(\delta)$,  $\{\epsilon  < 
\sigma :\epsilon \gamma (\delta ,\epsilon ) = \gamma(\delta )\}$  is unbounded
in  $\sigma $.  But  $\delta  \in  E_{\gamma (\delta)}$,
contradiction.]
\sn
\item "{$(*)_2$}"   Without loss of generality: if  $\beta  \in
\text{ nacc}(C_\alpha),\alpha < \lambda$ then $(\exists \xi  \in  
A_{\beta ,\gamma (\ast )})[\beta  > \xi  > \sup (\beta  \cap  C_\alpha ) \and 
\beta  \cap  C_\alpha  = C_\xi]$. \nl   
[Why?  Define  $C'_\alpha $  for  $\alpha  < \lambda$: \nl    
$C^0_\alpha  = \{ \beta :\beta  \in \text{ nacc}(C_\alpha) \text{ and }
(\exists \xi  \in  A_{\beta ,\gamma(\ast)})[\beta > \xi \ge \sup(\beta 
\cap C_\alpha) \and \beta \cap C_\alpha = C_\xi]\}$. \nl   
$C'_\alpha$ is:  $\emptyset$ if $\alpha \in S_2$,  
$\alpha > \sup(C^0_\alpha)$ \nl
$\alpha \cap \text{ closure of } C^0_\alpha$ otherwise.]    
Now  $\langle C_\alpha :\alpha  < \lambda \rangle $  can be replaced by  
$\langle C'_\alpha :\alpha < \lambda \rangle$.]
\sn
\item "{$(*)_3$}"   For some 
$\gamma_1 = \gamma_1(\ast) < \theta$ for every club 
$E$ of $\lambda$ for some $\delta \in E:\text{cf}(\delta ) = 
\text{ cf}(\delta (\ast ))$,  and there is a club  $e$ of 
$\delta $  satisfying:  $e \subseteq E$,  otp$(e)$  is  
$\delta (\ast )$,  and for arbitrarily large  $\beta  \in  
\text{ nacc}(e)$  we have  $e \cap  \beta  \in  \{C_\zeta :\zeta  \in  
A_{\delta ,\gamma_1}\}$. \nl
[Why?  If not, for each  $\gamma _1 < \theta $  there is a club  
$E_{\gamma _1}$ of  $\lambda $  for which there is no  $\delta $  as required. 
Let  $E =: \dbca_{\gamma_1 < \theta} E_{\gamma_1}$, so $E$ is a club
of  $\lambda $  hence for some  
$\alpha  \in \text{ acc}(E) \cap  S_2$,  $C_\alpha \subseteq E$.  
Letting again $C_\alpha = \{\beta_{\alpha,\epsilon}:\epsilon  < 
\sigma \}$  (increasing),  $C_\alpha  \cap  \beta _{\alpha ,\epsilon } = 
C_{\delta ,\beta ^\ast _{\delta ,\epsilon }}$ where  
$\beta ^\ast _{\delta ,\epsilon } \in  
A_{\beta _{\delta ,\epsilon +1},\gamma (\ast )}$ clearly  $\delta  =: 
\beta _{\alpha ,\delta (\ast )}$,  $e = \{\beta _{\delta ,\epsilon }:\epsilon  
< \delta (\ast )\}$  satisfies the requirements except the last.  
As cf$(\delta(\ast)) \ne \text{ cf}(\mu )$,  for some   
$\gamma _1(\ast ) < \theta $,  $\gamma _1(\ast ) \ge \gamma (\ast )$  and  
$\{\epsilon < \delta (\ast ):\beta ^\ast _{\delta ,\epsilon } \in  
A_{\beta _{\delta ,\delta (\ast )},\gamma_1(\ast)}\}$  is 
unbounded in  $\delta (\ast )$.  Clearly  $\delta  =: 
\beta _{\alpha ,\delta (\ast )}$,  $e =: C_\alpha  \cap  \delta $  satisfies 
the requirement.  Now this contradicts the choice of  $E_{\gamma
_1(\ast )}$.]
\sn
\item "{$(\ast)_4$}"   For some club $E^a$ of $\lambda$,  for every club  $E^b 
\subseteq E^a$ of $\lambda $,  for some  $\delta  \in  E^b$ we have:
{\roster
\itemitem{ $(a)$ }   cf$(\delta ) = \text{ cf}(\delta(\ast))$   
\sn
\itemitem{ $(b)$ }   for some club $e$ of $\delta:e \subseteq E^b$,
otp$(e) = \delta (\ast )$,  and for arbitrarily large  $\beta  \in  
\text{ nacc}(e)$  we have  $e \cap  \beta  \in  \{C_\xi :\epsilon  \in  
A_{\delta ,\gamma _1(\ast )}\}$   
\sn
\itemitem{ $(c)$ }   for every $\beta \in A_{\delta,\gamma_1(\ast)}$
we have:  $C_\beta \subseteq E^a \Rightarrow  C_\beta \subseteq
E^b$ (we could have demanded $C_\beta \cap  E^a = C_\beta  \cap
E^b)$. \nl     
[Why?  If not we choose $E_i$ for $i < \mu^+_{\gamma_1(\ast)}$ by 
induction on  $i$,  $[j < i \Rightarrow  E_i \subseteq E_j]$,  
$E_i$ a club of  $\lambda $,  and  $E_{i+1}$ exemplify the failure of  
$E_i$ as a candidate for  $E^a$.  So  $\dbca_i E_i$ is a club of 
$\lambda $  hence by $(\ast )_3$ there are  $\delta $  and  
$e$  as there.  Now  $\langle \{\beta  \in  
A_{\delta ,\gamma _1(\ast )}:C_\beta \subseteq E_i\}:
i < \mu ^+_{\gamma _1(\ast )}\rangle $  is a decreasing sequence of 
subsets of  $A_{\delta ,\gamma _1(\ast )}$ of length  
$\mu ^+_{\gamma _1(\ast )}$,  and  $|A_{\delta ,\gamma _1(\ast )}| \le  
\mu _{\gamma _1(\ast )}$,  hence it is eventually constant.  So for every  $i$ 
large enough,  $\delta $  contradicts the choice of  $E_{i+1}$.]
\endroster}
\endroster
\enddemo
\bn
\centerline {$* \qquad * \qquad *$}
\bn
Let $S = \{ \delta  < \lambda:\text{cf}(\delta ) = \text{
cf}(\delta(\ast))$, and there is a club  $e = e_\delta $        
of  $\delta $  satisfying:  $e \subseteq E^a$,  otp$(e) = \delta
(\ast)$,  $\alpha  \in \text{ nacc}(e) \Rightarrow  e \cap  
\alpha  \in  A_{\alpha ,\gamma (\ast )}$ and for arbitrarily large  $\beta  
\in \text{ nacc}(e)$  we have  $e \cap  \beta  \in  \{C_\xi :\xi  \in  
A_{\delta ,\gamma (\ast )}\}\}$.
\mn
So  $S$  is stationary, let for  $\delta  \in  S$,  $C^\ast_\delta$ be an  
$e$  as above.  For  $\alpha  < \lambda $ let  ${\Cal P}_{1,\alpha } = 
\{C_\beta:\beta \le \alpha,\beta \in A_{\alpha ,\gamma _2(\ast )}\}$
\mr 
\item "{$(\ast)_5(a)$}"   for every club $E$ of $\lambda$, for some $\delta  \in  
S$,  $C^\ast _\delta  \subseteq E$   
\sn
\item "{$(b)$}"   $C^\ast_\delta$ is a club of $\delta$, otp$(C^\ast_\delta) = 
\delta(\ast)$   
\sn
\item "{$(c)$}"   if $\beta \in \text{ nacc } C^\ast_\delta (\delta  \in  S)$  then  
$C^\ast _\delta  \cap  \beta  \in  {\Cal P}_{1,\beta }$   
\sn
\item "{$(d)$}"   $|{\Cal P}_{1,\beta }| \le \mu_{\gamma (\ast )}$,  
${\Cal P}_{1,\beta}$ is a family of closed subsets of  $\beta $  of order type
$< \delta (\ast)$, \nl
[Why?  This is what we have proved in $(\ast )_4$; noting that in $(\ast )_4$ 
in $(b)$, $(e)$ is not uniquely determined, but by $(c)$ every ``reasonable" 
candidate is O.K.]
\ermn
Now repeating $(\ast)_3$, $(\ast)_4$ of the proof of \scite{1.10}, and we finish. 
\hfill$\square_{\scite{1.12}}$
\bigskip

\proclaim{\stag{1.13} Claim}  1) Assume 
$\lambda = \mu^+$, $|\delta(\ast)| < \mu,\aleph_0 < { \text{\rm cf\/}}
(\delta(*)) = { \text{\rm cf\/}}(\mu)(< \mu)$; \ub{then} we can find 
stationary  $S \subseteq \{\delta < \lambda:{\text{\rm cf\/}}(\delta) =
{ \text{\rm cf\/}}(\delta(*))\}$ and $\bar {\Cal P}$ such that 
${}^* \otimes^{\lambda,\delta(*)}_{\bar{\Cal P},S}$, except when: 
\mr
\item "{$\oplus$}"   for every regular $\sigma < \mu$, we can find  $h:\sigma  
\rightarrow { \text{\rm cf\/}}(\mu)$ such that for no
$\delta,\epsilon$ do we have: 
if $\delta < \sigma,{\text{\rm cf\/}}(\delta) = { \text{\rm
cf\/}}(\mu),\epsilon < { \text{\rm cf\/}}(\mu)$ \ub{then} 
$\{\alpha < \delta:h(\alpha) < \epsilon \}$ 
is not a stationary subset of $\delta$.
\ermn
2)  In \scite{1.12} and \scite{1.13}(1) we can have $\mu > 
\sup\{|{\Cal P}_\alpha|:\alpha < \lambda\}$. \nl
3)  If \scite{1.13}(2) if $\mu$ is strong limit we can have  
$|{\Cal P}_\alpha| \le 1$ for each $\alpha$.
\endproclaim
\bigskip

\remark{Remark}   Compare with \cite[\S3]{Sh:186}.
\endremark
\bigskip

\demo{Proof}  Left to the 
reader (reread the proof of \scite{1.12} and \cite[\S3]{Sh:186}.
\enddemo
\bigskip

\proclaim{\stag{1.14} Claim}  1) Let 
$\kappa $  be regular uncountable and we have global 
choice (or restrict ourselves to  $\lambda < \lambda ^\ast )$.  We can choose 
for each regular  $\lambda  > \kappa ^+$,  $\bar {\Cal P}^\lambda  = 
\langle {\Cal P}^\lambda _\alpha :\alpha  < \lambda \rangle $  (assuming global
choice) such that:
\mr
\item "{$(a)$}"   for each $\lambda$, ${\Cal P}^\lambda_\alpha$ is a
family of $\le \lambda$ of closed subsets of  $\alpha $  of order type  $< \kappa$.
\sn
\item "{$(b)$}"   if  $\chi $  is regular,  $F$ 
is the function  $\lambda  \mapsto  
\bar {\Cal P}^\lambda$ (for  $\lambda$ regular $< \chi$),  $\aleph _0 < 
\kappa = { \text{\rm cf\/}}(\kappa),\kappa^{++} < \chi,x \in {\Cal H}(\chi)$ 
then we can find  $\bar N = \langle N_i:i \le  \kappa \rangle$,  an 
increasing continuous chain of elementary submodels of  
$({\Cal H}(\chi),\in,<^*_\chi,F),\langle N_j:j \le i \rangle \in
N_{i+1}$,  $\|N_i\| = \aleph _0 + |i|$,  $x \in  N_0$ such that: 
{\roster
\itemitem{ $(*)$ }  if $\kappa^+ < \theta = { \text{\rm cf\/}}(\theta) \in
N_i$, then for some club $C$ of $\sup(N_\kappa \cap \theta)$ of
order type $\kappa$;  for any $j^i_1 < j < \kappa$ we have: \nl
$C \cap \sup (N_j \cap \theta) \in N_{j+1},{\text{\rm otp\/}}
(C \cap  \sup (N_j \cap  \theta)) = j$.
\endroster}
\ermn
2)  We can above have $|{\Cal P}^\lambda_\alpha| < \lambda$.
\endproclaim
\bigskip

\demo{Proof}   1) Let 
$\langle C_\alpha :\alpha  \in  S\rangle $  be such that  $S 
\subseteq  \{\alpha  \le \kappa^{++}:\text{ cf}(\alpha) \le  \kappa\}$  is 
stationary, otp$(C_\alpha) \le \kappa$,  $[\beta  \in  C_\alpha  
\Rightarrow  C_\beta  = \beta  \cap C_\alpha],C_\alpha$ a closed subset 
of  $\alpha $,  $[\alpha$ limit $\Rightarrow \alpha = \sup(C_\alpha)]$,
$\{\alpha \in S:\text{cf}(\alpha) = \kappa\}$ 
stationary, and for every club  $E$  
of  $\kappa ^{++}$ there is  $\delta  \in  S$,  cf$(\delta ) = \kappa $,  
$C_\delta  \subseteq  E$.  For  $i \in  \kappa^{++}\backslash S$  let
$C_i = \emptyset $.  Now for every regular  $\lambda > \kappa^+$ 
and $\alpha  \le \lambda $,  let  $e^\lambda _\alpha  
\subseteq \alpha$ be a club of $\alpha$  of order type cf$(\alpha)$.  
For  $\lambda$ as above and for $\alpha \le \lambda$ limit let 
$\bar{\Cal P}^\lambda_\alpha  = \{\{i \in e_\delta :i <
\alpha,\text{otp}(e_\delta \cap i) \in  C_\beta \}:\delta  < \lambda$  has 
cofinality  $\kappa^{++}$, and $\beta \in  S\}$.
Given $x \in  H(\chi )$,  we choose by induction on  $i < \kappa ^{++}$,  
$M_i$,  $N_i$ such that: \nl   

$N_i \prec  M_i \prec ({\Cal H}(\chi ),\in,<^*_\chi,F)$  

$\|M_i\| = |i| + \aleph _0$  

$\|N_i\| = |C_i| + \aleph _0$  

$M_i(i < \kappa ^{++})$  is increasing continuous  

$x \in  M_0,$  

$\langle M_j:j \le i \rangle  \in  M_{i+1}$  

$N_i$ is the Skolem Hull of $\{\langle N_j:j \in C_\zeta \rangle:
\zeta \in C_i\}$.
\sn
We leave the checking to the reader. \nl
2)  We imitate the proof of \scite{1.4}. \hfill$\square_{\scite{1.14}}$
\enddemo
\newpage

\head {\S2 Measuring $[\lambda]^{< \kappa}$} \endhead  \resetall \sectno=2
 \spuriousreset
\bigskip

We prove here that two natural ways to measure  ${\Cal S}_{<\kappa }(\lambda)$ 
for $\kappa$ regular uncountable, give the same cardinal: the minimal 
cardinality of a cofinal subset; i.e. its cofinality
(i.e. cov$(\lambda,\kappa,\kappa,2))$ 
and the minimal cardinality of a stationary
subset.  The theorem is really somewhat stronger: for appropriate normal ideal 
on  ${\Cal S}_{<\kappa }(\lambda )$,  some member of the dual filter has the 
right cardinality. 

The problem is natural and I did not trace its origin, but until recent years 
it seems (at least to me) it surely is independent, and find it gratifying we 
get a clean answer.  I thank P. Matet and M. Gitik of reminding me of the 
problem. 

We then find applications to $\Delta$-systems and largeness of 
$\check I[\lambda]$.
\bigskip

\definition{\stag{2.1} Definition}  1) Let $(\bar C,\bar{\Cal P},Z) \in  
{\Cal T}^\ast[\theta,\kappa]$ \ub{when}: 
\mr
\widestnumber\item{$(viii)^-$}
\item "{$(i)$}"   $\aleph_0 < \kappa = \text{ cf}(\kappa) < \theta 
= \text{ cf}(\theta)$,
\sn
\item "{$(ii)$}"   $\bar C = \langle C_\delta :\delta \in S\rangle,
\bar{\Cal P} = \langle {\Cal P}_\delta:\delta \in S \rangle,Z =
\langle <_{{\Cal P}_\delta}:\delta \in S\rangle$ 
\sn
\item "{$(iii)$}"  $S \subseteq \theta$, $S$ is stationary 
(we shall write  $S = S(\bar C))$, 
\sn
\item "{$(iv)$}"  $C_\delta$ is an unbounded subset of $\delta$, (not necessarily 
closed) 
\sn
\item "{$(v)$}"   id$^a(\bar C)$ is a proper ideal (i.e. for every club  $E$  of  
$\theta $  for some  $\delta  \in  S$,  $C_\delta  \subseteq  E)$ 
\sn
\item "{$(vi)$}"    $\dsize \bigwedge_{\delta \in S} 
\text{ otp}(C_\delta) < \kappa$,  (hence $[\delta \in S \Rightarrow \text{ cf}(\delta)
< \kappa])$ 
\sn
\item "{$(vii)$}"  $(\alpha) \quad 
{\Cal P}_\delta$ is a family of bounded subsets of $C_\delta$, directed \nl

$\quad$ by the partial order $<_{{\Cal P}_\delta}$ which is a partial order
on
\nl

$\quad {\Cal P}^* = \{x \cap \alpha:x \in {\Cal P}_\delta$ for some $\delta
\in S$ and $\alpha < \theta\}$ satisfying 
\nl

$\quad y <_{{\Cal P}_\delta} z \Rightarrow y \subseteq z$, 
(but see parts (1A),(1B))
\sn
\item "{${{}}$}"  $(\beta) \quad
\dbcu_{x \in {\Cal P}_\delta }x = C_\delta$, and  
$|{\Cal P}_\delta | < \kappa$ 
\sn
\item "{$(viii)$}"   for some\footnote{a sufficient condition is:
\mr
\item "{$(viii)^+$}"  for every $\alpha < \theta$ the set   
${\Cal P}^\ast_\alpha  =: \{ a \cap  \alpha:\text{ for some } \delta
\in S$  we have  $\alpha < \delta  \in  S$,  $a \in {\Cal P}_\delta$ 
and $\alpha \in C_\delta \}$ has cardinality  $< \theta$ or at least
\endroster} 
list $\langle b^*_i:i < \theta \rangle$
of $\dbcu_{\alpha  \in S} {\Cal P}_\alpha \cup \{\emptyset\}$
satisfying $b^*_i \subseteq i$ we have: for every $\alpha \in S$ we have
${\Cal P}_\alpha \subseteq \{b^*_j:j < \alpha\}$
\sn
\item "{$(ix)$}"   for  $x \in \dbcu_{\delta \in S} {\Cal P}_\delta$
we have the set ${\Cal P}_x :=
\{y \in \dbcu_{\delta \in S} {\Cal P}_\delta:y <_{{\Cal P}_\delta}
x\}$ has cardinality $\kappa$.
\ermn
1A) If each $<_{{\Cal P}_\delta}$ is inclusion we may omit it.
\nl
1B) If $<_*$ is a partial order of $\dbcu_{\delta \in S} {\Cal
P}_\delta$ and $\delta \in S \Rightarrow <_{{\Cal P}_\delta} = <_*
\restriction {\Cal P}_\delta$ then we may write $<_*$ instead of $Z$.
\nl
2)  $\bar C \in {\Cal T}^0[\theta,\kappa]$, if $(\bar C,\bar{\Cal P}) 
\in {\Cal T}^\ast [\theta ,\kappa]$ where $\delta \in S(\bar C)
\Rightarrow {\Cal P}_\delta  = \{C_\delta \cap \alpha:\alpha \in C_\delta\}$. 
\nl
3)  $\bar C \in {\Cal T}^1[\theta,\kappa]$ if $(\bar C,\bar{\Cal P}) \in 
{\Cal T}^\ast [\theta ,\kappa]$ where $\delta \in S(\bar C) \Rightarrow  
{\Cal P}_\delta = [C_\delta]^{<\aleph_0}$.
\enddefinition
\bn
Note that:
\proclaim{\stag{2.2} Claim}   1) If $\theta = { \text{\rm cf\/}}(\theta)
 > \kappa = { \text{\rm cf\/}}(\kappa) > \sigma = { \text{\rm cf\/}}
(\sigma)$, \ub{then} there is  
$\bar C \in {\Cal T}^1[\theta,\kappa]$ such that:

$$
\{\delta \in S(\bar C):{\text{\rm cf\/}}(\delta) = \sigma\} \ne
\emptyset \, { \text{\rm mod id\/}}^a(\bar C).
$$
\medskip
\noindent
2)  If $S \subseteq \{\delta < \theta:{\text{\rm cf\/}}(\delta) < \kappa\}$ is 
stationary, $\bar C$ an $S$-club system, $|C_\delta| < \kappa$, and  
${\text{\rm id\/}}^a(\bar C)$ a proper ideal, \ub{then} 
$\bar C \in {\Cal T}^1[\theta,\kappa]$.
\nl
3) In (2) if in addition for each $\alpha < \theta$ we have 
$|\{C_\delta \cap \alpha:\alpha \in
C_\delta,\delta \in S\}| < \theta$ \ub{then} 
$\bar C \in {\Cal T}^0[\theta,\kappa]$.
\nl
4) If $\theta$ is a successor of regular \ub{then} in part
(2) we can demand  
$\bar C \in {\Cal T}^0[\theta,\kappa]$ each $C_\delta$ closed.
\nl
5)  If $\theta = { \text{\rm cf\/}}(\theta) > \kappa = { \text{\rm
cf\/}}(\kappa) > \sigma = { \text{\rm cf\/}}(\sigma)$, \ub{then} there is  
$\bar C \in {\Cal T}^0[\theta,\kappa]$ such that:  
$\{\delta \in S(\bar C):{\text{\rm cf\/}}(\delta) 
= \sigma\} \ne \emptyset \,{ \text{\rm mod id\/}}^a(\bar C)$.
\nl
6) If $\theta = \text{\rm cf}(\theta) = \text{ cf}(\kappa) > \sigma 
= \text{\rm cf}(\sigma)$ and $S \in \check I[\theta]$ is stationary \ub{then}
   there is $\bar C \in {\Cal T}^0[\theta,\kappa]$ such that $S(\bar
   C) = S$.
\endproclaim
\bigskip

\demo{Proof}  1) Let $S_0 \subseteq \{\delta < \theta:
\text{cf}(\delta) = \sigma\}$  
be stationary, $C^0_\delta$ a club of $\delta$ of order type $\sigma$.  
By \cite[\S2]{Sh:365}, for some club $E$ of $\lambda$
letting $S = S_0 \cap \text{ acc}(E)$ and letting, for $\delta \in S,
C_\delta  = g \ell(C^0_\delta,E) 
= \{\sup (\alpha \cap E):\alpha \in C_\delta\}$  
we have $S \notin \text{ id}^a(\langle C_\delta:
\delta \in S_0 \rangle)$, now use part (2). 
\nl
2)  Check.
\nl
3)  Check.
\nl
4)  By \cite[\S4]{Sh:351}, \cite[Ch.IV,3.4]{Sh:e}(2) 
or \cite[2.14]{Sh:365}(2)((c)+(d)) but see \cite{Sh:E12}.
\nl
5)  By \scite{1.5} and \scite{1.11} (so we use the non-accumulation
points). 
\nl
6) Similarly.   \hfill$\square_{\scite{2.2}}$
\sn
Remember (see \cite[\S3]{Sh:52}).
\enddemo
\bigskip

\definition{\stag{2.3} Definition}  1)  ${\Cal D}^\kappa_{<\kappa}
(\lambda)$ is the filter on $[\lambda]^{< \kappa}$ defined by:

for $X \subseteq [\lambda]^{<\kappa}$: \nl

\block $X \in {\Cal D}^\kappa_{<\kappa}(\lambda)$ \ub{iff} 
there is a function $F$ with domain the set of sequences of length $<
\kappa$ with elements from $[\lambda]^{<
\kappa}$ and $F$ is into $[\lambda]^{< \kappa}$
such that: if $a_\zeta \in [\lambda]^{<\kappa}$ for $\zeta < \kappa$,
is $\subseteq$-increasing continuous and for each  $\zeta < \kappa$  we have  
$F (\langle \ldots,a_\xi,\ldots \rangle)_{\xi \le \zeta} \subseteq  
a_{\zeta +1}$ then $\{\zeta < \kappa:a_\zeta \in X\} \in {\Cal
D}_\kappa$.
\endblock
\sn
(recall that ${\Cal D}_\kappa$ the filter generated by the family of clubs of
$\kappa$). 
\enddefinition
\bigskip
\noindent
Similarly
\definition{\stag{2.4} Definition}  For $\lambda \ge \theta = \text{
cf}(\theta) > \kappa  = \text{ cf}(\kappa) > \aleph_0$, $(\bar
C,\bar{\Cal P}) \in {\Cal T}^*[\theta,\kappa]$ and set $X$ of
cardinality $\ge \kappa$ we 
define a filter ${\Cal D}_{(\bar C,\bar{\Cal P})}(\lambda)$ on  
$[\lambda]^{<\kappa}$; (letting, e.g. $\chi = \beth_{\omega +1}(\lambda))$:
\sn 
$Y \in {\Cal D}_{(\bar C,\bar{\Cal P})}(\lambda)$ \ub{iff} $Y \subseteq 
[\lambda]^{<\kappa}$ and for some $\bold x \in {\Cal H}(\chi)$,  for every    
$\langle N_\alpha,N^\ast_a:\alpha < \theta,a \in \dbcu_{\delta \in S}
{\Cal P}_\delta \rangle$ satisfying $\otimes$ below, also 
there is $A \in \text{ id}^a(\bar C)$  such that:  
$\delta \in S(\bar C)\backslash A \Rightarrow \dbcu_{a \in {\Cal
P}_\delta}N^\ast_a \cap \lambda \in Y$ where, letting ${\Cal P} =
\cup\{{\Cal P}_\delta:\delta \in S\}$,
\mr
\item "{$\otimes(i)$}"  $N_\alpha \prec ({\Cal H}(\chi),\in,<^\ast_\chi)$   
\sn
\item "{$(ii)$}"  $\|N_\alpha\| < \theta$,
\sn
\item "{$(iii)$}"   $\langle N_\beta:
\beta \le \alpha \rangle \in N_{\alpha +1}$   
\sn
\item "{$(iv)$}"   $\langle N_\alpha:\alpha < \theta\rangle$ 
is increasing continuous   
\sn
\item "{$(v)$}"   $N^\ast_a \prec ({\Cal H}(\chi),\in,<^\ast_\chi)$ for $a
\in \dbcu_{\delta \in S} {\Cal P}_\delta $   
\sn
\item "{$(vi)$}"   $\|N^\ast_a\| < \kappa$, $N^\ast_a \cap \kappa$ an initial 
segment of $\kappa$
\sn
\item "{$(vii)$}"    $b \subseteq a$  (both in $\dbcu_{\delta \in S}
{\Cal P}_\delta)$  implies  $N^\ast _b \prec  N^\ast _a$   
\sn
\item "{$(viii)$}"  if  $\alpha \in a \in \dbcu_{\delta \in S} {\Cal
P}_\delta$ then  $\langle N_\beta,N^\ast_b:\beta \le \alpha,b
\subseteq a,b \in \{b^*_i:i \le \alpha\} \subseteq
{\Cal P} = \dbcu_{\delta \in S} {\Cal P}_\delta \rangle$
belongs to $N^\ast_a$   
\sn
\item "{$(ix)$}"   $\langle N_\beta,N^\ast_b:
\beta \le \alpha,b \subseteq \alpha + 1,b \in \{b^*_i:i \le \alpha
+1\} \subseteq \dbcu_{\delta \in S} {\Cal P}_\delta 
\rangle$ belongs to $N_{\alpha +1}$   
\sn
\item "{$(x)$}"   $a \subseteq N^\ast _a$ and $\alpha \in a
\Rightarrow \alpha \cap  a \in N^\ast_a$ 
\sn
\item "{$(xi)$}"    $a \subseteq \alpha,a \in {\Cal P}$ 
implies $N^\ast_a \in N_{\alpha +1}$ (follows from (ix) by clause
(viii) of Definition \scite{2.1}(1))
\sn
\item "{$(xii)$}"  $a \in {\Cal P}_\delta \and 
\delta \in S \and \alpha < \theta \Rightarrow \bold x \in  
N^\ast_a \and \bold x \in N_\alpha$.
\endroster
\enddefinition
\bn
Clearly
\proclaim{\stag{2.5} Claim}   1) Any $\chi$ that is ${\Cal H}(\chi)$ can serve,
and $\bold x = (Y,\lambda,\bar C,\bar{\Cal P})$  is enough. 
\nl
2) ${\Cal D}_{(\bar C,\bar{\Cal P})}(\lambda)$ is a (non-trivial)
fine $(< \kappa)$-complete filter on  
$[\lambda]^{<\kappa}$ when $(\bar C,\bar{\Cal P}) \in  
{\Cal T}^\ast [\theta,\kappa]$,  $\lambda  \ge \theta$,  hence it extends 
${\Cal D}_{<\kappa}(\lambda)$.  (Remember ${\text{\rm id\/}}^a(\bar C)$
is a proper ideal).
\endproclaim
\bigskip

\demo{Proof}  Should be clear.  \hfill$\square_{\scite{2.5}}$
\enddemo
\bigskip

\proclaim{\stag{2.6} Theorem}   Suppose $\lambda > \theta = 
\text{\rm cf}(\theta) > \kappa = { \text{\rm cf\/}}(\kappa) > \aleph_0$ and
$\theta = \kappa^+$.  \ub{Then} the 
following four cardinals are equal for any $(\bar C,\bar{\Cal P}) \in  
{\Cal T}^\ast[\theta,\kappa]$, recalling there are such $(\bar
C,\bar{\Cal P})$ by \scite{2.2}:
\mn
$\mu(0) = { \text{\rm cf\/}}([\lambda]^{< \kappa},\subseteq)$ 
\mn
$\mu(1) = { \text{\rm cov\/}}(\lambda,\kappa,\kappa,2) = { \text{\rm
Min\/}}\{|{\Cal P}|:{\Cal P} \subseteq [\lambda]^{< \kappa}$, 
and for every  $a \subseteq \lambda,|a| < \kappa$
there is $b \in {\Cal P}$ satisfying $a \subseteq b\}$
\mn
$\mu(2) = { \text{\rm Min\/}}\{|S|:S \subseteq [\lambda]^{<\kappa}$ is
stationary$\}$
\mn
$\mu(3) = \mu_{(\bar C,\bar{\Cal P})} = \text{\rm Min}\{|Y|:Y \in  
{\Cal D}_{(\bar C,\bar P)}(\lambda)\}$.
\endproclaim
\bigskip

\remark{\stag{2.6A} Remark}  0) We thank M. Shioya for asking for a
correction of an inaccuracy in the proof in a meeting in 
the summer of 1999 in which we answer him; this and 
other minor changes are done here.  I thank P. Komjath for helpful
comments and S. Garti for help in proofreading.
\nl
1) It is well known that if $\lambda > 2^{<\kappa}$ then 
the equality holds as they are all equal to $\lambda^{< \kappa}$. \nl
2)  This is close to ``strong covering".
\nl
3) Note that only $\mu(3)$ has $(\bar C,\bar{\Cal P})$  in its 
definition, so actually  $\mu (3)$  does not depend on  
$(\bar C,\bar{\Cal P})$, recalling that by Claim \scite{2.2} we know
that ${\Cal T}^*[\theta,\kappa]$ is not empty. 
\nl
4) $\mu(0),\mu(1)$ are equal trivially.
\endremark
\bigskip

\remark{\stag{2.6B} Remark}   0) We can concentrate on the case $(\bar
C,\bar{\Cal P}) \in {\Cal T}^1[\theta,\kappa]$ or ${\Cal
T}^0[\theta,\kappa)$.  This somewhat simplifies and is enough.
\nl
1) We can weaken in Definition
\scite{2.1}(1) demand (ix) as follows: 
\mr
\item "{$(ix)'$}"   there is a sequence $\langle a_i,{\Cal P}^\ast_i:i
< \lambda \rangle$ such that   
{\roster
\itemitem{ $(a)$ }   $|a_i| < \kappa,{\Cal P}^\ast_i$ 
is a family of  $< \kappa $  subsets of  $a_i$   
\sn
\itemitem{ $(b)$ }   for every $\delta \in S$ and 
$x \in {\Cal P}_\delta$ for some  $i <  \delta,a_i = x$ and \nl
$(\forall b)[b \in {\Cal P}_\delta \and b \subseteq 
a \Rightarrow  b \in {\Cal P}^\ast_i]$.
\endroster}
\ermn
In this case \scite{2.6}, \scite{2.6A}(4) (and \scite{2.5}) 
remains true and we can strengthen \scite{2.2}. 
\nl
2)  We can even use ${\Cal P}_\delta$ with another order (not $\subseteq$).
\endremark
\bigskip

\demo{Proof}   Clearly $\lambda \le \mu(0) = \mu(1) \le \mu (2) \le \mu(3)$
(the last --- by \scite{2.5}(2)).  So we shall finish by proving
$\mu(3) \le \mu (1)$,  
and let ${\Cal Q}$ exemplify $\mu(1) =
\text{ cov}(\lambda ,\kappa ,\kappa ,2)$. Let $S=S(\bar C)$, etc.

Let $\chi$ be e.g. $\beth_3(\lambda)^+$ and let $M^\ast_\lambda$
be the model with universe  $\lambda  + 1$
and all functions definable in  
$({\Cal H}(\chi),\in,<^\ast_\chi,\lambda,\kappa,\mu(1))$.  Let $M^\ast$ be an
elementary submodel of $({\Cal H}(\chi),\in,<^\ast_\chi)$  of cardinality  
$\mu(1)$ such that ${\Cal Q} \in M^*,M^*_\lambda \in  M^\ast$,  
$(\bar C,\bar{\Cal P}) \in M^\ast$ and $\mu(1) + 1 \subseteq M^\ast$ 
hence ${\Cal Q} \subseteq M^\ast$.  It is enough to prove 
that $M^\ast  \cap  [\lambda]^{< \kappa}$ 
belongs to ${\Cal D}_{(\bar C,\bar P)}(\lambda)$. 

So let $N_i$ (for $i < \theta$), $N^*_x$ 
(for $x \in \dbcu_{\delta \in S} {\Cal P}_\delta)$ be such that: 
they satisfy  $\otimes$ of Definition \scite{2.4} for $\bold x :=
\langle M^*_\lambda,M^\ast,{\Cal P},\lambda,\kappa $, 
$(\bar C,\bar{\Cal P})\rangle$ so it belongs to 
every  $N_\alpha $,  $N^\ast_x$.  It is 
enough to prove that $\{\delta \in S:\lambda \cap \dbcu_{x \in
{\Cal P}_\delta} N^\ast_x \in M^\ast\} = \theta$ mod id$^a(\bar C)$.
For $i \in S$ clearly $x \subseteq y$ (or 
$x <_{{\Cal P}_i} y) \Rightarrow N^*_x \prec
N^*_y$ and ${\Cal P}_i$ is directed (by the partial order $\subseteq$
or $<_{{\Cal P}_i}$ recalling clause (vii) of $\otimes$ of Definition
\scite{2.4}) hence $N'_i := \cup\{N^*_x:x \in
{\Cal P}_i\}$ is $\prec ({\Cal H}(\chi),\in,<^*_\chi)$ and even $\prec
N_i$ and $N'_i$ has
cardinality $< \kappa$ (as $|{\Cal P}_i| < \kappa$ and each $N^*_x$
has cardinality $< \kappa$ and $\kappa$ is regular) 
and we have to show that $\{i \in S:\lambda
\cap N'_i \in M^*\} = \theta$ mod id$^a(\bar C)$.

For each $i \in S$ by the choice of ${\Cal Q}$,
there is a set  $a_i$ such that $N'_i \cap \lambda = 
(\dbcu_{y \in {\Cal P}_i} N^\ast_y)
\cap \lambda \subseteq a_i \in {\Cal P}$;  
so as ${\Cal P}$ and $\langle N^*_y:y \in {\Cal P}_i\rangle$ belong
to $N_{i+1}$, see clause (ix) of Definition \scite{2.4}
 without loss of generality  $a_i \in N_{i+1}$.  Let ${\frak a}_i 
=: \text{ Reg } \cap  a_i \cap \lambda^+ \backslash \theta^+$, so  
${\frak a}_i$ is a set of $< \kappa$  
regular cardinals $\ge \theta^+$ and ${\frak a}_i \in N_{i+1}$ too, 
so there is a generating sequence 
$\langle {\frak b}_\lambda[{\frak a}_i]:\lambda \in \text{ pcf}({\frak
a}_i) \rangle$ as in \cite[VII,2.6]{Sh:g} = \cite[2.6]{Sh:371}, 
without loss of generality it is definable from ${\frak a}_i$ (in
$({\Cal H}(\chi),\in,<^\ast_\chi)$ say the $<^*_\chi$-first such object).  
Also ${\frak a}_i \in {\Cal P} \subseteq M^\ast$ so ${\frak a}_i \in M^\ast$.
As ${\frak a}_i \in N_{i+1}$ we have   
$\langle {\frak b}_\lambda[{\frak a}_i]:\lambda \in 
\text{ pcf}({\frak a}_i) \rangle  \in  N_{i+1} \cap  
M^\ast $,  and also there is $\langle f^{{\frak
a}_i}_{\partial,\alpha}:\alpha < \partial,\partial 
\in \text{ pcf}({\frak a}_i) \rangle$ as in \cite[VIII,1.2]{Sh:g} = 
\cite[1.2]{Sh:371}, and again without 
loss of generality it belongs to  $N_{i+1} \cap M^\ast$.  
As max pcf$({\frak a}_i) \le 
\text{ cov}(\lambda,\kappa,\kappa,2) = \mu(1)$,  (first inequality by 
\cite[II,5.4]{Sh:g} = \cite[5.4]{Sh:355}) 
clearly each $f^{{\frak a}_i}_{\partial,\alpha} \in M^\ast$.   

Let
\mr
\item "{$\odot_1$}"   $h$ be the function with domain ${\frak a} :=
\dbcu_{i \in S} {\frak a}_i$ defined by $h(\sigma) =
\sup(\sigma \cap \dbcu_{i < \theta} N_i)$.
\ermn
So by \cite[VIII,2.3]{Sh:g}(1) = \cite[2.3]{Sh:371}(1) 
\mr
\item "{$\odot_2$}"  if $i \in S$ \ub{then} 
$h \restriction {\frak a}_i$ has the 
form Max$\{f^{{\frak a}_i}_{\partial_\ell,\alpha_\ell}:\ell < n\}$ for
some $n < \omega,\partial_\ell \in \text{ pcf}({\frak a}_\ell)$ and
$\alpha_\ell < \partial_\ell$ for $\ell < n$
\nl
hence
\sn
\item "{$\odot_3$}"   if $i \in S$ \ub{then} $h \rest {\frak a}_i$
belongs to $M^\ast$
\nl
and obviously (as $\sigma \in {\frak a}_i \wedge i < j_1 < j_2
\Rightarrow \sup(\sigma \cap N_{j_1}) < \sup(\sigma \cap N_{j_2})$
\sn
\item "{$\odot_4$}"   $\sigma \in \text{ Dom}(h) \Rightarrow 
\text{ cf}(h(\sigma)) = \theta$.  
\ermn
Let  $e$  be a definable function in
$({\Cal H}(\chi),\in,<^\ast_\chi,\lambda,\kappa)$ with Dom$(e) = 
\lambda + 1$ such that $e(\alpha) = e_\alpha $ is a club 
of $\alpha$ of order type  
cf$(\alpha)$,  enumerated as $\langle e_\alpha(\zeta):\zeta < \text{
cf}(\alpha) \rangle$.  Now for each $\sigma \in \dbcu_{i < \theta} 
{\frak a}_i$ let
\mr
\item "{$\odot_5$}"   $E_\sigma =: \{i < \theta:(\forall \zeta < \theta)
[e_{h(\sigma)}(\zeta) \in N_i \Leftrightarrow \zeta < i],i$ is a limit ordinal
and sup$(N_i \cap \sigma) = \sup \{e_{h(\sigma)}(\zeta):\zeta < i\}\}$. 
\ermn
Clearly $E_\sigma$ is a club of $\theta$,  hence (on $\langle b^*_j:j
< \theta\rangle$, see clause (viii)$^-$ of Definition \scite{2.1})
 
$$
\align
E = \{\delta < \theta:&\delta \text{ is a limit ordinal and }
\sigma \in \cup\{{\frak a}_i:i < \delta\} \subseteq \\
 &\text{ Reg } \cap \lambda^+ \backslash \theta^+ 
\Rightarrow \delta \in \text{ acc}(E_\sigma) \text{ and } 
N_\delta \cap \theta = \delta\}
\endalign
$$
\mn
is a club of $\theta$.
For each $\delta \in E \cap S$ such that $C_\delta \subseteq E$,  let  
$\delta^\ast := \text{ sup}(\kappa \cap N'_\delta)
 = \sup(\kappa \cap \dbcu_{y \in {\Cal P}_\delta}
N^\ast_y)$ so $\delta^\ast < \kappa$,  and we 
define by induction on $n$  models $M_{y,\delta,n}$ for every  
$y \in  {\Cal P}_\delta$, (really, they do not depend on $\delta$).  

First, $M_{y,\delta,0}$ is 
the Skolem Hull in $M^\ast_\lambda$ of $\{i:i \in y\}
\cup (N'_\delta \cap \kappa)$.  
\nl
Second, $M_{y,\delta,n+1}$ is the Skolem Hull in  
$M^\ast_\lambda$ of $M_{y,\delta,n} \cup \{e_{h(\sigma)}(\zeta):
\sigma \in ({\text{\rm Reg\/}} \cap \lambda^+ \backslash \theta^+) 
\cap M_{y,\delta,n}$ and $\zeta \in y\}$.  Now we note
\mr
\item "{$(*)_0$}"  if $y \in \{b^*_i:i < \zeta\},\zeta \in C_\delta$
and $\delta \in E$ \ub{then} $N^*_y \in N_\zeta$ hence $N^*_y \prec
N_\zeta$.
\ermn
[Why?  By clause (ix) of $\otimes$ of Definition \scite{2.4} we
 have $N^*_y \in N_\zeta$; as $\|N^*_y\| < \kappa < \theta$ and
 $N_\zeta \cap \theta \in \theta$ as $\zeta \in C_\zeta \subseteq E$
 we have $N^*_y \subseteq N_\zeta$ hence $N^*_y \prec N_\zeta$.]
\mr
\item "{$(*)_1$}"  if $\zeta \in E (\subseteq \theta)$ and $\sigma \in
\text{ Reg } \cap N_\zeta \cap \lambda^+ \backslash \theta$ then
$e_{h(\sigma)}(\zeta) = \sup(N_\zeta \cap \sigma)$.
\ermn
[Why?  By the choice of $E$.]
\mr
\item "{$(*)_2$}"  assume $\delta \in S$ satisfies $\delta \in E$,
moreover $C_\delta \subseteq E$; if $y \in {\Cal P}_\delta$ and $\sigma \in
N^*_y \cap \text{ Reg } \lambda^+ \backslash \theta^+$ \ub{then}
($h(\sigma)$ has cofinality $\theta$, the sequence $\langle
e_{h(\sigma)}(\zeta):\zeta < \theta\rangle$ is increasing continuous
with limit $h(\sigma)$ and):
{\roster
\itemitem{ $(i)$ }  if $y \in \{b^*_i:i < \zeta\}$ and 
$\zeta \in C_\delta$ then
$\sup(N_\zeta \cap \sigma) = e_{h(\sigma)}(\zeta)$
\sn
\itemitem{ $(ii)$ }  if $y \in \{b^*_i:i < \zeta\},
\zeta \in z \in {\Cal P}_\delta$ and
$y <_{{\Cal P}_\delta} z$  \ub{then} $y \in N_z,N^*_y \in N^*_z,N^*_y \prec
N^*_z$ and $e_{h(\sigma)}(\zeta) \in N^*_z$
\sn
\itemitem{ $(iii)$ }  $\{e_{h(\sigma)}(\zeta):\zeta \in C_\delta\}$ 
is a subset of $N'_\delta = \dbcu_{z \in {\Cal P}_\delta} N^*_z$
\sn
\itemitem{ $(iv)$ }  the set above is an unbounded subset of
$N'_\delta \cap \sigma$.
\endroster}
\ermn
[Why?  \ub{Clause $(i)$}:  So we assume $\zeta \in C_\delta$ and $y
\in \{b^*_i:i < \zeta\}$. 

By $(*)_0$ we have $N^*_y \prec N_\zeta$.
By the definition of $E_\sigma$ as $\sigma \in N^*_y \prec
N_\zeta \wedge \zeta \in E$ clearly $\zeta \in E_\sigma$ hence
$\sup(N_\zeta \cap \theta) = e_{h(\sigma)}(\zeta)$ by $(*)_1$.
\sn
\ub{Clause $(ii)$}:  So assume $y \in \{b^*_i:i < \zeta\},\zeta\in z$ and $y
<_{{\Cal P}_\delta} z$ (so $y,z \in {\Cal P}_\delta$) hence
${\Cal P}_{z,\zeta} = \{x \in \dbcu_{\alpha \in S} {\Cal P}_\alpha:x
\subseteq z \cap \zeta\}$ has cardinality $< \kappa$ and $z \cap \zeta
\in N^*_z$ by clause (x) of \scite{2.4}, so ${\Cal P}_{z,\zeta} = \{x
\in \cup\{{\Cal P}_\alpha:\alpha < \delta\}:x \subseteq z \cap \zeta\}
\in N^*_z$, so (as $N^*_z \cap \kappa \in \kappa$,
$|{\Cal P}_{z,\zeta}| < \kappa$) clearly ${\Cal P}_{z,\zeta} 
\subseteq N^*_z$ hence $y \in N^*_z$.  By 
clause (viii) of $\otimes$ of Definition
\scite{2.4} it follows that $N^*_y \in N^*_z$.  But $|N^*_y| < \kappa
\wedge N^*_z \cap \kappa \in \kappa$ hence $N^*_y \subseteq N^*_z$ so
$N^*_y \prec N^*_z$.  But $\sigma \in N^*_y$
hence $\sigma \in N^*_z$.  Also $N_\zeta \in N^*_z$ as $\zeta \in z
\subseteq N^*_z$ recalling $(viii)$ of \scite{2.4} hence 
$e_{h(\sigma)}(\zeta) = \sup(N_\zeta \cap \sigma) \in N^*_z$ recalling
$(*)_1$ so we have shown all clauses of $(ii)$.
\sn
\ub{Clause $(iii)$}:  So let $\zeta \in C_\delta$; by clause
(vii)$(\beta)$ of Definition \scite{2.1} we know that $C_\delta =
\cup\{y:y \in {\Cal P}_\delta\}$ hence for some $y_1 \in {\Cal
P}_\delta$ we have $\zeta \in y_1$.  By clause $(x)$ of $\otimes$ from
Definition \scite{2.4} we have $y_1 \subseteq N^*_{y_1}$ hence $\zeta
\in N^*_{y_1}$.  Also we are assuming in $(*)_2$ that $\sigma \in
N^*_y,y \in {\Cal P}_\delta$, so recalling ${\Cal P}_\delta$ is
directed, we can find $y_2 \in {\Cal P}_\delta$ which is a common
$\subseteq$-upper bound of $y_1,y_2$ hence $N^*_y \prec N^*_{y_2},
N^*_{y_1} \prec N^*_{y_2}$ hence $\sigma,\zeta \in N^*_{y_2}$.  

By the choice of the function $e$ and the model $M^*_\lambda$ clearly
$e(-,-)$ is a function of $M^*_\lambda$, but the object $\bold x$
belongs to $N^*_{y_2}$ and by its choice this implies that $e \in
N^*_{y_2}$.  By clause (viii) of \scite{2.4} recalling $\zeta \in N^*_{y_2}$
we know that $N_\zeta \in N^*_{y_2}$ but $\sigma
\in N^*_{y_2}$ hence by $(*)_1$ we have sup$(N_\zeta \cap \sigma) \in
N^*_{y_1}$.   But we are assuming in $(*)_2$ that $C_\delta \subseteq
E$ and, see above, $\zeta \in C_\delta$ so $\zeta \in E$ and $\zeta \in
C_\delta \subseteq N_\zeta,\sigma \in N^*_{y_2} \subseteq N'_\delta
\subseteq N_\zeta$ so sup$(N_\zeta \cap \sigma) = e_{h(\sigma)}(\zeta)$ so
by the previous sentence $e_{h(\sigma)}(\zeta) \in N^*_{y_2}$, hence
$e_{h(\sigma)}(\zeta) \in \cup\{N^*_x:x \in {\Cal P}_\delta\} = N'_\delta$
as required.
\sn
\ub{Clause $(iv)$}:  By clause $(iii)$ it is $\subseteq N'_\delta$,
and by the choice of the function $e$ it is $\subseteq \sigma$ hence
it is $\subseteq N'_\delta \cap \sigma$.  Now $N'_\delta =
\cup\{N^*_z:z \in {\Cal P}_\delta\}$ and $z \in {\Cal P}_\delta
\Rightarrow N^*_z \prec N_\delta$ by $(*)_0$ hence $N'_\delta
\subseteq N_\delta$.  Now we know that $\langle
e_{h(\sigma)}(\zeta):\zeta < \delta\rangle$ is increasing with limit
$e_{h(\sigma)}(\delta) = \sup(N_\delta \cap \sigma)$ hence is
unbounded in it and even $\langle e_{h(\sigma)}(\zeta):\zeta \in
C_\delta\rangle$ is an unbounded subset of $e_{h(\sigma)}(\delta)$ and
it is included in $N'_\delta$ as required.

So $(*)_2$ indeed holds.

Now (A), $(B)$, $(C)$, $(D)$, $(E)$ below clearly suffice to finish.
\mr
\item "{$(A)$}"  $(a) \quad$ for $\delta \in S,y \in {\Cal P}_\delta$
and $n < \omega$ we have $M_{y,\delta,n} 
\subseteq N'_\delta = \dbcu_{z \in {\Cal P}_\delta} N^\ast_z$.
\ermn
[Why?  We prove this by induction on $n$.  First assume
$n=0,M_{y,\delta,n}$ is the Skolem hull of $y \cup (N'_\delta \cap
\kappa)$ in the model
$M^*_\lambda$, well defined as $y \subseteq \lambda$ hence $y
\subseteq M^*_\lambda$ and $N' \cap \kappa \subseteq \kappa \subseteq
\lambda$.  As $y \subseteq N^*_y \subseteq N'_\delta$ and $M^*_\lambda \in
N^*_y \subseteq N'_\delta$ clearly $M_{y,\delta,n} \subseteq
N'_\delta$.  Second, assume $n=m+1$ and $M_{y,\delta,m} \subseteq
N'_\delta$.  Now $M_{y,\delta,n}$ in the Skolem hull of
$M_{y,\delta,m} \cup \{e_{h(\sigma)}(\zeta):\sigma \in M_{y,\delta,m}
\cap \text{ Reg} \cap (\lambda^+ \backslash \theta^+)$ and $\zeta \in
y\}$, so it is enough to show that: if
$\sigma \in M_{y,\delta,m}$ (hence $\sigma \in N'_\delta$) and 
$\sigma \in \text{ Reg } \cap \lambda^+ \backslash \theta^+$ and
$\zeta \in y$ then $e_{h(\sigma)}(\zeta) \in N'_\delta$.  But by
$(*)_2(iii)$ this holds.
\mr
\item "{${{}}$}"  $(b) \quad$ for  $z \subseteq y$ in  
${\Cal P}_\delta $ we have $M_{z,\delta,n} \subseteq M_{y,\delta ,n}$.
\ermn
[Why?  Just by their choice, i.e. we prove this by induction on $n < \omega$.]
\mr
\item "{${{}}$}"   $(c) \quad$  for $y \in {\Cal P}_\delta$ 
and $m \le n$ we have $M_{y,\delta,m} \subseteq M_{y,\delta,n}$.
\ermn
[Why?  Just by their choice, i.e. we prove this by induction on $n$.] 
\mr
\item "{${{}}$}"  $(d) \quad M'_\delta := \cup\{M_{y,\delta,n}:y \in
{\Cal P}_\delta$ and $n < \omega\}$ is $\prec N'_\delta$.
\ermn
[Why?  By the above.]
\mr
\item "{${{}}$}"  $(e) \quad$ if $\zeta \in z$ (hence $\zeta \in C_\delta
\subseteq E$),  $\{y,z\} \subseteq {\Cal P}_\delta$, $\sup(y) < \zeta,
y <_{{\Cal P}_\delta} z$
\nl

\hskip25pt  and $\sigma \in \text{ Reg } \cap \lambda^+ 
\backslash \theta^+$ then: $\sigma \in N^\ast_y \prec  
N_\zeta \Rightarrow e_{h(\sigma)}(\zeta)$
\nl

\hskip25pt $ = \sup(\sigma \cap N_\zeta) \in N^\ast _z$.
\ermn
[Why?  By $(*)_2(i) + (ii)$ this holds.]
\mr
\item "{$(B)$}"  We can also prove that 
$\langle M_{y,\delta ,n}:n < \omega ,y \in  
{\Cal P}_\delta \rangle$ is definable in $({\Cal H}(\chi),\in,
<^\ast_\chi)$  
from the parameters $\delta,M^\ast_\lambda,(\bar C,\bar{\Cal P})$ 
and $h \restriction a_i$, all of them belong to $M^\ast$, hence 
the sequence, and $M'_\delta = \cup\{M_{y,\delta,n}:n < \omega,y 
\in {\Cal P}_\delta\}$, belongs to $M^\ast$
\sn
\item "{$(C)$}"   $M'_\delta \cap \text{ Reg } \cap 
(\theta,\lambda^+)$ is a subset of ${\frak a}_\delta$.
\ermn
[Why?  Use (A)(a) and definition of $a_i,{\frak a}_i)$.]
\mr
\item "{$(D)$}"    if $\sigma \in M'_\delta$ 
and $\sigma  \in \text{ Reg } \cap  \lambda^+ \backslash
\kappa$  then $\sigma \cap M'_\delta$ is unbounded in $\sigma \cap N'_\delta$.
\ermn
[Why?  When $\sigma  > \theta$ use $(*)_2(iii),(iv)$.  For $\sigma = \theta$ we
have $N'_\delta \cap \theta \subseteq N_\delta \cap \theta = \delta$ as
$\delta \in E$ and $C_\delta \subseteq \delta = \sup(C_\delta)$ so it
is enough to show $C_\delta \subseteq N'_\delta$, but $C_\delta$ is 
equal to $\dbcu_{y \in {\Cal P}_\delta} y$.  For $\sigma = \kappa$  
see the choice of $M_{y,\delta,0}$.  So as $\theta = \kappa^+$ we are done.]
\mr
\item "{$(E)$}"  $M'_\delta \cap \lambda = N'_y \cap \lambda$.  
\ermn
[Why?  By (A)(a) we have one inclusion, the $\subseteq$.  By the
choice of $M^*_\lambda$ and clause (D) the result follows by 
\cite[3.3A,5.1A]{Sh:400} recalling $N'_\delta \cap \kappa \in \kappa$.] 
${{}}$  \hfill$\square_{\scite{2.6}}$
\enddemo
\bn
But to get normality of the filter we better define
\definition{\stag{2b.4A} Definition}  Assume $\theta = 
\text{ cf}(\theta) > \kappa = \text{ cf}(\kappa) > \aleph_0,(\bar C,\bar{\Cal
P}) \in {\Cal T}^*[\theta,\kappa]$ and $X$ is a set, of cardinality $\ge
\theta$ for simplicity and let $\chi$ be large enough.  We define a
filter ${\Cal D}_{(\bar C,\bar{\Cal P})}[X]$ on $[X]^{< \kappa}$ as
the set of $Y \subseteq [X]^{< \kappa}$ such that for some $\bold x
\in {\Cal H}(\chi)$, for every sequence $\langle N_\alpha,N^*_a:\alpha
< \theta,a \in \dbcu_{\delta \in S} {\Cal P}_\delta\rangle$ satisfying
$\otimes$ below, there is $A \in \text{ id}^a(\bar C)$ such that
$\bold x \in \dbcu_{a \in {\Cal P}_\delta} N^*_a \and \delta \in S(\bar
C) \backslash A \Rightarrow \dbcu_{a \in {\Cal P}_\delta} N^*_a \cap X
\in Y$ where
\mr
\item "{$\otimes$}"  as in Definition \scite{2.4} omitting $\bold x
\in N_\alpha$.
\endroster
\enddefinition
\bigskip

\proclaim{\stag{2b.5A} Claim}  Let $(\bar C,\bar{\Cal P}) \in 
{\Cal T}^*[\theta,\kappa]$.
\nl
1) An $\chi$ such that ${\Cal P}(X) \subseteq {\Cal H}(\chi)$ can
   serve in Definition \scite{2b.4A}, and $\bold x = Y$ can serve.
\nl
2) If $X_1,X_2$ are sets of cardinality $\lambda \ge \chi$ and $f$ is
   a one-to-one function from $X_1$ onto $X_2$, \ub{then} $f$ maps
   ${\Cal D}_{(\bar C,\bar{\Cal P})}(X_1)$ onto ${\Cal D}_{(\bar
   C,\bar{\Cal P})}(X_2)$.
\nl
3) If $X_1 \subseteq X_2$ has cardinality $\ge \theta$ then $Y \in
   {\Cal D}_{(\bar C,\bar{\Cal P})}[X_1] \Rightarrow \{u \in [X_2]^{<
   \kappa}:u \cap X_1 \in Y\} \in {\Cal D}_{(\bar C,\bar{\Cal
   P})}[X_2]$ and $Y \in {\Cal D}_{(\bar C,\bar{\Cal P})}(X_2)
   \Rightarrow \{u \cap X_1:u \in Y\} \in {\Cal D}_{(\bar C,\bar{\Cal
   P})}(X_1)$.
\nl
2) For any set $X$ of cardinality $\ge \kappa$, really ${\Cal
   D}_{(\bar C,\bar{\Cal P})}(X)$ is a fine normal filter on $X$,
   i.e.:
\mr
\item "{$(a)$}"  fine: $t \in X \Rightarrow \{u \in [X]^{< \kappa}:t
   \in u\} \in {\Cal D}_{(\bar C,\bar{\Cal P})}(X)$
\sn
\item "{$(b)$}"  normal: if $Y_t \in {\Cal D}_{(\bar C,\bar{\Cal
P})}(X)$ for $t \in X$ then $Y := \Delta\{Y_t:t \in X\} = \{u \in 
[X]^{< kappa}:u \ne \emptyset$ and $t \in u \Rightarrow u \in Y_t\}$.
\endroster
\endproclaim
\bigskip

\demo{Proof}  1),2) Easy.
\nl
3) The ``fine" is trivial and for normal let $\bold x_t$ be a witness for
   $Y_t \in {\Cal D}_{(\bar C,\bar{\Cal P})}[X]$ now $\bold x =
   \langle \bold x_t:t \in x\rangle$ witness that $Y \in {\Cal
   D}_{(\bar C,\bar{\Cal P})}[X]$.
\enddemo
\bigskip

\proclaim{\stag{2b.5D} Claim}  Let $(\bar C,\bar{\Cal P}) \in 
{\Cal T}^*[\theta,\kappa]$.
\nl
1) ${\Cal D}_{(\bar C,\bar{\Cal P})}(\lambda) \supseteq
{\Cal D}_{(\bar C,\bar{\Cal P})}[\lambda]$.
\nl
2) In \scite{2.6} we can replace ${\Cal D}_{(\bar C,\bar{\Cal
   P})}(\lambda)$ by ${\Cal D}_{(\bar C,\bar{\Cal P})}[\lambda]$.
\nl
3) Assume that {\rm cf}$(\lambda) \ge \kappa$ and $\beta < \alpha
   \Rightarrow \lambda > \text{\rm cov}(|\beta|,\kappa,\kappa,2)$.
   \ub{Then} there is $S \in {\Cal D}_{(\bar C,\bar{\Cal
   P})}(\lambda)$ such that $\alpha < S \Rightarrow \lambda > |\{u \in
   S:u \subseteq \alpha\}|$.
\endproclaim
\bigskip

\demo{Proof}  1) Trivial.
\nl
2) Repeat the proof, the change is minor.
\nl
3) We can find ${\Cal Q} = \{u_i:i < \lambda\} \subseteq [\lambda]^{<
   \kappa}$ which is cofinal such that $\forall \alpha <
   \lambda(\beta)(\alpha)[\alpha \le \beta < \lambda \wedge [\{u_i:i <
   \beta,u_i \subseteq \alpha\}]$ is cofinal in $[\alpha]^{< \kappa}$.
\enddemo
\bigskip

\remark{\stag{2.6D} Remark}  In \scite{2.6} we can replace $\theta =
\kappa^+$ by $\theta > \kappa_\sigma > \sigma = \text{\rm cf}(\sigma)$
and $\alpha < \theta \Rightarrow |\alpha|^{<\sigma>_{\text{tr}}} <
\theta$ and $\delta \in S(\bar C) \Rightarrow \text{\rm cf}(\delta) = \sigma$.
\endremark
\bigskip

\demo{Proof}  Fill.
\enddemo
\bigskip

\demo{\stag{2.7} Conclusion}   Suppose $\lambda > \kappa > 
\aleph _0$ are regular cardinals and  $(\forall \mu  < \lambda)
[\text{cov}(\mu ,\kappa ,\kappa ,2) < \lambda ].$ 
\nl
1) If for  $\alpha  < \lambda $,  $a_\alpha $ is a subset of  $\lambda $  of 
cardinality $< \kappa$ and $S \in {\Cal D}_{<\kappa}(\lambda)$ and
$T_1 \subseteq \{\delta < \lambda:\text{cf}(\delta) \ge \kappa\}$ is
stationary, \ub{then} we can find a stationary $T_2 \subseteq T_1,
c \subseteq \lambda$ 
and $\langle b_\delta:\delta \in T\rangle$ such that:

$$
a_\delta \subseteq  b_\delta \in  S \text{ for } \delta  \in  T_2
$$

$$
b_\delta  \cap  \delta  = c \text{ for } \delta  \in  T_2.
$$
\mn
2) If in addition $(\bar C,\bar{\Cal P}) \in {\Cal
   T}^*[\kappa^+,\kappa]$ and $S \in ({\Cal D}_{(\bar C,\bar{\Cal
   P})}(\lambda))^+$ \ub{then} part (1) holds for this $S$.
\enddemo
\bigskip

\remark{Remark}  See on this and on \scite{2.9} Rubin Shelah
\cite[4.12,pg.76]{RuSh:117} and \cite[\S6]{Sh:371}.
There we do not know that $(\forall \mu <
\lambda)[\text{cov}(\mu,\kappa,\kappa,2) < \lambda]$ implies (as
proved ehre) that
\mr
\item "{$\boxtimes_{\lambda,\kappa}$}"  for each $\alpha < \lambda$ we
can find $S_\alpha$ a stationary $S_\alpha \subseteq [\alpha]^{<
\lambda}$ of cardinality $< \lambda$; moreover such that $\{\alpha\}
\cup u:u \in S_\alpha,\alpha < \lambda\} \subseteq [\lambda]^{<
\kappa}$ is stationary, (if $\lambda$ is a successor cardinal, the
moreover follows.  So the assumption there seems just what was used now.
 So we could just quote. 
\endroster
\endremark
\bigskip

\demo{Proof}  1) By part (2).
\nl
2) For each $\alpha < \lambda$ let $S_\alpha \in {\Cal D}_{(\bar
   C,\bar{\Cal P})}[\alpha]$ be of cardinality
   cov$(|\alpha|,\kappa,\kappa,2)$.

Let $S = \{u \in [\lambda]^{< \kappa}$: if $\alpha \in u \backslash
\kappa^+$ then $u \cap \alpha \in S_\alpha\}$, so by \scite{2b.5A} we
know that $S \in {\Cal D}_{(\bar C,\bar{\Cal P})}[\lambda]$; and by
\scite{2b.5D}(3) \wilog
\mr
\item "{$(*)$}"  $\alpha < \lambda \Rightarrow \{u \in S:u \subseteq
\alpha\}$ has cardinality $< \lambda$.
\ermn
Now for each $\alpha < \lambda$ let $b_\alpha \in S$ be such that
$a_\alpha \subseteq b_\delta$, clearly exist and let $h:T_1
\rightarrow \lambda$ be defined by $h(\delta) = \sup(b_\delta \cap
\delta)$ so $\delta \in T_2 \Rightarrow h(\delta) < \delta$ as
cf$(\delta) \ge \kappa > |b_\delta|$.  So for some $\gamma_* < \gamma$
the set $T'_2 := \{\delta \in T_1:h(\delta) = \gamma_*\}$ is
stationary and by $(*)$ for some $c$ the set $T_2 := \{\delta \in
T'_2:b_\delta \cap \delta = c\}$ is stationary.  Let $E = \{\delta <
\lambda:\delta$ is a limit ordinal such that $\alpha < \delta
\Rightarrow b_\alpha \subseteq \delta\}$, it is a club of $\lambda$.
\enddemo
\bigskip

\demo{\stag{2.8} Conclusion}   If  $\lambda > \kappa > \aleph_0,\lambda$  and  
$\kappa $  are regular cardinals and  $[\kappa  < \mu  < \lambda  \Rightarrow  
\text{ cov}(\mu ,\kappa ,\kappa ,2) < \lambda ]$ \ub{then} $\{\delta  < 
\lambda :\text{cf}(\delta ) < \kappa \} \in \check I[\lambda ].$
\enddemo
\bigskip

\demo{Proof}  Use  $\mu(3)$  of \scite{2.6}.
\enddemo
\bigskip

\proclaim{\stag{2.9} Claim}  Let $(\ast )_{\mu ,\lambda ,\kappa }$
mean: if  $a_i \in [\lambda]^{<\kappa}$ for $i \in S$ and 
$S \subseteq \{\delta < \mu:{\text{\rm cf\/}}(\delta) = \kappa\}$ 
is stationary, \ub{then} for some  $b \in [\lambda]^{<\kappa}$ the set
$\{i \in S:a_i \cap  i \subseteq  b\}$  is 
stationary.  Let $(\ast )^-_{\mu ,\lambda ,\kappa }$ be defined similarly but  
$\{i \in  S:a_i \subseteq  b\}$  only unbounded. \nl
\ub{Then} for $\aleph_0 < \kappa  < \lambda  < \mu$ regular we have:  

$$
\align
{\text{\rm cov\/}}(\lambda,\kappa,\kappa,2) < \mu &\Rightarrow  
(\ast )_{\mu ,\lambda ,\kappa } \Rightarrow 
(\ast )^-_{\mu,\lambda,\kappa} \\
  &\Rightarrow (\forall \lambda')[\kappa < \lambda ' \le  
\lambda \and \,{\text{\rm cf\/}}(\lambda') < \kappa \Rightarrow \,
{\text{\rm pp\/}}_{<\kappa}(\lambda') < \mu].
\endalign
$$
\endproclaim
\bigskip

\remark{Remark}   So it is conceivable that the $\Rightarrow$ are  
$\Leftrightarrow$.  See \cite[\S3]{Sh:430}.
\endremark
\bigskip

\demo{Proof}   Straightforward.  \hfill$\square_{\scite{2.9}}$ 
\enddemo
\bn
\ub{Exercise}:  Generalize to the following filter.

Let $\theta = \text{ cf}(\theta) \ge \kappa = \text{ cf}(\kappa)$ and
$S_* \subseteq [\theta]^{< \kappa}$ be stationary.  For any set $X$ of
cardinality $\ge \theta$ we define a filter ${\Cal D}^1_{S_*}[X]$ as
follows: $Y \in {\Cal D}_{S_*}[X]$ \ub{iff} $Y \subseteq [X]^{<
\kappa}$ and for any $\chi$ large enough there is $\bold x \in {\Cal
H}(\chi)$ such that if $\langle N_\alpha,f_\alpha:\alpha \le
\theta\rangle$ satisfy $\circledast$ below, then for some $S' \in
{\Cal D}_{< \kappa}(\theta)$ for every $u \in S_* \cap S'$ we have:

if $\bold x \in f''_\theta(u)$ \ub{then} $f''(u) \in Y$, when:
\mr
\item "{$\circledast$}"  $(a) \quad N_\alpha \prec ({\Cal
H}(\chi),\in,<^*_\chi)$
\sn
\item "{${{}}$}"  $(b) \quad N_\alpha$ is $\prec$-increasing
continuous
\sn
\item "{${{}}$}"  $(c) \quad \|N_\alpha\| < |\alpha|^+ + \theta$
\sn
\item "{${{}}$}"  $(d) \quad \langle N_\beta:\beta \le \alpha\rangle
\in N_{\alpha +1}$ if $\alpha < \theta$
\sn
\item "{${{}}$}"  $(e) \quad$ can add $\langle
\kappa,\theta,X,S_*\rangle \in N_0$.
\endroster
\newpage

\head {\S3  Nice Filters Revisited} \endhead  \resetall \sectno=3
 \spuriousreset
\bigskip

This generalizes \cite{Sh:386} (and see there). 
\nl
See \cite[\S5]{Sh:410} on this generalization of normal filters.
\bigskip

\demo{\stag{3.0} Convention}  1) $\bold n$ is a niceness context; we
use $\kappa$, FILL, etc., for $\kappa_{\bold n}$, Fil$_{\bold n} =
\text{ FIL}(\bold n)$ when dealing from the content.
\enddemo
\bigskip

\definition{\stag{3.1} Definition}  We say the $\bold n$ is a niceness
context or a $\kappa$-niceness context or a $(\kappa,\mu)$-niceness
context if it consists of the following objects satisfying the
following conditions:
\mr
\item "{$(a)$}"  $\kappa$ is a regular uncountable cardinal
\sn
\item "{$(b)$}"  $I \subseteq {}^{\omega >} \omega$ is non-empty
$\triangleleft$-downward closed with no $\triangleleft$-maximal member
\footnote{For ${\Cal T}$ the two interesting cases are ${\Cal T}=
{}^{\omega >}\omega$ and ${\Cal T} = \{<>\}$ and ${}^{\omega
>}\{0\}$.  The default value will be ${}^{\omega >}\omega$.}
default value is $\{0_n:n < \omega\}$
\sn
\item "{$(c)$}"  let $\mu$ be $ > \kappa$ and $\langle{\Cal Y}:i <
\kappa\rangle$ is a sequence of pairwise disjoint sets and 
${\Cal Y} \cup \{{\Cal Y}_i:i < \omega_1\}$ so 
$i < \omega_1 \Rightarrow |{\Cal Y}|,|{\Cal Y}_i|$ 
\sn
\item "{$(d)$}"  the function $\iota$ with domain ${\Cal Y}$ is
defined by $\iota(y) = i$ when $y \in {\Cal Y}_i$
\sn
\item "{$(e)$}"  $\bold e$ is a 
set of equivalence relations $e$ on ${\Cal Y}$  
refining  $\dbcu_{i < \omega_1} {\Cal Y}_i \times {\Cal Y}_i$ with  
$< \mu^\ast$ equivalence classes, each class of cardinality 
$|{\Cal Y}|$
\sn
\item "{$(f)$}"  for $e \in \bold e$, FIL$(e) = \text{ FIL}(e,\bold
n)$ is a set of $D$  such that: 
{\roster
\itemitem{ $(\alpha)$ }     $D$  is a filter on  ${\Cal Y}/e$, 
\sn
\itemitem{ $(\beta)$ }    for any club $C$ of $\kappa$ we have
$\dbcu_{i \in C} {\Cal Y}_i/e \in D$,
\sn
\itemitem{ $(\gamma)$ }   normality: if $X_i \in D$ for $i < \omega_1$
\ub{then} the following set belongs to $D$: \nl  

$\quad \{(\delta,j)/e:(\delta,j) \in {\Cal Y},\delta$ limit and  
$i < \delta \Rightarrow (\delta,j) \in X_i\}$
\endroster}
\item "{$(g)$}"  Suc $\in \{(D_1,D_2):e(D_1) \le e(D_2)\}$.
\endroster
\enddefinition
\bigskip

\remark{Remark}  For $\bold e$ an important case is when it is a
singleton $\{\cup\{{\Cal Y}_i \times {\Cal Y}_i:i < \kappa\}\}$, so we
are dealing with normal filters on the old case.
\endremark
\bigskip

\definition{\stag{3.2} Definition}  Let $\bold n$ be a
$\kappa$-niceness context. 
\nl
1) We say  $e_1 \le e_2$ if $e_2$ refines $e_1$.  
If not said otherwise, every  $e$  is from $\bold e$.  Let 
$\bold e_\mu$ be the set of all such equivalence relations with  $< \mu $  
equivalence classes.  Let  $\iota(x/e) = \iota(x)$.
\nl
2) FIL $ = \text{ FIL}(\bold n)$ is $\dbcu_{e \in \bold e} 
\text{ FIL}(e,\bold n)$.  For  $D \in \text{ FIL}$, let 
$e = e[D]$ be the unique $e \in \bold e$ such that $D \in \text{ FIL}
(e,\bold n)$. \nl
3)  For $D \in \text{ FIL}(e)$ let $D^{[\ast]} = \{ X \subseteq {\Cal
Y}:X^{[*]} \in D\}$; see (5) below.
\nl
4)  For $D \in \text{ FIL}(\bold n)$ and $e(1) \ge e(D)$,  
let  $D^{[e(1)]} = \{ X \subseteq {\Cal Y}/e(1):X^{[\ast]} \in
D^{[\ast]}\}$, see (5) below.
\nl
5)  For  $A \subseteq {\Cal Y}/e,A^{[\ast]} = \{(x /e):(x/e) \in A\}$,  and 
for $e(1) \ge e$ let $A^{[e(1)]} = \{y/e(1):y/e \in A\}$.
\enddefinition
\bigskip

\definition{\stag{3.2A} Definition}  1) For $D \in \text{ FIL}(e,\bold
n)$, let $D^+$ be $\{Y \subseteq {\Cal Y}/e:Y \ne \emptyset$  mod
$D\}$. \nl
2) $\bold n$ is 1-closed if $D \in \text{ FIL}(\bold n),A \in D^+
\Rightarrow D+A \in \text{ FIL}(\bold n)$. \nl
3) $\bold n$ is 0-closed if for every $D_1 \in \text{ FIL}_{\bold n}$
and $A \in D^+_1$ there is $D_2 \in \text{ FIL}_2$ such that $(D_1 +A)
\in (D_2) \subseteq D_2$.
\nl
4) A niceness context $\bold n$ is full \ub{if}
\mr
\item "{$(a)$}"  for every $e \in \bold e_{\bold n}$, every filter on
${\Cal Y}_{\bold n}/e$ which is normal (with respect to the function
$\iota_{\bold n}$) belong to FIL$_{\bold n}(e)$.
\ermn
4A) A niceness content $\bold n$ is semi-full when: for every $e_1 \in
\bold e_{\bold n}$ and $D_1 \in \text{ FIL}_{\bold n}(e_1)$ and
$e_2,e_1 \le e_2 \in \bold e_{\bold n}$ and ${\Cal A} \subset {\Cal
P}({\Cal Y}_{\bold n}/e_2)$ lift$(W) \in \text{ FIL}(e_2)$ whenever
\mr
\item "{$(*)_{e_1,e_2,D_1,W}$}"  $(a) \quad e_1 \le e_2$ in $\bold
e_{\bold n}$
\sn
\item "{${{}}$}"  $(b) \quad D_1 \in \text{ FIL}_n(e_2)$
\sn
\item "{${{}}$}"  $(c) \quad \mu \ge 2^{({\Cal Y}/e_2)}$ (or more ???)
\sn
\item "{${{}}$}"  $(d) \quad W \subseteq [\mu]^{\le \aleph_0}$ is
stationary
\sn
\item "{${{}}$}"  $(e) \quad D_2 = \text{ lift}(W,D^{[e_2]}_1)$ is
normal (i.e. $\emptyset \in \text{ lift}(W,D_1))$.
\ermn
5) A niceness context $\bold n$ is \ub{thin} when

$$
\align
\text{Suc}_{\bold n} = \{(D_1,D_2):&D_1 = D_2 \in \text{ FIL}_{\bold n}
\text{ and} \\
  &D_2 = D^{[e_1]}_1 +A \text{ for some } A \in (D^{[e_1]}_1)^+\}.
\endalign
$$  
\mn
6) A niceness context $\bold n$ is thick if: Suc$_{\bold n} =
\{(D_1,D_2):D_1,D_2 \in \text{ FIL}_{\bold n},e(D_1) \le e(D_2)$ and
$D^{[e_2]}_1 \subseteq D_2$ and if $\mu = 2^{|{\Cal Y}_{\bold n}/e_2)},
W_1 \subseteq [\mu]^{\le \aleph_0}$ is stationary and
lift$(W,D_1) = D_1$ then for some stationary $W_2 \subseteq W_1$ we
have lift$(W_2,D_2) = D_2\}$.
\enddefinition
\bigskip

\remark{Remark}  1) On lift see Definition \scite{3.9B}, HERE?? \nl
2) We can use more freedom in the higher objects.
\endremark
\bigskip

\proclaim{\stag{3.2K} Claim}  Assume
\mr
\item "{$(a)$}"  the $\kappa$-niceness context is thick
\sn
\item "{$(b)$}"  $D_1 \in { \text{\rm FIL\/}}_{\bold n}(e_1)$
\sn
\item "{$(c)$}"  $e_1 \le e_2 \in \bold e_{\bold d}$
\sn
\item "{$(d)$}"  for each $y \in {\Cal Y}_{\bold n}/e_1,\langle
z_{y,\varepsilon}:\varepsilon < \varepsilon_y \rangle$ list $\{z/e_2:z
\in y_1\},d_{y,\varepsilon}$ is a $\kappa$-complete filter on
$\varepsilon_y$
\sn
\item "{$(e)$}"  $D_2 \in { \text{\rm FIL\/}}_{\bold n}(e_2)$
\sn
\item "{$(f)$}"  if $A \in D_2$ then $\{y \in {\Cal Y}_{\bold
n}/e_1:\{\varepsilon < \varepsilon_y:z_{y,\varepsilon} \in A\} \in
d_{y,\varepsilon}\}$ belongs to $D_1$.
\ermn
Then $D_2 \in { \text{\rm Suc\/}}_{\bold n}(D_1)$.
\endproclaim
\bn
\ub{Discussion}:  We may consider allowing player $I$, in the
beginning of each move to choose $W_n$ as above.
\bigskip

\definition{\stag{3.3} Definition}  (0) For $f:{\Cal Y}/e \rightarrow  X$  let  
$f^{[\ast ]}:{\Cal Y} \rightarrow  X$ be $f^{[\ast ]}(x) = f(x/e)$.  We say  
$f:{\Cal Y} \rightarrow  X$  is supported by  $e$  if it has the form  
$g^{[\ast]}$ for some  $g:{\Cal Y}/e \rightarrow  X$.  If $e_1,e_2 \in  
\bold e$ and $f_\ell :{\Cal Y}/e_\ell  \rightarrow X$ for $\ell =1,2$ then:  
we say $f_1 = f^{[e_1]}_2$ 
if  $f^{[\ast]}_1 = f^{[\ast]}_2$.  Writing  $f^{[\ast ]}$ for  $f \in  
{}^{\omega_1}X$  we identify  $\{i\}$, $i < \omega _1$ with  ${\Cal Y}_i.$ 
\nl
(1) Let  $F_c({\Cal T},e) = F_c({\Cal T},e,{\Cal Y})$  be the 
family of  $\bar g$,  a sequence of the form  $\langle g_\eta :\eta  \in  
u \rangle $,  $u \in  f_c({\Cal T}) =$ the family of non-empty finite 
subsets of  ${}^{\omega >}\omega$ closed under taking initial segments,  and 
for each  $\eta  \in  u$  we have $g_\eta \in {}^{\Cal Y}\text{Ord}$  
is supported by $e$.  Let Dom$(\bar g) = u$,  Range$(\bar g) 
= \{g_\eta:\eta \in u\}$.  We let  $e = e(\bar g)$,  for the 
minimal possible  $e$  assuming it exists and we
shall say  $g_\eta <_D g_\nu$ instead $g_\eta <_{D^{[\ast]}} g_\nu$ and 
not always distinguish between $g \in {}^{{\Cal Y}/e}\text{Ord}$ and 
$g^{[\ast]}$ in an abuse of notation. \nl
(2) We say $\bar g$ is decreasing for $D$ or $D$-decreasing (for  $D \in  
\text{ FIL}(e,I))$ if $\eta \triangleleft \nu 
\Rightarrow g_\nu <_D g_\eta$. \nl
(3) If  $u = \{<>\}$,  $g = g_{<>}$ we may write  $g$  instead  
$\langle g_\eta :\eta  \in  u \rangle$.
\enddefinition
\bigskip

\definition{\stag{3.4} Definition}  1) For $e \in \bold e,D \in 
\text{ FIL}(e)$ and $D$-decreasing $\bar g \in F_c({\Cal T},e)$ we 
define a game $\Game^\ast (D,\bar g,e) = \Game^*(D,\bar g,e,\bold n)$.
In the nth move (stipulating  $e_{-1} = e$, $D_{-1} = D,
\bar g_{-1} = \bar g)$: 
\sn
\ub{the case $\bold n$ is then}
\block
player I chooses $e_n \ge e_{n-1}$ and $A_n \subseteq 
{\Cal Y}/e_n$,  $A_n \ne \emptyset$ mod $D^{[e_n]}_{n-1}$ and he chooses  
$\bar g^n \in F_c({\Cal T},e_n)$ extending $\bar g_{n-1}$ (i.e.  
$\bar g^{n-1} = \bar g^n \restriction \text{ Dom}(\bar g_{n-1})),
\bar g^n$ supported by $e_n$ and 
$\bar g^n$ is $(D^{[e_n]}_n + A_n)$-decreasing, 
player II chooses 
$D_n \in \text{ FIL}(e_n)$ extending $D^{[e_n]}_{n-1} + A_n$. 
\endblock
\mn
\ub{In the general case}:
\sn
Player I chooses $e_n$ and $D_{n,1} \in \text{ Duc}_{\bold n}(D_{n-1})$
and let $e_n = e(D_{n-1})$ and he chooses $\bar g^n \in F \subset
({\Cal T},e(D_{n-1})$ which is extending $\bar g^{n-1}$ then $\eta \in
\text{ Dom}(\bar g^n)$ (i.e. $\bar g^{n-1} = \bar g^n \restriction
\text{ Dom}(\bar g^{n-1}),\bar g^n$ supported by $e(D_{n,1})$ and
$\bar g^n$ is $D_{n,1}$-decreasing.
\sn
Player II chooses $D_n = D_{n,2} \in \text{ FIL}(\bold e_n)$
extending $D_{n,1}$.
\sn
In the end, the second player wins if 
$\dbcu_{n < \omega} \text{ Dom}(\bar g^n)$ has no infinite branch.
\nl
2) Let $\bar \gamma$ be such that Dom$(\bar \gamma) = \text{ Dom}(\bar
g)$ and each $\gamma_\eta$ is an ordinal decreasing with $\eta$.  Now
$\Game^{\bar \gamma}(D,\bar g,e)$ is defined similarly to  
$\Game^\ast(D,\bar g,e)$  but the second player has in 
addition, to choose an ordinal  $\alpha_\eta$ for  
$\eta \in \text{ Dom}(\bar g^n) \backslash 
\dbcu_{\ell < n} \text{ Dom}(\bar g^\ell)$ such that 
$[\eta \triangleleft \nu \and \nu \in \text{ Dom}(\bar g^{n-1}) 
\Rightarrow \alpha_\nu < \alpha_\eta]$ we let  
$\alpha _\eta  = \gamma_\eta$ for $\eta \in \text{ Dom}(\bar g)$. \nl
3)  $w \Game^\ast (D,\bar g,e)$ and  
$w \Game^{\bar \gamma}(D,\bar g,e)$ 
are defined similarly but $e$ is not changed during a play. 
(If e.g. $\bold e = \{e\}$ then this makes not difference.)
\nl
4)  If  $\bar \gamma  = \langle \gamma _{<>}\rangle $,  $\bar g = 
\langle g_{<>}\rangle $  we write  $\gamma _{<>}$ instead  $\bar \gamma $,  
$g_{<>}$ instead  $\bar g.$
\nl
5)  If  $E \subseteq \text{ FIL}$ the games $\Game^\ast _E$,  
$\Game^{\bar \gamma }_E$ are defined similarly, 
but player II can choose filters only from $E$  
(so we naturally assume to have $A \in D^+,D \in  E \Rightarrow  D + A 
\in  E)$.
\enddefinition
\bigskip

\remark{\stag{3.4A} Remark}   Denote the above games $\Game^\ast_0,
\Game^{\bar \gamma}_0,w \Game^\ast_0$.  Another variant is \nl
3)  For  $e \in \bold e,D \in \text{ FIL}(e)$ and 
$D$-decreasing  $\bar g \in F_c({\Cal T})$ we define a game 
$\Game^\ast_1(D,\bar g,e)$.  We stipulate  $e_{-1} = e$,  $D_{-1} = D.$ 

In the nth move first player chooses $e_n,e_{n-1} \le e_n \in {\Cal T}$ and
$D'_n \in \text{ FIL}(e_n)$ and $D'_n$-decreasing  
$\bar g^n$ extending  $\bar g^{n-1}$ such that $(D_{n-1} +
A_n)^{[e_n]} \subseteq D_n$ and: 
\mr
\item "{$(\ast)$}"   for some  $A_n \subseteq {\Cal Y}/e_{n-1}$,  $A_n
\ne \emptyset$ mod $D_{n-1}$ we have:
{\roster
\itemitem{ $(i)$ }  $D'_n$  is the normal filter on ${\Cal Y}/e_n$ 
generated by  $(D_{n-1} + A_n)^{[e_n]} \cup  \{A^n_\zeta:\zeta <
\zeta^\ast_n\}$  where for some  
$\langle C_\zeta:\zeta < \zeta_n \rangle$ we have:     
\sn
\itemitem{ $(a)$ }  each $C_\zeta $ is a club of  $\omega _1,$     
\sn
\itemitem{ $(b)$ }   if $\zeta_\ell < \zeta^\ast_n$ for $\ell <
\omega$,  $i \in  \dbca_{\ell < \omega} C_{\zeta_\ell}$,  
$x \in  {\Cal Y}/e_{n-1}$,  and  $\iota (x) = i$,  
then for some  $x' \in  {\Cal Y}/e_n$,  we have  $x' \subseteq x$,  
$x' \in \dbca_{\ell < \omega} A^n_{\zeta _\ell }$.
\endroster}
\ermn
The first player also chooses $\bar g^n$ extending $\bar
g^{n-1},D'_n$-decreasing. 
Then second player chooses $D_n$ such that $D'_n \subseteq D_n \in
\text{ FIL}(e_n)$.
\nl
2)  We define $\Game^{\bar \gamma}_1(D,\bar g,e)$ as in (2) using  
$\Game^\ast_1$ instead of $\Game^\ast_0$.
\nl
3)  If player II wins, e.g. $\Game^{\bar \gamma}_E(D,\bar f,e)$ this 
is true for \nl
$E' =: \{D' \in G:\text{ player II wins } 
\Game^{\bar \gamma}_{E^\ast }(D',\bar f,e)\}$.
\endremark
\bigskip

\definition{\stag{3.5} Definition}  1) We say $D \in \text{ FIL}$ is
nice to $\bar g \in F_c({\Cal T},e,{\Cal Y})$,  $e =
e(D)$, if player II wins the game $\Game^\ast (D,\bar g,e)$  
(so in particular $\bar g$ is $D$-decreasing, $\bar g$ supported by  $e$).
\nl
2)  We say $D \in \text{ FIL}$ is nice \ub{if} it is nice to $\bar g$ for 
every $\bar g \in F_c({\Cal T},e)$.
\nl
3)  We say $D$ is nice to $\alpha$ \ub{if} it is nice to the constant function 
$\alpha $.  We say $D$ is nice to $g \in {}^\kappa\text{Ord}$ \ub{if}
it is nice to $g^{[e(D)]}.$
\nl
4)  ``Weakly nice" is defined similarly but  $e$  is not changed. \nl
5) Above replacing $D$ by $\bold n$ means: for every $D \in 
\text{ FIL}_{\bold n}$.
\enddefinition
\bigskip

\remark{\stag{3.5A} Remark}  ``Nice" in \cite{Sh:386} 
is the weakly nice here, but 
\mr
\item "{$(a)$}"  we can use $\bold n$ with $\bold e_{\bold n} = \{e\}$
\sn
\item "{$(b)$}"  formally they act
on different objects; but if  $x e y \Leftrightarrow \iota(x) =
\iota(y)$  we get a situation isomorphic to the old one.
\endroster
\endremark
\bigskip

\proclaim{\stag{3.6} Claim}  Let $D \in { \text{\rm FIL\/}}$ and $e = e(D).$
\nl
1) If  $D$  is nice to  $f$,  $f \in  F_c({\Cal T},e),g \in F_c
({\Cal T},e)$ and $g \le f$ \ub{then} $D$ is nice to  $f$.
\nl
2)  If  $D$  is nice to  $f$,  $e = e(D) \le e(1) \in \bold e$ then  
$D^{[e(1)]}$ is nice to  $f^{[e(1)]}$.
\nl
3)  The games from \scite{3.4}(2) 
are determined and winning strategies do not need memory.
\nl
4)  $D$  is nice to $\bar g$ \ub{iff} $D$ is nice to $g_{<>}$ (when  $\bar g 
\in F_c({\Cal T},e)$ is $D$-decreasing).
\nl
5)  If $\bold e \subseteq \bold e$ and for simplicity  
$\dbcu_{i < \omega_1} \{i\} \times {\Cal Y}_i \in \bold e$  
and for every  $e \in \bold e,e \le e(1) \in \bold e$ for some
permutation $\pi$ of $\bar{\Cal Y}$ (i.e. a permuation of ${\Cal Y}$
mapping each ${\Cal Y}_i \, (i < \omega_1)$ onto itself) (and $\bold
n$ is full for simplicity) we have $\pi (e) = e,
\pi (e(1)) \le e(2) \in \bold e$ \ub{then} we can replace $\bold e$ by
$\bold e$.
\nl
6) For $\bold e = \bold e_\mu$ (where $\mu \le \mu^\ast)$  
there is  $\bold e$  as above  with:  $|\bold e|$  
countable if  $\mu $  is a successor cardinal  $(> \aleph _1)$,  
$|\bold e| = { \text{\rm cf\/}}(\mu)$ if $\mu $  is a limit cardinal.
\endproclaim
\bigskip

\demo{Proof}  Left to the reader.  (For part (4) use \scite{3.7}(2) below).
\enddemo
\bigskip

\proclaim{\stag{3.7} Claim}   1) Second player wins $\Game^\ast(D,\bar
g,e)$  \ub{iff} for some $\bar \gamma$ second player wins 
$\Game^{\bar \gamma }(D,\bar g,e)$.
\nl
2) If second player wins $\Game^\gamma(D,f,e)$ then for any $D$-decreasing  
$\bar g \in F_c({\Cal T},e),\bar g$ supported by  $e$  
and  $\dsize \bigwedge_{\eta,y} g_\eta(y) \le f(y)$, 
the second player wins in $\Game^{\bar \gamma }(D,\bar g,e)$,  when we let

$$
\gamma_\eta  = \gamma + [\max\{(\ell g(\nu)-\ell g(\eta)+1):
\nu \text{ satisfies } \eta \trianglelefteq \nu \in \text{ Dom}(\bar g)\}].
$$
\mn
3)  If $u_1$, $u_2 \in F_c({\Cal T}),h:u_1 \rightarrow  u_2$ 
satisfies  $[\eta \nu \Leftrightarrow h(\eta)h(\nu)]$ and for $\ell =
1,2$ we have $\bar g^\ell \in F_c({\Cal T},e_2),g^1_\eta \ge 
g^2_{h(\eta )}$ (for  $\eta \in u_1)$,  
$\bar \gamma ^\ell  = \langle \gamma ^\ell _\eta :\eta  
\in  u_\ell \rangle$ is a $\triangleleft$-decreasing sequence of ordinals,  
$\gamma ^2_\eta  \ge \gamma^2_{h(\eta )}$ and the second player wins in  
$\Game^{\bar \gamma ^2}(D,\bar g^2,e)$  then the second 
player wins in $\Game^{\bar \gamma ^1}(D,\bar g^1,e)$.
\endproclaim
\bigskip

\demo{Proof}  1) The ``if part" is trivial, the ``only if part" [FILL]
is as in \cite{Sh:386}.
\nl
2), 3)  Left to the reader.
\enddemo
\bn
The following is a consequence of a theorem of Dodd and Jensen \cite{DoJe81}:
\proclaim{\stag{3.8} Theorem}   If $\lambda$ is a cardinal,  
$S \subseteq \lambda$  \ub{then}: 
\nl
(1) $\bold K[S]$,  the core model, is a model of  $ZFC + (\forall \mu \ge 
\lambda )2^\mu  = \mu ^+$. \nl
(2) If in  $\bold K[S]$ there is no Ramsey cardinal  
$\mu  > \lambda $  (or much 
weaker condition holds) \ub{then} ($\bold K[S],\bold V$)  
satisfies the $\mu $-covering lemma 
for  $\mu  \ge \lambda + \aleph_1$;  i.e. if $B \in \bold V$  is a set of 
ordinals of cardinality $\le \mu$  then there is $B' \in \bold K[S]$
satisfying  $B \subseteq B'$ and $\bold V \models |B'| \le \mu$. \nl
(3) If  $\bold V \models (\exists \mu  \ge \lambda )(\exists \kappa ) 
[\mu^\kappa > \mu^+ > 2^\kappa]$ \ub{then} in $\bold K[S]$  there is a 
Ramsey cardinal  $\mu > \lambda$.
\endproclaim
\bigskip

\proclaim{\stag{3.9} Lemma}   Suppose
\mr
\item "{$(a)$}"  $\bold n$ is a semi-full niceness content thin or
medium $\kappa = \aleph_1$
\sn
\item "{$(b)$}"   $f^* \in {}^\kappa${\rm Ord},  
$\lambda > \lambda_0 =: \sup\{(2^{|{\Cal Y}/e|^{\aleph_0}}):e \in
\bold e_{\bold n}\}$
\sn
\item "{$(c)$}"   for every $A \subseteq \lambda_0$, in $K$ there is a
Ramsey cardindal $> \lambda_0$, \ub{then} for every filter 
$D \in { \text{\rm FIL\/}}_{\bold n}(e)$ is nice to $f^*$.
\endroster
\endproclaim
\bigskip

\remark{Remark}  1) The point 
in the proof is that via forcing we translate the filters 
from FIL$(e,{\Cal Y})$  to normal filters on $\kappa$ [for higher  
$\kappa$'s cardinal restrictions are better]. \nl
2) At present we do not care too much what is the value of
$\lambda_0$, i.e., equivalently, how much we like the set $S$ to code. \nl
Saharon: compare with \cite[V]{Sh:g}, i.e., improve as there!  But if
we use $\bold e = \{e\}$, the proofs are more similar to \cite[V]{Sh:g}
we can consider just Levy$(\aleph_1),|D|)$, now in some proofs we may
consider filters generated by $|\text{pcf}({\frak a})|$ set $|{\frak a}| 
< aleph_\omega$.
\endremark
\bigskip

\demo{First Proof}  Without loss of generality $(\forall i)f(i) \ge 2$.
Let  $S \subseteq \lambda _0$ be such that $[\alpha < \mu \and A 
\subseteq 2^{|\alpha|^{\aleph_0}} \Rightarrow  A \in \bold L[S]],\bold
e \in \bold L[S]$  (see \scite{3.6}(6)) and: if  
$g \in {}^\kappa\text{Ord},(\forall i < \kappa_1)g(i) \le  
f(i)$ then $g \in \bold L[S]$  (possible as  $\dsize \prod_{i < \omega_1}
|f(i) + 1| \le \lambda_0$.  We work for awhile in $\bold K[S]$.  In  
$\bold K[S]$ there is a Ramsey cardinal $\mu > \lambda_0$ 
(see \scite{3.8}(3)).  Let in $\bold K[S]$.
\nl
Let

$$
\align
Y_0 = \{X:&X \subseteq \mu,X \cap \kappa \text{ a countable ordinal } > 0,  
\{\kappa,\lambda_0\} \subseteq X, \\
  &\text{ moreover } X \cap \lambda_0 \text{ is countable}\}.
\endalign
$$
\mn
Let

$$
Y_* = Y_1 = \{X \in Y_0:X \text{ has order type } \ge f(X \cap \kappa)\}.
$$
\mn
Now for $g \in {}^\kappa\text{Ord}$ such that  $\dsize \bigwedge_{i <
\omega_1} g(i) < f(i)$  let  $\hat g$  be the function with domain $Y_1$,  
$\hat g(X) =$ the $g(X \cap \kappa)$-th member of  $X$. 

Let $D_* = \{ A_i:\kappa \le i \le 2^{|{\Cal Y}/e|}\}$ and we arrange
$\langle A^D_i:\kappa \le i < 2^{|{\Cal Y}/e|} \rangle \in \bold L[S]$,  
(as  ${\Cal Y}/e$  has cardinality  $< \mu^\ast$, so  $2^{|{\Cal
Y}/e|} \le \lambda_0).$ 

Let $J$  be the minimal fine normal ideal on Y (in $\bold K[S])$  to which  
$Y \backslash Y_D$ belongs where

$$
Y_D = \{X:X \in Y_* \text{ and } i \in (\kappa,2^{|{\Cal Y}/e|}) \cap
X \Rightarrow  X \cap  \omega _1 \in  A_i\}.
$$
\mn
Clearly it is a proper filter as $\bold K[S] \models$ 
``$\mu$ is a Ramsey cardinal".
\enddemo
\bigskip

\demo{\stag{3.9Y} Observation}  Assume
\mr
\item "{$(a)$}"   $\Bbb P$ is a proper forcing notion of 
cardinality  $\le |\alpha|^{\aleph_0}$ for some $\alpha < \mu^\ast$ (or
just $\Bbb P,MAC(\Bbb P) \in \bold K[S]$ and $\{X \in Y_1:X \cap 
(MAC(\Bbb P)|$ is countable$\} \in = Y_*$ mod $J$ where $MAC(\Bbb P)$ is the 
set of maximal antichains of  $\Bbb P)$ and let 
$J^{\Bbb P}$ be the normal fine ideal which $J$ generates in $\bold
V^{\Bbb P}$.  
\ermn
(1) $F$-positiveness is preserved; i.e. if $X \in \bold K[S],X 
\subseteq Y_1,F \in \text{ FIL}$  and $\bold V \models  ``X \ne  
\emptyset$ mod  $F"$  \ub{then}  $\Vdash_{\Bbb P} ``X \ne \emptyset$ mod 
$F^{\Bbb P}$. \nl
(2)  Moreover, if $\Bbb Q \lessdot \Bbb P$, ($\Bbb Q$ proper and)  
$\Bbb P/\Bbb Q$ is proper \ub{then} forcing with $\Bbb P/\Bbb Q$ 
preserve $F^{\Bbb Q}$-positiveness.
\enddemo
\bigskip

\demo{Continuation of the proof of \scite{3.9}}
\sn
\ub{Case 1}: $\bold e = \{e\}$.  Here only \scite{3.9A}(1) is needed
and then it is as in the old case.
\mn
\ub{Case 2}:  General.

Let ${\Cal P}({\Cal Y}/e) = \{ A^e_\zeta:\zeta < 2^{|{\Cal Y}/e|}\}$. 

Now we describe a winning strategy for the second player.  In the side we 
choose also $(p_n,\Gamma_n,{\underset\tilde {}\to f_n})$,  $\bar
\gamma^n,{\underset\tilde {}\to W_n}$ such that
\footnote {For the forcing 
notions actually used below by the homogeneity of the 
forcing notion the value of  $p_n$ is immaterial}
(where  $e_n$,  $A_n$ are chosen by the second player): 
\mr
\item "{$(A)(i)$}"   $\Bbb P_n = \dsize \prod_{\ell \le n} \Bbb
Q_\ell$ where $\Bbb Q_\ell$ is Levy$(\aleph_1,{\Cal Y}/e_n)$      
\nl
(we could use iterations, too, here it does not matter).   
\sn
\item "{$(ii)$}"   $p_n \in \Bbb P_n$   
\sn
\item "{$(iii)$}"  $p_n$ increasing in  $n$   
\sn
\item "{$(iv)$}"  ${\underset\tilde {}\to f_n}$ is a 
$\Bbb P_n$-name of a function from $\omega_1$ to ${\Cal Y}/e_n$   
\sn
\item "{$(v)$}"    $p_n \Vdash_{{\Bbb P}_n}
``{\underset\tilde {}\to f_n}(i) \in  {\Cal Y}_i/e_n"$   
\sn
\item "{$(vi)$}"   $p_{n+1} \Vdash ``{\underset\tilde {}\to f_{n+1}}(i) 
\le {\underset\tilde {}\to f_n}(i)$  for every  $i < \omega_1"$,
\sn
\item "{$(vii)$}"  ${\underset\tilde {}\to f_n}$ is given 
naturally --- it can be interpreted as the generic object of 
$\Bbb Q_n$ except trivialities. 
\sn
\item "{$(B)(i)$}"   $\bar \gamma^n,\bar g^n$ have the same domain,  
$\gamma^n_\eta < \mu$   
\sn
\item "{$(ii)$}"    $p_n \Vdash_{{\Bbb P}_n}
``{\underset\tilde {}\to W_n} \subseteq Y_D$,  
${\underset\tilde {}\to W_{n+1}} \subseteq {\underset\tilde {}\to W_n}"$
\sn
\item "{$(iii)$}"   $\bar \gamma^n = \bar \gamma^{n+1} \restriction
\text{ Dom}(\bar \gamma^n)$, Dom$(\bar \gamma^n) = 
\text{ Dom}(\bar g^n)$ and $\bar \gamma^n$ is $\triangleleft$-decreasing   
\sn
\item "{$(iv)$}"   $p_n \Vdash_{{\Bbb P}_n} 
``\{X \in Y_D:\text{ for } \ell \in \{0,...,n\},
{\underset\tilde {}\to f_\ell} (X \cap \omega_1) \in A_\ell$ and           
$\dsize \bigwedge_{\eta \in \text{ Dom}(\bar g^n)} 
\hat g_\eta(X) = \gamma_\eta$ 
and for  $\ell \in \{-1,0,...,n-1\},\zeta \in X \cap 
2^{|{\Cal Y}/e_\ell|}$ we have:
\nl          
$A^{e_\ell }_\zeta \in D_\ell \Rightarrow {\underset\tilde {}\to
f_\ell} (X \cap \omega_1) \in A^{e_\ell}_\zeta\} \supseteq
{\underset\tilde {}\to W_n} \ne \emptyset$ mod $F^{{\Bbb P}_n}"$   
\sn
\item "{$(v)$}"   $\bar g^n = \bar g^{n+1} \restriction 
\text{ Dom}(\bar g^n)$ [difference] 
\sn
\item "{$(C)(i)$}"   $D_n = \{Z \subseteq {\Cal Y}/e_n:p_n
\Vdash_{{\Bbb P}_n} ``\{X \in J_D:{\underset\tilde {}\to f_n}
(X \cap \omega_1) \notin Z\} = \emptyset$ mod $(D^{\Bbb P_n}_n + 
{\underset\tilde {}\to W_n})"\}$
\sn
\item "{$(ii)$}"   $\bar g^n$ is $D_n$-decreasing. [Saharon: diff]
\ermn
Note that $D_n \in \bold K[S]$,  so 
every initial segment of the play (in which the
second player uses this strategy) belongs to $\bold K[S]$. \nl
By $(B)(iii)$ this is a winning strategy.   \hfill$\square_{\scite{3.9}}$
\enddemo
\bn
Recall all normal filters on ${\Cal Y}/e$ belong to FIL$(e)$.
\mn
\ub{Alternate}:  We split the proof to a series of claims and
definitions.
\bigskip

\definition{\stag{3.9A} Definition}  1) $W_* = \{u \subseteq
\mu:\text{otp}(u) \ge f^*(u \cap w_1)$ and $u \cap \lambda$ is
countable$\}$.
\nl
2) Let $J$ be the following ideal on $Y_0$:

$W \in J$ \ub{iff} for some model $M$ on $\mu$ with countable
vocabulary (with Skolem function) we have

$$
W_* \supseteq W \subseteq \{w \in W_*:w = c \ell_M(w)\}.
$$
\mn
3) For $g \in \dsize \prod_{i < \kappa} (f(i) +1))$ let $\hat g$ be
the function with domain $Y_*$ and $\hat g(A)$ is the $g(i)$-the
member of $A$.
\nl
4) For $W \in J^+$ let proj$(W) = \{A \subseteq w_1:\{w \in W:w \cap
w_1 \notin A\} \in J\}$.
\enddefinition
\bigskip

\demo{\stag{3.9B} Fact}  1) $Y_* \notin J$.
\nl
2) $J$ is a fine normal filter on $W_*$ (and $W_* \notin J$) in fact
the ideal of non-stationary subsets of $W_*$.
\nl
3) $Y_{\bar A} \in J^+$ if $\bar A = \langle A_i:i < 0
\rangle,2^{\aleph_1}$ list the subset of some normal filter $D$ on
$\omega_1$ (see \scite{3.9G}'s proof. \nl
4) If $\bar A',\bar A''$ list the same normal filter on $w_1$ then
$Y_{\bar A'} = Y_{\bar A'}$ mod $J$. 
\nl
5) For $g \in \dsize \prod_{i < \omega} (f^*(i)+1),\hat g$ is well
defined, is a choice function of $Y_*$.
\nl
6) If $g_1 <_D g_2$ then $\hat g_1 \restriction J_D < \hat g_2
\restriction J_D$ mod $J + Y_*$.
\enddemo
\bigskip

\demo{Proof}  1) As $\mu$ is a Ramsey cardinal $> \lambda_0$.
\nl
2) By the definitions.
\nl
3) Easy.
\enddemo
\bigskip

\proclaim{\stag{3.9C} Claim}  Assume $\Bbb Q$ is an $\aleph_1$-complete
forcing notion with $\le \lambda_0$ maximal antichains.
\nl
1) Forcing with $\Bbb Q$ preserves all our assumptions:
\mr
\item "{$(a)$}"  $\mu$ is a Ramsey cardinal$^+$
\sn 
\item "{$(b)$}"  $W_*$ is a family of subsets of $\mu$ such that {\rm
otp}$(w) \ge f(w \cap \omega_1)$ and $J$, defined above, is a fine
normal ideal on $Y_*$ satisfying \scite{3.9B}(3)...then we can forget
$(a)$.
\ermn
2) Forcing with $\Bbb Q$ preserves $``y \in J^+"$ (i.e. if $W \in J^+$
then $\Vdash_{\Bbb Q} ``W \in J^+"$.
\endproclaim
\bigskip

\demo{Proof}  Easy, fill.
\enddemo
\bigskip

\definition{\stag{3.9D} Definition}  Assume $e \in \bold e_{\bold n}$
and $D \in \text{ FIL}_{\bold n}(e)$.
\nl
1) $\Bbb Q = \Bbb Q_e = \{f:f$ is a function with domain a countable
ordinal such that $i \in \text{ Dom}(f) \Rightarrow f(i) \in 
{\Cal Y}^{\bold n}_i\}$.
\nl
2) ${\underset\tilde {}\to f_e}$ is the $\Bbb Q$-name $\cup\{f:f \in
{\underset\tilde {}\to G_{\Bbb Q_e}}\}$.
\nl
3) Let $D/{\underset\tilde {}\to f_e}$ be the $\Bbb Q_e$-name of $\{A 
\subseteq \omega_1$: for every $B \in D$ for stationarily many $i <
\omega_1,{\underset\tilde {}\to f_e}(i) \in B\}$ and
nor$(D,{\underset\tilde {}\to f_e})$ the normal filter which
$D/{\underset\tilde {}\to f_e}$ generates.
\nl
4) For $W \in J^+$ let lift$(W,D) = \{A \subseteq {\Cal Y}/e$ for some
$B \in D: \Vdash_{\Bbb Q_e} ``\{w \in W:{\underset\tilde {}\to f_e}(w
\cap \omega_1) \in B \backslash A \in J"$ (note that we have enough
homogeneity for $\Bbb Q_e$.
\enddefinition
\bigskip

\proclaim{\stag{3.9d1} Claim}  Assume $e \in \bold e_{\bold n}$ and $D
\in { \text{\rm FIL\/}}_{\bold n}(e)$.
\nl
1) $\Vdash_{\Bbb Q} ``\underset\tilde {}\to D/{\underset\tilde {}\to
f_e}$ is a normal filter on $\omega_1"$, (i.e. $w_1 \notin
\underset\tilde {}\to D$).
\nl
2) $|\Bbb Q_e| \le |{\Cal Y}^{\bold n}/e|^{\aleph_0}$ so 
$Z^{|\Bbb Q_e|} \le \lambda_0$ hence $\Bbb Q_e$ has $\le \lambda_0$
maximal antichains; in fact, equality holds as we have demand $|{\Cal
Y}/e| = |\cup\{{\Cal Y}_i:i \in [i_0,\omega_1)\}/e|$ for every $e \in
\bold e$. 
\nl
3) Combine \ scite{3.2A}(4) + \scite{3.9D} - FILL.
\endproclaim
\bigskip

\definition{\stag{3.9E} Definition}  1) We say that ${\frak x} =
(e,D,\bar g,\bar \alpha,f,W)$ is a good position (in the content of
proving \scite{3.9}) if
\mr
\item "{$(a)$}"  $e \in \bold e_{\bold n}$
\sn
\item "{$(b)$}"  $D \in \text{ FIL}_{\bold n}(e)$
\sn
\item "{$(c)$}"  $\bar g = \langle g_\eta:\eta \in u \rangle \in
\text{ Fc}({\Cal T},e)$, so $u=u^{\frak x}$
\sn
\item "{$(d)$}"  $\bar \alpha = \langle \alpha_\eta:\eta \in u
\rangle,\alpha_\eta < \mu$
\sn
\item "{$(e)$}"  $p \in \Bbb Q_e$
\sn
\item "{$(f)$}"  $W = \{w \in W^*:\hat g_\eta(w) = \alpha_\eta$ for
$\eta \in u\} \in J^+$
\sn
\item "{$(g)$}"  $p \Vdash_{\Bbb Q_e} ``W^{\frak x} \cap W_{D,f_e} \in
J^+"$ and proj$(W^{\frak x} \cap W_{D,f_e}) = D$ nor$(D,f_e)$ [FILL].
\endroster
\enddefinition
\bigskip

\demo{\stag{3.9F} Observation}  1) If ${\frak x} = (e,D,\bar g,\bar
\alpha,p,\underset\tilde {}\to W)$ is a good position then
\mr
\item "{$(a)$}"  $\bar \alpha$ is decreasing
\sn  
\item "{$(b)$}"  $D_{\underset\tilde {}\to W}$.
\endroster
\enddemo
\bigskip

\proclaim{\stag{3.9G} Claim}  If $e \in \bold e_{\bold n},D \in 
{ \text{\rm FIL\/}}_{\bold n}(e)$ and $\bar g = \langle g_\eta:
\eta \in u \rangle \in { \text{\rm Fc\/}}({\Cal T},e)$ and $g_\eta \le
f[e]$ for every $\eta \in \text{ Dom}(\bar g)$ \ub{then} we can find a
good position ${\frak x}$ with $\bar g^{\frak x} = e^{\frak x} =
e,\bar g^{\frak x} = g$ and $D \subseteq D^{\frak x}$.
\endproclaim
\bigskip

\demo{Proof}  Let $\bold G \in \Bbb Q_e$ be generic over $\bold V$ and
$f_e = {\underset\tilde {}\to f_e}[G]$.  So in $\bold V[\bold G]$ the
set $W_{D,{\underset\tilde {}\to f_e}[\bold G]}$ belongs to $J^+$ (by
\scite{3.9B}(3)), i.e., let $\langle A^{D_1}_\zeta:\zeta < \zeta^*
\rangle$ list $D_1$ and $W,D,f_e = \{w \in W$: if $\zeta \in w \cap
\zeta^*$ then $f_e(i) = {\underset\tilde {}\to f_e}[\bold G](i) \in
A_\zeta\}$.

Also $\hat g_\eta$ defined in \scite{3.9A}(3) is a choice function on
$W_{D,f_e}$ (see \scite{3.9B}(4)), so as $J$ is a normal ideal and $u$
finite, we can find $\bar \alpha = \langle \alpha_\eta:\eta \in u
\rangle$ such that $W = \{w \in W_{D,f_e}:\hat g_\eta(w) =
\alpha_\eta$ for $\eta \in u\}$ belongs to $J^+$.  As all this holds
in $\bold V[\bold G]$.  So $\bar \alpha$ there is a condition $p \in
\Bbb Q_e$ which forces this, and we are done.
\enddemo
\bigskip

\proclaim{\stag{3.9H} Claim}  Assume that
\mr
\item "{$(a)$}"  ${\frak x}_1 = (e_1,D_1,\bar g_1,\bar
\alpha_1,p,{\underset\tilde {}\to W_1})$ is a good position
\sn
\item "{$(b)$}"  $\bar g_2 = \langle g^2_\eta:\eta \in u_2 \rangle \in
{ \text{\rm Fc\/}}({\Cal T},\bold n)$ and $\bar g_2 \restriction u_1 =
\bar g_2$
\sn
\item "{$(c)$}"  $e_1 \le e_2$ in $\bold e_n$ and $D_2 \in \text{
FIL}_{\bold n}(e_2)$ or just ${\Cal A} \subseteq {\Cal P}({\Cal
Y}_{\bold n}/e_2),{\Cal A} = \{A_\zeta:\zeta < \zeta^*\}$
\sn
\item "{$(d)$}"  $p_1 \Vdash_{\Bbb Q_{e_1}} ``\{w \in {\underset\tilde
{}\to W_1}:{\Cal Y}_{w \cap w_1} \nsubseteq \cup\{A_\zeta:\zeta \in
\zeta^* \cap w\}\}$ does not belong to $J^{\bold V[\Bbb Q_{e_1}]}"$.
\ermn
\ub{Then} we can find a good position ${\frak x}_2$ such that
$e^{{\frak x}_2} = e_2,\bar g^{{\frak x}_2} = \bar g^2$ and $D_2
\subseteq D^{{\frak x}_2}$.
\endproclaim
\bigskip

\demo{Proof}  Let $\bold G$ be a subset of $\Bbb Q_{e_1[{\frak x}_1]}$
generic over $\bold V$ such that $p^{{\frak x}_1} \in \bold G_1$.  Now
$\Bbb Q_{e_2}$ is an $\aleph_1$-complete forcing of cardinality $\le
|{\Cal Y}_{\bold n}/e_2|^{\aleph_0} \le \lambda_0$ and $\Bbb Q_{e_1}$
is $\aleph_1$-complete $|\Bbb Q_{e_1}| \le |{\Cal Y}_{\bold
n}/e_1|^{\aleph_0} \le |{\Cal Y}_{\bold n}/e_2|^{\aleph_0} \le
\lambda_0$, so $\Bbb Q_{e_2}$ satisfies the same conditions in $\bold
V[\bold G_1]$ (if $\lambda_0$ is no longer a cardinal it does not
matter).

Note that by assumption (c)
\mr
\item "{$\circledast$}"  in $\bold V[\bold G_1],\Bbb Q_{e_2} \Vdash$
``the set $\{{\underset\tilde {}\to W^1_2} =: \{w \in {\underset\tilde
{}\to W_1}[\bold G_1]$: the set $(({\underset\tilde {}\to
f_{e_1}}[\bold G_1])(w \cap \omega_1))^{[e_2]} \in {\Cal Y}_{w \cap
\omega_1}/e_2$ is not included in $\cup\{A_\zeta:\zeta \in w\}\}$ is
stationary (i.e. $\notin J$)".
\ermn
We continue as in the previous claim.
\enddemo
\bigskip

\proclaim{\stag{3.9J} Claim}  If clauses (a) + (b) of \scite{3.9G}
holds, \ub{then} a sufficient condition for clause (c) is

(c)' $\quad$ FILL.
\endproclaim
\bigskip

\demo{\stag{3.9I} Proof of \scite{3.9}}  During the play, the player II
chooses also a good position ${\frak x}_n$ and maintains $\bar
g^{{\frak x}_n} = \bar g_n,\bar \alpha^{{\frak x}_n} = \bar \alpha$.
\enddemo
\bigskip

\remark{\stag{3.9K} Remark}   1) From the proof, 
instead $\bold K[S] \models ``\lambda$ is 
Ramsey",  $\bold K[S] \models  ``\mu  \rightarrow (\alpha)^{<\omega
}_{\lambda_0}$ for $\alpha < \lambda_0"$ is enough for showing for \scite{3.9}.
\nl
2)  Also if $\dsize \prod_{i < \omega_1}(|f(i)| + 1) < \mu _0,
[\alpha  < \mu _0 \Rightarrow  
|\alpha |^{\aleph _0} < \mu _0]$,  it is enough:  
$S \subseteq \alpha < \mu_0 \Rightarrow$ in  
$\bold K[S]$ there is $\mu \rightarrow (\alpha)^{<\omega}_2$.
\endremark
\bigskip

\proclaim{\stag{3.10} Theorem}  Assume $\bold n$ is a $\kappa$-niceness context.
Let $D^\ast \in \text{ FIL}(e,{\Cal Y})$ be a normal ideal 
on ${\Cal Y}_{\bold n}/e$.  If for every $f:{\Cal Y} \rightarrow  
(\sup\{{\text{\rm Suc\/}}(D'):D' \in { \text{\rm FIL\/}}_{\bold
n}\})^+$ supported by some $e \in \bold e_{\bold n}$.
$D^\ast_{\bold n}$ is nice to $f$, \ub{then} for every 
$f \in {}^\kappa${\rm Ord}, $\bold n$ is nice to $f$.
\endproclaim
\bigskip

\demo{Proof}  By determinacy of the games (and the LS argument).
\enddemo
\bigskip

\remark{\stag{3.10A} Remark}  0) The value $|\text{FIL}_{\bold e}|$
really should be an upper bound.
\nl
1)  So, the existence of $\mu,\mu \rightarrow  
(\alpha)^{<\omega}_{\aleph_0}$ for every $\alpha < 
(\dsize \sum_{\chi < \mu} \chi^\kappa)^+$,  
is enough for $``D^\ast$ is nice". 
\nl
2) If there 
is a nice $D$'s in the plays from \scite{3.4}, the second player winning 
strategy can be chosen such that all subsequent filters are nice: just by 
renaming have $g_{<>}$ constant large enough. [Saharon: diff]
\endremark
\bigskip

\proclaim{\stag{3.10B} Claim}  In claim \scite{3.9} we can omit
``$\kappa_{\bold n} = \aleph_1"$.
\endproclaim
\bigskip

\demo{Proof}  Let $\Bbb P = \text{ Levy}(\aleph_0,\kappa_{\bold n})$.
Now
\mr
\item "{$(*)$}"  also in $\bold V^{\Bbb P}$ the object $\bold n$ is a
successor content, if we do not distinguish between $D \in \text{
FIL}_{\bold n}$ and $\{A \in \bold V^{\Bbb P}:A \subseteq {\Cal
Y}/e(D)$ and $(\exists B \in D)(B \subseteq A)\}$.
\endroster
\enddemo
\bigskip

\demo{\stag{3.11} Conclusion}:  Let  $\lambda_0 = 
(\sup\{|\text{Suc}_{\bold n}(D')|:D' \in \text{ FIL}_{\bold n}\})^+
\cup \{2^{|{\Cal Y}/e|^{< \kappa}}:e \in \bold e_{\bold n}\})^+,
\mu^\ast \ge  \aleph_2$; if for every $S \subseteq \lambda_0$ 
there is a Ramsey cardinal in $\bold K[S]$ above $\lambda_0$ 
then $\bold n$ is nice.
\enddemo
\bigskip

\demo{Proof}  By \scite{3.9}, \scite{3.10}.
\enddemo
\bigskip

\remark{\stag{3.12} Concluding Remark}  1) We could have 
used other forcing notions, not Levy$(\kappa,|{\Cal Y}/e_n|)$.  
E.q., if $\kappa = \aleph_1,\mu  = \kappa^+$ we could use 
finite iterations of the forcing of Baumgartner to add a club of  $\omega _1$, 
by finite conditions.  (So this forcing notion has cardinality  $\aleph _1)$.  
Then in \scite{3.9} we can weaken the demands on
$\lambda_0: \lambda_0 = \dsize \sum_{\chi < \mu_0} 2^\chi + \dsize
\prod_{i < \omega_1}|1 + f(i)| + |\bold e|$,  hence also in 
\scite{3.11},  $\lambda _0 = \dsize \sum_{\chi < \mu_\ast} 2^\chi$ is O.K. \nl
2)  Concerning  $|\bold e|$  remember \scite{3.6}(5),(6).
\nl
3)  Similarly to (1).  If $\theta < \mu \Rightarrow \text{ cov}(\theta,\aleph_1,\aleph_1,2) < \mu$ then by \scite{2.6} we can use forcing
notions of Todorcevic for collapsing  $\theta  < \mu $  which has cardinality  
$< \mu$. \nl
4)  If we want to have $\lambda_0 =: \dsize \prod_{i < \omega_1}|f(i)
+ 2|$  (or even $T_D(f + 2))$,  we can get this by weakening
further the first player letting him choose only  $A_n$ which are easily 
definable from the  $\bar g^{n-1}$,  we shall return to it in a subsequent 
paper.
\endremark
\newpage

\head {\S4 Ranks} \endhead  \resetall \sectno=4
 \spuriousreset
\bigskip

\demo{\stag{4.1} Convention}  1) Like \scite{3.1} and:
\nl
2)  $\bar g^\ast \in F_c({\Cal T},e^\ast,{\Cal Y}),\eta^* \in 
\text{ Dom}(\bar g^\ast),\nu^\ast$ an immediate successor of 
$\eta^\ast$ not in  Dom $g^\ast $,  
$D^\ast  \in \text{ FIL}(e^\ast,{\Cal Y})$  is such that in  
$\Game^{\bar \gamma^\ast}(D^\ast ,\bar g^\ast ,e^\ast )$  second player wins  
(all constant for this section).  FIL$^\ast(e)$  will be the set of  
$D \in \text{ FIL}(e,{\Cal Y})$  such that $e \ge e^\ast$,  $(D^\ast )^{[e]} 
\subseteq D$ and in 
$\Game^{\bar \gamma^\ast}(D^\ast,\bar g^\ast,e^\ast)$  second 
player wins.  (So actually FIL$(e^\ast,{\Cal Y})$  depends on  
$D^\ast ,\bar g^\ast ,e^\ast $,  too).
\enddemo
\bigskip

\definition{\stag{4.2} Definition}   1) rk$^5_D(f)$ for 
$D \in \text{ FIL}^\ast(e,{\Cal Y})$,  $f \in {}^{{\Cal
Y}/e}\text{Ord},
f <_D \bar g^\ast _{\eta ^\ast }$ will be: the minimal 
ordinal  $\alpha $  such that for some  $D_1$, $e_1$, $\bar \gamma ^1$ we have 
$D^{[e_1]} \subseteq D_1 \in \text{ FIL}(e_1,{\Cal Y})$,  $\bar \gamma^1 = 
\bar \gamma^\ast \char 94 \langle \nu^\ast,\alpha \rangle$ (i.e.  
dom$(\bar \gamma^1) = (\text{dom}(\bar \gamma^\ast)) \cup  \{\nu ^\ast \}$,  
$\bar \gamma^1 \restriction \text{ dom}(\bar \gamma^\ast) 
= \bar \gamma ^\ast $,  
$\gamma^1_{\nu^\ast} = \alpha)$ and in $\Game^{\bar \gamma ^1}(D,\bar g^\ast 
\char 94  <\nu^\ast,f>)$  second player wins and  $\infty $  if there is no 
such $\alpha$. \nl
2) rk$^4_D(f)$ is $\sup\{\text{rk}^5_{D+A}(f):A \in  D^+\}.$
\enddefinition
\bigskip

\proclaim{\stag{4.3} Claim}  1) {\rm rk}$^5_D(f)$ is (under the
circumstances of \scite{4.1}, \scite{4.2}) 
an ordinal $< \gamma^\ast_{\eta^\ast }.$
\nl
2) {\rm rk}$^4_D(f)$  is an ordinal  $\leq\gamma^\ast _{\eta ^\ast }.$
\endproclaim
\bigskip

\proclaim{\stag{4.4} Claim}   If $D \in { \text{\rm FIL\/}}^*(e,{\Cal
Y}),h <_D f <_D g^\ast_{\eta^\ast}$ then {\rm rk}$^5_D(h) < 
{ \text{\rm rk\/}}^5_D(f).$
\endproclaim
\bigskip

\demo{Proof}   Let $e_1,D_1$ witness rk$^5_D(f) = \alpha$ so $e(D)\le  
e_1$,  $D \subseteq D_1 \in \text{ FIL}^*(e_1)$  
and in $G^{\bar \gamma \char 94 < \nu^\ast,\alpha >}(D_1,\bar g^\ast
\char94 < \nu^*,f>,e)$  second player wins.  
We play for the first player:  $e = e_1$,  
$A_0 = {\Cal Y}/e_1$,  $\bar g^0 = \bar g^\ast \char 94 \langle
\nu^\ast,f \rangle \char 94 \langle \nu^\ast \char 94 <0>,g \rangle$,
now the first player should be able to answer say  
$e_2$,  $D_2$,  $\bar \gamma ^2$.  So  $\gamma^2_{\nu^\ast \char 94 <0>} < 
\gamma^2_{\nu^\ast} = \alpha$,  and by \scite{3.7}(3), we know that in  
$G^{\bar \gamma^{'}}(D_2,\bar g^\ast \char 94 < \nu^\ast,g>,e_2)$  where  
$\bar \gamma^{'} = \bar \gamma^\char 94 \langle \nu^\ast,
\gamma^2_{\nu ^\ast \char94 <0>} \rangle$,  
second player wins.  \nl
${{}}$     \hfill$\square_{\scite{4.4}}$
\enddemo
\bigskip

\proclaim{\stag{4.5} Claim}  Let $e \ge e^*,
D \in { \text{\rm FIL\/}}^*(e,{\Cal Y})$.
\nl
1) For $e \ge e(D),A \in (D^{[e]^+},f \in {}^{{\Cal Y}/e}${\rm Ord},
$f <_D g^\ast_{\eta^\ast}$ we have: 

$$
{\text{\rm rk\/}}^5_D(f) \le 
{ \text{\rm rk\/}}^5_{D^{[e]}+A}(f) \le { \text{\rm
rk\/}}^4_{D^{[e]}+A}(f) \le { \text{\rm rk\/}}^4_D(f).
$$
\mn
2) If $e_2 \ge e_1 \ge e(D),f_\ell \in {}^{\Cal Y}${\rm Ord} is supported 
by $e_\ell$,  $f_1 \le_D f_2 <_D g^\ast_{\eta ^\ast }$ \ub{then}  
{\rm rk}$^\ell_D(f_1) \le { \text{\rm rk\/}}^\ell_D(f_2)$ for $\ell = 4,5$.
\endproclaim
\bigskip

\demo{Proof}  Left to the reader.
\enddemo
\newpage

\head {\S5  More on Ranks and Higher Objects} \endhead  \resetall \sectno=5
 \spuriousreset
\bigskip

\demo{\stag{5.1} Convention}
\mr
\item "{$(a)$}"   $\mu^\ast$ is a cardinal  $> \aleph_1$ (using $\aleph_1$ rather than 
an uncountable regular  $\kappa $  is to save parameters)
\sn
\item "{$(b)$}"   ${\Cal Y}$ a set of cardinality  $\dsize
\sum_{\kappa < \mu_\ast} \kappa$
\sn
\item "{$(c)$}"   $\iota$ a function from  ${\Cal Y}$  onto  $\omega _1$,  
$|\iota^{-1}(\{\alpha \})| = |{\Cal Y}|$  for  $\alpha  < \omega ,$
\sn
\item "{$(d)$}"  Eq the set of equivalence relation  $e$  on  ${\Cal Y}$ such that: 
{\roster
\itemitem{ $(\alpha)$ }   $y e z \Rightarrow \iota (y) = \iota (z)$ 
\sn
\itemitem{ $(\beta)$ }  each equivalence class has cardinality  $|{\Cal Y}|$ 
\sn
\itemitem{ $(\gamma)$ }  $e$ has $< \mu^\ast$ equivalence classes
\endroster}
\item "{$(e)$}"   $D$  denotes a normal filter on some  ${\Cal Y}/e (e
\in \text{ Eq})$,  we write  $e = e(D)$.  The set of such  $D$'s is FIL$({\Cal Y})$.
\sn
\item "{$(f)$}"  $E$ denotes a set of $D$'s as above, such that: 
{\roster
\itemitem{ $(\alpha)$ }  for some $D = \text{ Min } E \in E$    
$(\forall D')[D' \in  E \Rightarrow  (e,D) \le (e(D'),D')]$ 
\sn
\itemitem{ $(\beta)$ }  if  $D \in  E$,  $A \subseteq {\Cal Y}/e_1,e_1
\ge e(D)$, $A \ne \emptyset$ mod $D$  then     
$D^{[e_1]} + A \in E$
\endroster}
\item "{$(g)$}"   $E^{[e]} =: \{D \in  E:e(D) = e\}$
\sn
\item "{$(h)$}"   ${\Cal E}$ denotes a set of $E$'s as above, such that: 
{\roster
\itemitem{ $(\alpha)$ }  there is $E = \text{ Min }{\Cal E} 
\in {\Cal E}$  satisfying     
$(\forall E')(E' \in  E \Rightarrow  E' \subseteq E)$ 
\sn
\itemitem{ $(\beta)$ }  if  $D \in E \in {\Cal E}$ then    
$E_{[D]} = \{D':D' \in E$ and $(e(D),D) \le (e(D'),D')\} \in {\Cal
E}$.
\endroster}
\endroster
\enddemo
\bigskip

\definition{\stag{5.2} Definition}   1) We say $E$ is $\lambda$-divisible when: 
for every  $D \in  E$,  and  $Z$,  a set of cardinality  $< \lambda $ 
there is $D$'s such that:  
\mr
\item "{$(\alpha)$}"    $D' \in  E$  
\sn
\item "{$(\beta)$}"   $(e(D),D) \le (e(D'),D')$  
\sn
\item "{$(\gamma)$}"   $\bold j:{\Cal Y}/e(D') \rightarrow  Z$  
\sn
\item "{$(\delta)$}"   for every function $h:{\Cal Y}/e(D) \rightarrow Z$  we have     
$\{y/e(D'):h(y/e(D)) =  (y/e(D'))\} \ne \emptyset$ mod $D'$.
\ermn
2)  We say $E$ has $\lambda$-sums when: 
for every  $D \in  E \in  {\Cal E}$  and sequence  $\langle Z_\zeta :\zeta  < 
\zeta^\ast < \lambda \rangle$ of subsets of ${\Cal Y}/e(D)$ 
there is  $Z^\ast \subseteq {\Cal Y}/(e/(D)$,  such that:  
$Z^\ast \cap Z_\zeta = \emptyset$ mod $D$ and:  
[if $(e(D),D) \le (e',D'),e' = e(D'),D' \in E_{[D]}$ and   
$\dsize \bigwedge_\zeta Z^{[e']}_\zeta = \emptyset$ mod $D'$ then 
$Z^\ast \in D'$]. \nl
3)  We say  $E$ has weak $\lambda$-sum if for every $D \in E(\in {\Cal
E})$ and sequence $\langle Z_\zeta:\zeta < \zeta^\ast < \lambda
\rangle$  of subsets of  ${\Cal Y}/e(D)$  there is  $D^\ast $,  
$D^\ast  \in  E_{[D]}$ such that: 
\mr
\item "{$(\alpha)$}"   if $(e(D),D) \le 
(e',D'),D' \in E_{[D]}$ and  $Z_\zeta  = 
\emptyset$ mod $D'$  for $\zeta < \zeta^\ast$ and $e(D^\ast) \le
e(D')$  then  $D^\ast  \subseteq  D'$  (more exactly  $D^{\ast ^{[\ast ]}} 
\subseteq  D^{[\ast ]}$ and) 
\sn
\item "{$(\beta)$}"   $Z_\zeta = \emptyset$ mod $D^\ast$ for $\zeta <
\zeta^\ast$.
\ermn
4)  If $\lambda = \mu^\ast$ we omit it.  We say ${\Cal E}$  is 
$\lambda$-divisible if every  $E \in {\Cal E}$ has.  We say  ${\Cal E}$  has
weak $\lambda$-sums if: [rest diff] for every  $E \in {\Cal E}$  and sequence  
$\langle Z_\zeta:\zeta < \zeta^\ast < \lambda \rangle$ of subsets of  
${\Cal Y}/e(E)$ there is $E^\ast$, $E^\ast \in {\Cal E}_{[E]}$ such 
that:
\mr
\item "{$(\alpha)$}"   if $(e(E),E) \le (e',E')$,  $E' \in  {\Cal E}$
and $Z_\zeta = \emptyset$ mod Min$(E')$ for $\zeta < \zeta^\ast$ and $e(E^\ast ) 
\le e(E')$ then $E^\ast \subseteq E'$  
\sn
\item "{$(\beta)$}"  $Z_\zeta = \emptyset$ mod Min$(E^\ast)$ for $\zeta  < 
\zeta^\ast$.
\ermn
We now define variants of the games from \S3.
\enddefinition
\bigskip

\definition{\stag{5.3} Definition}  For a given ${\Cal E}$, for every
$E \in {\Cal E}$:
\nl
1)  We define a game $G^\ast_2(E,\bar g)$.
In the $n-th$ move first player chooses $D_n \in E_{n-1}$ (stipulating  
$E_{-1} = E)$  and choose  $\bar g_n \in  F_c({}^\omega \omega ,e(D_n),{\Cal Y})$
extending  $\bar g_{n-1}$ (stipulating  $\bar g_{-1} = \bar g)$  such that  
$\bar g_n$ is $D_n$-decreasing.  \ub{Then} the second player chooses  $E_n$,  
$(E_{n-1})_{[D_n]} \subseteq E_n \in {\Cal E}$. 

In the end the second player wins if $\dbcu_{n < \omega} \text{ Dom } 
\bar g_n$ has no infinite branch.
\nl
2)  We define a game $G^{\bar \gamma}_2(E,\bar g)$  
where Dom$(\bar \gamma) = \text{ Dom}(\bar g)$, each  
$\gamma _\eta $ an ordinal,  $[\eta \triangleleft \nu  
\Rightarrow  \gamma _\eta  > \gamma _\nu ]$  similarly to  
$G^\ast _2(D,\bar g)$  but the second player in addition chooses an indexed set
$\bar \gamma_n$ of ordinals,  Dom$(\bar \gamma_n) = 
\text{ Dom}(\bar g_n)$,  
$\bar \gamma_n \restriction \text{ Dom}(\bar \gamma_{n-1}) = \bar
\gamma_{n-1}$ and $[\eta \triangleleft \nu 
\Rightarrow  \gamma _{n,\eta } > \gamma _{n,\nu }].$
\enddefinition
\bigskip

\definition{\stag{5.4} Definition}   1) We say ${\Cal E}$ is nice to 
$\bar g \in F_c({\Cal T},e,{\Cal Y})$  if for every  $E \in  {\Cal E}$  with  $e 
\le e(E)$ the second player wins the game $\Game^\ast _2(E,\bar g).$
\nl
2)  We say ${\Cal E}$  is nice if it is nice to  $\bar g$  whenever  $E \in  
{\Cal E}$,  $e \le e(E)$,  $\bar g \in F_c({\Cal T},e)$,  
$\bar g$  is (Min $E$)-decreasing, we have:  ${\Cal E}_{[E]}$ is nice to  
$\bar g.$
\nl
3)  If Dom$(\bar g) = \{<>\}$ we write $g_{<>}$ instead $\bar g$.
\nl
4)  We say  ${\Cal E}$ is nice to $\alpha$ if it is nice to the constant 
function  $\alpha$.
\enddefinition
\bigskip

\proclaim{\stag{5.5} Claim}  1) If ${\Cal E}$ is nice to  $f$,  $f \in  
F_c({\Cal T},e,{\Cal Y}),g \in F_c({\Cal T},e,{\Cal Y})$,  
$g \le f$ \ub{then} ${\Cal E}$  is nice to  $f.$
\nl
2)  The games from 
\scite{5.4} are determined, and the winning side has winning 
strategy which does not need memory.
\nl
3) The second player wins $G^\ast_2(E,\bar g)$ iff for some  $\bar \gamma $
second player wins $G^{\bar \gamma}_2(E,g).$
\nl
4)  If the second player wins  $G^\gamma _2(E,f)$,  $\bar g \in  
F_c({\Cal T},e(E))g_\eta \le f$  for  $\eta \in {\text{\rm Dom\/}}(\bar g)$  
then the second player wins in  $G^{\bar \gamma }_2(E,\bar g)$ when we let

$$
\gamma_\eta = \gamma + \,{\text{\rm max\/}}\{(\ell g(\nu) - \ell g(\eta) +
1):\nu \text{ satisfies } \eta \trianglelefteq \nu \in 
{ \text{\rm Dom\/}}(\bar g)\}.
$$
\endproclaim
\bigskip

\proclaim{\stag{5.6} Lemma}   Suppose 
$f_0 \in {}^{({\Cal Y}/e)}\text{Ord},e \in \text{ Eq}$ and $\lambda_0 =: 
\sup\{\dsize \prod_{x \in Y} {\Cal Y}_e (f^{[e]}_0(x)+1:e$ satisfies
$e_0 \le e \in \bold e\}$. \nl
1)  If there is a Ramsey cardinal  
$\ge \cup \{f(x)+1:x \in { \text{\rm Dom\/}}(f_0)\}$ \ub{then} 
there is a $\mu^\ast$-divisible ${\Cal E}$ nice to $f_0$ having weak 
$\mu^\ast$-sums.
\nl
2)  If for every  $A \subseteq  \lambda _0$ there is in $\bold K[A_0]$
a Ramsey cardinal $> \lambda_0$, \ub{then} 
there is a $\mu ^\ast $-divisible  ${\Cal E}$  
which has weak $\mu ^\ast $-sums and is nice to  $f.$
\nl
3)  In part 2 if $\lambda_0 = 2^{<\mu_0}$ \ub{then} there is a 
$\mu^\ast$-divisible nice ${\Cal E}$  which has weak $\mu ^\ast $-sums.
\endproclaim
\bigskip

\remark{\stag{5.6A} Remark}   This enables us to 
pass from $``\text{pp}_{\Gamma(\theta,\aleph_1)}$ large" to  
``pp$_{\text{normal}}$ is large".
\endremark
\bigskip

\demo{Proof}  1) Define  $f_1 \in {}^{(\aleph_1)}\text{Ord},f_1(i) = 
\sup \{f_0(y/e):\iota (y) = i\}$,  let  $\lambda $  be such that:  $\lambda  
\rightarrow (\sup\{f_1(i))^{<\omega }_2:i < \aleph_1\}$ (or 
just $\emptyset \notin D^\ast_n$ - see below)  let  
$\lambda_n = (\lambda^{\mu^\ast})^{+n}$,  

$$
I_n = \{s:s \subseteq \lambda_n,s \cap \omega_1 \text{ a countable
ordinal}\}
$$  

$$
J_n = \{s \in I_n:s \cap \lambda \text{ has order type } \ge f_0(s \cap  
\omega _1)\}.
$$ 
\mn
Let $D^\ast_n$ be the minimal fine normal filter on  $J_n.$ 

Let for  $n < \omega$ and $e \in \text{ Eq}$,
$H_{n,e} = \{h:h$ a function from  $J_n$ into ${\Cal Y}/e$  such that  
$\iota(h(s)) = s \cap \omega_1\}$. 

Let $\Bbb P_n = \{p:p \subseteq J_n,p \ne \emptyset$ mod $D^\ast_n\},
\Bbb P = \dbcu_{n < \omega} \Bbb P_n$ and for $p \in \Bbb P$ let  
$n(p)$  be the unique  $n$  such that $p \in \Bbb P_n.$ 

Let $p \le q$ (in $\Bbb P$) if 
$n(p) \le n(q)$ and $\{s \cap \lambda_{n(p)}:s \in  q\} 
\subseteq p$. \nl
Now for every 
$e \in \text{ Eq}$,  $n < \omega $,  $p \in  P_n$,  $h \in  H_{n,e}$ we
let:

$$
D^{n,e,h}_p = \{ A \subseteq {\Cal Y}/e:h^{-1}(A) 
\supseteq p \text{ mod } D^\ast_{n(p)}\}
$$

$$
E^{n,e,h}_p = \{ D^{n^1,e^1,h^1}_q:p \le q \in P,n^1 = 
n(q) \text{ and } (n^1,e^1,h^1) \ge (n,e,h)\}
$$
\mn
where $(n^1,e^1,h^1) \ge (n,e,h)$ means: $n \le n^1 < \omega $,  $e 
\le e^1 \in \text{ Eq}$,  $h^1 \in  H_{n^1,e^1}$ and for  $s \in  J_{(n^1)}$,  
$h^1(s)^{[e]} = h(s \cap \lambda_n)$ and we define  $(p^1,n^1,e^1,h^1) \ge 
(p,n,e,h)$  similarly.  Let

$$
{\Cal E}^{n,e,h}_p = \{ E^{n^1,e^1,h^1}_q:p \leq q \in P,n^1 = 
n(q),(n^1,e^1,h^1) \ge (n,e,h)\}.
$$
\mn
Note:  $(p^1,n^1,e^1,h^1) \ge (p,e,n,h)$ implies  
$D^{n^1,e^1,h^1}_{p^1} \supseteq D^{n,e,h}_p$,
$E^{n^1,e^1,h^1}_{p^1} \subseteq  E^{n,e,h}_p$ and  
${\Cal E}^{n^1,e^1,h^1}_{p^1} \subseteq  {\Cal E}^{n,e,h}_p$.  Now any  
${\Cal E} = {\Cal E}^{n,e,h}_p (p \in P)$ is as required. 

A new point is  $``{\Cal E}$ is $\mu^\ast$-divisible".  So suppose  $E \in  
{\Cal E} = {\Cal E}^{n,e,h}_p$ so $E = E^{n^1,e^1,h^1}_q$ for some  
$(q,n^1,e^1,h^1) \ge (p,n,e,h)$.  Let  $Z$  be a set of cardinality  $< 
\mu^\ast$, so $(\lambda_{n^1})^{|Z|} = \lambda _{n_1}$;  let 
$\{ h_\zeta:\zeta < \zeta^\ast = |{\Cal Y}/e_1|^{|Z|} \le  
2^\mu \le \lambda_{n^1}\}$ list all function  $h$  from  
${\Cal Y}/e_1$ to  $Z$.  Let  $\langle S_\zeta :\zeta  < 
|{\Cal Y}/e_1|^{|Z|}\rangle $  list a sequence of pairwise disjoint stationary 
subsets of  $\{\delta < \lambda_{n^1+1}:\text{cf}(\delta) =
\aleph_0\}$.  Let $e_2 \in \text{ Eq}$ be such that $e_1 \le e_2$ 
and for every $y \in {\Cal Y}$,  
$\{z/e_2:ze_1y\} = \{ x(y/e,t):t \in Z \}$,  we let  
$q_2$,  $q \le q_2 \in P$ be:  $q_2 = \{ s \in J_{n^1+1}:s \cap 
\lambda_{n^1} \in q$ and sup $s \in \dbcu_\zeta S_\zeta\}$,  
lastly we define  $h^2:J_{n^1+1} \rightarrow  
{\Cal Y}/e_1$ by:  $h^2(s) = x(h^1(s \cap \lambda_{n^1}),h_\zeta(s \cap \lambda_{n^1}))$  
if $s \in q_2$, sup $s \in S_\zeta$ (for  $s \in J_{n^1+1} \backslash
q_2$ it does not matter).
The proof that  $q_2$,  $e_2$,  $h^2$ are as required is as in \cite{RuSh:117} and 
more specifically \cite{Sh:212}.
As for proving  $``{\Cal E}^{n,e,h}_p$ has weak $\mu^\ast$-sums" the point is
that the family of fine normal filters on  $\mu $  has $\mu^\ast$-sum.
\nl
2)  Similar to \scite{3.9}(and \scite{3.6}(5),(6)).
\nl
3)  Similar to \cite[1.7]{Sh:386}.  \hfill$\square_{\scite{5.6}}$
\enddemo
\newpage

\head {\S6  Hypotheses: Weakening of GCH} \endhead  \resetall \sectno=6
 \spuriousreset
\bigskip

We define some hypotheses; except the first we do not know now whether their 
negations are consistent with  ZFC.
\definition{\stag{6.1} Definition}  We define a series of hypothesis:  
\mn
$(A)$  pp$(\lambda) = \lambda^+$ for every singular  $\lambda .$
\nl
$(B)$ If ${\frak a}$ is a set of regular cardinals, $|{\frak a}| 
< \text{ Min}({\frak a})$ \ub{then} $|\text{pcf}({\frak a})| \le 
|{\frak a}|$.
\nl
$(C)$  If ${\frak a}$ is a set of regular cardinals, $|{\frak a}| 
< \text{ Min}({\frak a})$ \ub{then} pcf$({\frak a})$ 
has no accumulation point which is inaccessible 
(i.e. $\lambda$ inaccessible $\Rightarrow  \sup 
(\lambda \cap \text{ pcf}({\frak a}) < \lambda).$
\nl
$(D)$   For every $\lambda$, $\{\mu < \lambda:\mu$ singular and  
pp$(\mu) \ge \lambda\}$  is countable.
\nl
$(E)$  For every $\lambda$, $\{\mu < \lambda:\mu$ singular and cf$(\mu) =
\aleph_0$ and pp$(\mu) \ge  \lambda\}$  is countable.
\nl
$(F)$  For every  $\lambda,$ $\{\mu <\lambda:\mu$ singular of uncountable 
cofinality, pp$_{\Gamma(\text{cf}(\mu))}(\mu) \ge \lambda\}$  is finite.
\nl
$(D)_{\theta,\sigma,\kappa}$   For every $\lambda $,  $\{\mu  < \lambda :\mu 
> \text{ cf}(\mu) \in [\sigma,\theta)$ and pp$_{\Gamma (\theta ,\sigma )}(\mu )
\ge \lambda\}$ has cardinality  $< \kappa .$
\nl
$(A)_\Gamma$    If $\mu > \text{ cf}(\mu)$ then pp$_\Gamma (\mu) = \mu^+$ (or in 
the definition of pp$_\Gamma(\mu)$ the supremum is on the empty set).
\nl
$(B)_\Gamma,(C)_\Gamma$   Similar versions (i.e. use pcf$_\Gamma$).
\enddefinition
\bigskip

We concentrate on the parameter free case.
\proclaim{\stag{6.2} Claim}:  In \scite{6.1}, we have:
\mr
\item   $(A) \Rightarrow  (B) \Rightarrow  (C)$
\sn
\item   $(A) \Rightarrow  (D) \Rightarrow  (E)$,  $(A) \Rightarrow  (F)$
\sn
\item   $(E) + (F) \Rightarrow  (D) \Rightarrow  (B)$.  
[Last implication --- by the localization theorem \cite[\S2]{Sh:371}]
\sn
\item  if $(\forall \mu)(\mu > \text{ cf}(\mu) = \aleph_0$ the
hypothesis (A) of \scite{6.1} holds. \nl
[Why?  By \cite[xx]{Sh:g}.]
\endroster
\endproclaim
\bigskip

\proclaim{\stag{6.3} Theorem}   Assume Hypothesis \scite{6.1}(A).
\nl
1)  For every  $\lambda > \kappa$,

$$
{\text{\rm cov\/}}(\lambda,\kappa^+,\kappa^+,2) = 
\cases \lambda^+ &{\text{\rm if cf\/}}(\lambda) \le \kappa \\
  \lambda &{\text{\rm if cf\/}}(\lambda) > \kappa. \endcases
$$
\mn
2)  For every  $\lambda > \kappa = { \text{\rm cf\/}}(\kappa) 
> \aleph_0$,  there is a 
stationary $S \subseteq [\lambda]^{\le \kappa},|S| = 
\lambda^+$ if ${\text{\rm cf\/}}(\lambda) \le \kappa$ and 
$|S| = \lambda$ if ${\text{\rm cf\/}}(\lambda) > \kappa$.
\nl
3)  For $\mu$ singular, there is a tree with ${\text{\rm cf\/}}(\mu)$ 
levels each level of cardinality $< \mu$, and with 
$\ge \mu^+ ({\text{\rm cf\/}}(\mu))$-branches.
\nl
4)  If $\kappa \le { \text{\rm cf\/}}(\mu) < \mu 
\le 2^\kappa$ \ub{then} there is an 
entangled linear order  ${\Cal T}$ of cardinality $\mu^+$.
\endproclaim
\bigskip

\demo{Proof}  1) By \cite[\S1]{Sh:400}.
\nl
2)  By part (1) and \scite{2.6}.
\nl
3, 4)  By \cite[\S4]{Sh:355}.
\enddemo
\newpage


\nocite{ignore-this-bibtex-warning} 
\newpage
    
REFERENCES.  
\bibliographystyle{lit-plain}
\bibliography{lista,listb,listx,listf,liste}

\def\germ{\frak} \def\scr{\cal} \ifx\documentclass\undefinedcs
  \def\bf{\fam\bffam\tenbf}\def\rm{\fam0\tenrm}\fi 
  \def\defaultdefine#1#2{\expandafter\ifx\csname#1\endcsname\relax
  \expandafter\def\csname#1\endcsname{#2}\fi} \defaultdefine{Bbb}{\bf}
  \defaultdefine{frak}{\bf} \defaultdefine{=}{\B} 
  \defaultdefine{mathfrak}{\frak} \defaultdefine{mathbb}{\bf}
  \defaultdefine{mathcal}{\cal}
  \defaultdefine{beth}{BETH}\defaultdefine{cal}{\bf} \def\bbfI{{\Bbb I}}
  \def\mbox{\hbox} \def\text{\hbox} \def\om{\omega} \def\Cal#1{{\bf #1}}
  \def\pcf{pcf} \defaultdefine{cf}{cf} \defaultdefine{reals}{{\Bbb R}}
  \defaultdefine{real}{{\Bbb R}} \def\restriction{{|}} \def\club{CLUB}
  \def\w{\omega} \def\exist{\exists} \def\se{{\germ se}} \def\bb{{\bf b}}
  \def\equivalence{\equiv} \let\lt< \let\gt>
  \def\implies{\Rightarrow}\def\mathfrak{\bf}\def\germ{\frak} \def\scr{\cal}
  \ifx\documentclass\undefinedcs
  \def\bf{\fam\bffam\tenbf}\def\rm{\fam0\tenrm}\fi 
  \def\defaultdefine#1#2{\expandafter\ifx\csname#1\endcsname\relax
  \expandafter\def\csname#1\endcsname{#2}\fi} \defaultdefine{Bbb}{\bf}
  \defaultdefine{frak}{\bf} \defaultdefine{=}{\B} 
  \defaultdefine{mathfrak}{\frak} \defaultdefine{mathbb}{\bf}
  \defaultdefine{mathcal}{\cal}
  \defaultdefine{beth}{BETH}\defaultdefine{cal}{\bf} \def\bbfI{{\Bbb I}}
  \def\mbox{\hbox} \def\text{\hbox} \def\om{\omega} \def\Cal#1{{\bf #1}}
  \def\pcf{pcf} \defaultdefine{cf}{cf} \defaultdefine{reals}{{\Bbb R}}
  \defaultdefine{real}{{\Bbb R}} \def\restriction{{|}} \def\club{CLUB}
  \def\w{\omega} \def\exist{\exists} \def\se{{\germ se}} \def\bb{{\bf b}}
  \def\equivalence{\equiv} \let\lt< \let\gt>
\begin{thebibliography}{RuSh 117}
\makeatletter \renewcommand{\@biblabel}[1]{[#1]} \makeatother
\def\eprintfootnotetext{References of the form {\tt math.XX/$\cdots$}
 refer to {\tt arXiv.org} }
\ifx\documentstyle\undefinedcontrolsequence
   \def\anyfootnote{\footnote{*}}
   \else\def\anyfootnote{\footnote}\fi
\def\eprintfn{\ifEprint\anyfootnote{\eprintfootnotetext}\fi\Eprintfalse }
\newif\ifEprint  \Eprinttrue

\bibitem[DoJe81]{DoJe81}A.~Dodd and Ronald~B. Jensen.
\newblock {The core model}.
\newblock {\em Annals of Mathematical Logic}, {\bf 20}:43--75, 1981.

\bibitem[RuSh 117]{RuSh:117}Matatyahu Rubin and Saharon Shelah.
\newblock {Combinatorial problems on trees: partitions, $\Delta$-systems and
  large free subtrees}.
\newblock {\em {Annals of Pure and Applied Logic}}, {\bf 33}:43--81, 1987.

\bibitem[Sh:E12]{Sh:E12}Saharon Shelah.
\newblock {Analytical Guide and Corrections to \cite{Sh:g}.}
\newblock math.LO/9906022.

\bibitem[Sh:e]{Sh:e}Saharon Shelah.
\newblock {\em {Non--structure theory}}, accepted.
\newblock {Oxford University Press}.

\bibitem[Sh 52]{Sh:52}Saharon Shelah.
\newblock {A compactness theorem for singular cardinals, free algebras,
  Whitehead problem and transversals}.
\newblock {\em {Israel Journal of Mathematics}}, {\bf 21}:319--349, 1975.

\bibitem[Sh 108]{Sh:108}Saharon Shelah.
\newblock {On successors of singular cardinals}.
\newblock In {\em {Logic Colloquium '78 (Mons, 1978)}}, volume~97 of {\em
  {Stud. Logic Foundations Math}}, pages 357--380. {North-Holland,
  Amsterdam-New York}, 1979.

\bibitem[Sh 186]{Sh:186}Saharon Shelah.
\newblock {Diamonds, uniformization}.
\newblock {\em {The Journal of Symbolic Logic}}, {\bf 49}:1022--1033, 1984.

\bibitem[Sh 212]{Sh:212}Saharon Shelah.
\newblock {The existence of coding sets}.
\newblock In {\em {Around classification theory of models}}, volume 1182 of
  {\em {Lecture Notes in Mathematics}}, pages 188--202. {Springer, Berlin},
  1986.

\bibitem[Sh 88a]{Sh:88a}Saharon Shelah.
\newblock {Appendix: on stationary sets (in ``Classification of nonelementary
  classes. II. Abstract elementary classes'')}.
\newblock In {\em {Classification theory (Chicago, IL, 1985)}}, volume 1292 of
  {\em {Lecture Notes in Mathematics}}, pages 483--495. {Springer, Berlin},
  1987.
\newblock {Proceedings of the USA--Israel Conference on Classification Theory,
  Chicago, December 1985; ed. Baldwin, J.T.}

\bibitem[Sh 351]{Sh:351}Saharon Shelah.
\newblock {Reflecting stationary sets and successors of singular cardinals}.
\newblock {\em {Archive for Mathematical Logic}}, {\bf 31}:25--53, 1991.

\bibitem[Sh 410]{Sh:410}Saharon Shelah.
\newblock {More on Cardinal Arithmetic}.
\newblock {\em {Archive for Mathematical Logic}}, {\bf 32}:399--428, 1993.
\newblock math.LO/0406550.

\bibitem[Sh 371]{Sh:371}Saharon Shelah.
\newblock {Advanced: cofinalities of small reduced products}.
\newblock In {\em {Cardinal Arithmetic}}, volume~29 of {\em {Oxford Logic
  Guides}}, chapter {VIII}. {Oxford University Press}, 1994.

\bibitem[Sh 355]{Sh:355}Saharon Shelah.
\newblock {$\aleph _{\omega +1}$ has a Jonsson Algebra}.
\newblock In {\em {Cardinal Arithmetic}}, volume~29 of {\em {Oxford Logic
  Guides}}, chapter~II. {Oxford University Press}, 1994.

\bibitem[Sh 386]{Sh:386}Saharon Shelah.
\newblock {Bounding $pp(\mu )$ when $cf(\mu ) > \mu > \aleph _0$ using ranks
  and normal ideals}.
\newblock In {\em {Cardinal Arithmetic}}, volume~29 of {\em {Oxford Logic
  Guides}}, chapter~VI. {Oxford University Press}, 1994.

\bibitem[Sh:g]{Sh:g}Saharon Shelah.
\newblock {\em {Cardinal Arithmetic}}, volume~29 of {\em {Oxford Logic
  Guides}}.
\newblock {Oxford University Press}, 1994.

\bibitem[Sh 400]{Sh:400}Saharon Shelah.
\newblock {Cardinal Arithmetic}.
\newblock In {\em {Cardinal Arithmetic}}, volume~29 of {\em {Oxford Logic
  Guides}}, chapter~IX. {Oxford University Press}, 1994.
\newblock Note: See also [Sh400a] below.

\bibitem[Sh 365]{Sh:365}Saharon Shelah.
\newblock {There are Jonsson algebras in many inaccessible cardinals}.
\newblock In {\em {Cardinal Arithmetic}}, volume~29 of {\em {Oxford Logic
  Guides}}, chapter III. {Oxford University Press}, 1994.

\bibitem[Sh 430]{Sh:430}Saharon Shelah.
\newblock {Further cardinal arithmetic}.
\newblock {\em {Israel Journal of Mathematics}}, {\bf 95}:61--114, 1996.
\newblock math.LO/9610226.

\end{thebibliography}

\enddocument